\documentclass{article}
\usepackage{latexsym}
\usepackage{graphics}
\usepackage{amsmath}
\usepackage{amsfonts}
\begin{document}
\bibliographystyle{plain}
\title{The decomposition of global conformal invariants I:
On a conjecture of Deser and Schwimmer.}
\author{Spyros Alexakis\thanks{University of Toronto, alexakis@math.toronto.edu.
\newline
This work has absorbed
the best part of the author's energy over many years. 
This research was partially conducted during 
the period the author served as a Clay Research Fellow, 
an MSRI postdoctoral fellow,
a Clay Liftoff fellow and a Procter Fellow.  
\newline
The author is immensely indebted to Charles
Fefferman for devoting twelve long months to the meticulous
proof-reading of the present paper. He also wishes to express his
gratitude to the Mathematics Department of Princeton University
for its support during his work on this project.}} \date{}
\maketitle
\newtheorem{proposition}{Proposition}
\newtheorem{theorem}{Theorem}
\newcommand{\Sum}{\sum}
\newtheorem{lemma}{Lemma}
\newtheorem{observation}{Observation}
\newtheorem{formulation}{Formulation}
\newtheorem{definition}{Definition}
\newtheorem{conjecture}{Conjecture}
\newtheorem{corollary}{Corollary}
\numberwithin{equation}{section}
\numberwithin{lemma}{section}
\numberwithin{theorem}{section}
\numberwithin{definition}{section}
\numberwithin{proposition}{section}

\begin{abstract}
This is the first in a series of papers where we  prove a conjecture of Deser and Schwimmer
regarding the algebraic structure of ``global conformal
invariants''; these  are defined to
be conformally invariant integrals of geometric scalars.
 The conjecture asserts that the integrand of
any such  integral can be expressed as a linear
combination of a local conformal invariant, a divergence and of
the Chern-Gauss-Bonnet integrand. 

In this paper we set up an iterative procedure that proves 
the decomposition. We then derive the iterative step in the 
first of two cases, subject to a purely algebraic result 
which is proven in \cite{alexakis4,alexakis5,alexakis6}.
\end{abstract}

\tableofcontents
\section{Introduction.}

\par This is the first in a series of papers, 
\cite{alexakis2,alexakis3,alexakis4,alexakis5,alexakis6}, where we provide a rigorous proof 
of a conjecture of Deser and Schwimmer, originally formulated in \cite{ds:gccaad}.
This series is a continuation of the previous work of the author, 
 \cite{a:dgciI,a:dgciII} which established the conjecture in a 
special case and developed some useful tools which we will use below. 

 The purpose of this introduction is to firstly provide a formulation of
 the conjecture, and then to  give a very brief synopsis of 
some of the main ideas in the proof, followed by a more detailed outline 
of the present paper.

\subsection{Formulation of the problem.}

\par We start by recalling the conjecture of Deser and Schwimmer.
Firstly, we recall a classical notion from Riemannian geometry, that
of a ``scalar Riemannian invariant'':

In brief,  given a Riemannian manifold $(M,g)$, scalar Riemannian invariants 
are {\it intrinsic, scalar-valued} functions of the metric $g$. More precisely:

\begin{definition}
\label{scalar.invariants}
Let $L(g)$ be 
a formal polynomial expression in the in the (formal) variables 
$\partial^{(k)}_{r_1\dots r_k}g_{ij}, k\ge 0$ and $(det g)^{-1}$ (here the 
indices ${}_{r_1},\dots,{}_{r_k},{}_i,{}_j$ take values ${}_1,\dots ,{}_n$).
 Given any coordinate neighborhood $U\subset \mathbb{R}^n$ and 
any Riemannian metric $g$ expressed in the form $g_{ij}dx^idx^j$ in 
terms of the coordinates $\{x^1,\dots,x^n\}\in U$, let $L_g^U$ 
 stand for the function that arises by plugging in the values 
$\partial^{(k)}_{r_1\dots r_k}g_{ij}$,
$(det g)^{-1}$ into the formal expression $L(g)$. We say that $L(g)$ is
a Riemannian invariant of weight $K$ if:

\begin{enumerate}

\item Let $g,g'$ be two Riemannian metrics defined over neighborhoods
 $U,U'\subset \mathbb{R}^n$, and let $L^U_g, L^{U'}_{g'}$ be the scalar-valued
functions defined over $U,U'$ that we obtain by substituting $g,g'$ into
the formal expression $L(g)$. Then we require that if $g,g'$ are isometric via
the map $\Phi:U\rightarrow U'$ then $L^U_g (x)=L^{U'}_{g'}(\Phi(x))$ for every $x\in U$.
(This property is called the intrinsicness property of $L(g)$).

\item Let  $g$ be a Riemannian metric defined over
$U\subset \mathbb{R}^n$ and let $t> 0$. Let $g'$ be the 
Riemannian metric $t^2\cdot g$. Let $L^U_g, L^U_{g'}$ be the scalar-valued
functions defined over $U$ that we obtain by substituting $g,g'$ into
the formal expression $L(g)$. We then require that 
$L^U_{g'}(x)=t^{K}L^U_g(x)$ for every $x\in U$. (We then say that $L(g)$ has weight $K$).
\end{enumerate}
\end{definition}

In view of the first property, a Riemannian invariant $L(g)$ assigns a
well-defined,\footnote{(Meaning coordinate-independent).} 
scalar-valued function to any Riemannian manifold $(M,g)$.

We next review  a classical theorem which essentially goes back to Weyl, \cite{w:cg},
which states that any scalar Riemannian invariant can be expressed in terms of
 complete contractions of covariant derivatives of the curvature tensor.
To state this result precisely, let us recall some basic 
facts from Riemannian geometry:

Given a Riemannian metric $g$ defined over a manifold $M$, consider the 
curvature tensor $R_{ijkl}$ and its covariant derivatives 
$\nabla^{(m)}_{r_1\dots r_m}R_{ijkl}$ (these are thought of as $(0,m+4)$-tensors). 
This gives us a list of {\it tensors} 
defined over $M$.

A natural way to form intrinsic {\it scalars} out of this list of intrinsic
{\it tensors} is by taking tensor products and then contracting indices using
 the metric $g^{ab}$: Firstly we take a (finite) number of tensor products, say:

\begin{equation}
 \label{pity}
\nabla^{(m_1)}_{r_1\dots r_{m_1}}R_{i_1j_1k_1l_1}\otimes\dots\otimes
\nabla^{(m_s)}_{r_1\dots r_{m_s}}R_{i_sj_sk_sl_s},
\end{equation}
thus obtaining a tensor of rank $(m_1+4)+\dots+(m_s+4)$. Then,
 we can repeatedly pick out pairs of indices in the above
expression and contract them against each other using the metric $g^{ab}$.
In the end we obtain a scalar. We will denote such complete contractions
by $C(g)$.\footnote{A rigorous, if somewhat abstract, definition of a
complete contraction appears in the introduction of \cite{a:dgciI}.}
Observe that any such complete contraction will be a scalar Riemannian
invariant of weight $-[(m_1+2)\dots+(m_s+2)]$.
Thus, taking linear combinations of complete contractions of a given weight $w$
we can construct local Riemannian invariants of weight $w$.
 We will denote such linear combinations by $\sum_{r\in R} a_r C^r(g)$ (here $R$
is the index set of the complete contractions, $C^r(g), r\in R$ are
the different complete
contractions appearing and $a_r$ are their coefficients).

Now, a classical result in Riemannian geometry (essentially due to Weyl, \cite{w:cg})
is
that the converse is also true:
For any Riemannian invariant $L(g)$ there exists a (non-unique)
linear combination of complete contractions in the
form (\ref{pity}), $\sum_{r\in R} a_r C^r(g)$ so that
for every manifold $(M,g)$ the value of $L(g)$ is equal to the value of the
linear combination
$\sum_{r\in R} a_r C^r(g)$. Thus from now on we will be identifying
Riemannian invariants with linear combinations of the form:

\begin{equation}
\label{zaxaropoulos}
L(g)=\sum_{l\in L} a_l C^l(g),
\end{equation}
where each $C^l(g)$ is a complete contraction (with respect to the
metric $g$) in the form:

\begin{equation}
\label{zaxaropoulos2}
C^l(g)=contr(\nabla^{(m_1)}R\otimes\dots\otimes\nabla^{(m_a)}R).
\end{equation}
(We do not write out the indices of the tensors involved for brevity).
 We also remark that a complete contraction is determined by
the {\it pattern} according to which different indices contract
against each other. Thus, for example, the complete contraction
$R_{abcd}\otimes R^{abcd}$ is different from the complete 
contraction ${R^a}_{bad}\otimes {R_s}^{bsd}$. The notation
(\ref{zaxaropoulos2}) of course does not encode this pattern of which index
is contracting against which etc.
\newline

{\it The Deser-Schwimmer conjecture:} The conjecture deals with conformally invariant {\it integrals} of 
Riemannian scalars: 

\begin{definition}
\label{globconfinv} Consider a Riemannian invariant $P(g)$ of weight $-n$
($n$ even). We will say that the integral $\int_{M^n}P(g)dV_g$ is 
a ``global conformal invariant'' if 
 the value of $\int_{M^n} P(g)dV_g$ remains invariant
under conformal re-scalings of the metric $g$.

In other words, $\int_{M^n}P(g)dV_g$ is 
a ``global conformal invariant'' if  for any
 $\phi\in C^\infty (M^n)$ we have 
 $\int_{M^n} P(e^{2\phi}g)dV_{e^{2\phi}g}=\int_{M^n}P(g)dV_g$.
\end{definition}

\par In order to state the Deser-Schwimmer conjecture, we recall that a 
{\it local conformal invariant} of weight $-n$
is a Riemannian invariant $W(g)$ for which
$W(e^{2\phi}g)=e^{-n\phi}W(g)$ for every Riemannian metric $g$ and
every function $\phi\in C^\infty(M^n)$. Furthermore, a Riemannian
vector field $T^i(g)$ is a linear combination  $T^i(g)=\sum_{q\in
Q} a_q C^{q,i}(g)$, where each $C^{q,i}(g)$ is a partial
contraction (with one free index) in the form:
$$C^{q,i}(g)=pcontr(\nabla^{(m_1)}R\otimes\dots\otimes\nabla^{(m_a)}R)$$
with $\sum_{t=1}^a (m_t+2)=n-1$. (Notice that for each such vector
field, the divergence $div_iT^i(g)$ is a Riemannian invariant of
weight $-n$). Finally, we recall that
$\operatorname{Pfaff}(R_{ijkl})$ stands for the Pfaffian of the
curvature tensor.\footnote{Recall the Chern-Gauss-Bonnet theorem
which says that for any compact orientable Riemannian $n$-manifold
$(M^n,g)$ we must have $\int_{M^n}
\operatorname{Pfaff}(R_{ijkl})dV_{g^n}=
\frac{2^n\pi^{\frac{n}{2}}(\frac{n}{2}-1)!}{2(n-1)!}\chi(M^n)$.}
\newline

\par The Deser-Schwimmer conjecture \cite{ds:gccaad} asserts:

\begin{conjecture}
\label{desschwim} Let $P(g)$ be a Riemannian invariant of weight $-n$ 
such that the integral $\int_{M^n}P(g)dV_g$ 
is a global conformal invariant. Then
there exists a local conformal invariant $W(g)$, a Riemannian
vector field $T^i(g)$ and a constant $(Const)$ so that $P(g)$ can
be expressed in the form:
\begin{equation} \label{dseqn} P(g)=W(g)+div_i T^i(g)+(Const)\cdot
\operatorname{Pfaff}(R_{ijkl}).
\end{equation}
\end{conjecture}

\par We recall the theorem we proved in \cite{a:dgciI} and
\cite{a:dgciII}:

\begin{theorem}
\label{Iproved}[A] Let $\int_{M^n} P(g)dV_g$ be a global conformal
invariant, where $P(g)$ is in the special form:
\begin{equation}
 \label{protrelos}
P(g)=\sum_{l\in L} a_l contr^l(R_{i_1j_1k_1l_1}\otimes\dots\otimes R_{i_{\frac{n}{2}}j_{\frac{n}{2}}
k_{\frac{n}{2}}l_{\frac{n}{2}}});
\end{equation}
 (i.e.~each of the
complete contractions above has $\frac{n}{2}$ {\it undifferentiated}
factors $R_{ijkl}$). Then $P(g)$ can be expressed in the form:
$$P(g)=W(g)+(const)\cdot \operatorname{Pfaff}(R_{ijkl}).$$
\end{theorem}

\par In this series of papers we will build on the work in \cite{a:dgciI},
\cite{a:dgciII} to prove the whole Deser-Schwimmer conjecture:

\begin{theorem}
\label{thetheorem} Conjecture \ref{desschwim} is true.
\end{theorem}

{\it Related questions:} The motivation for the above theorem, along with its
implications to the notions of $Q$-curvature and renormalized
volume have been discussed in the introduction of \cite{a:dgciI}.
We refer to that paper for that discussion. We just 
 wish to mention the recent work of A.~Juhl
\cite{juhl}, where he obtains new remarkable insight on the 
significance of $Q$-curvature, from an entirely fresh 
point of view. For now, we 
remark that an analogous problem arises in
the context of understanding the asymptotic expansion of the
Szeg\"o kernel of strictly pseudo-convex domains in $\mathbb{C}^n$
(or alternatively of abstract CR-manifolds). In particular, the
leading term of the logarithmic singularity of the Szeg\"o kernel
exhibits a global invariance which is very similar to the one we
discuss here, see \cite{h:lsskgissd}. A further problem related to 
the Deser-Schwimmer conjecture arises 
in K\"ahler geometry: The problem is to understand 
the algebraic structure of the coefficients in 
the Tian-Yau-Zelditch expansion; this is a {\it local version} 
of the classical Riemann-Roch theorem regarding the dimension 
of the space of holomorphic sections of high powers of
 ample line bundles over complex manifolds, see \cite{zelditch} 
for a detailed discussion. The analogy with the Deser-Schwimmer conjecture lies in the fact 
that these coefficients are local invariants of a K\"ahler metric, 
whose {\it integral} over the base manifold remains invariant 
under K\"ahler deformations of the metric. 

Finally, we wish to point
out that alternative  notions of ``global conformal invariants''
have been introduced and studied in the context of general
relativity, see \cite{bs:gciids}.
\newline

Before proceeding to outline the proof of Theorem \ref{thetheorem}
and to synopsize the present paper, we  
briefly digress in order to discuss the relationship
of this work \cite{alexakis1}--\cite{alexakis6} with the study of local invariants 
of geometric structures (mostly Riemannian and conformal) 
and with certain questions motivated by index theory.

{\bf Broad Discussion:} The theory of {\it local} invariants of Riemannian structures 
(and indeed, of more general geometries,
e.g.~conformal, projective, or CR)  has a long history. As stated above, the original foundations of this 
field were laid in the work of Hermann Weyl and \'Elie Cartan, see \cite{w:cg, cartan}. 
The task of writing out local invariants of a given geometry is intimately connected
with understanding which polynomials in a space of tensors with  given symmetries 
 remain invariant under the action of a Lie group. 
In particular, the problem of writing down all 
 local Riemannian invariants\footnote{The {\it scalar-valued} 
invariants considered in Definition \ref{scalar.invariants} are particular cases of such 
local invariants.} reduces to understanding 
the invariants of the orthogonal group. 

 In more recent times, a major program was laid out by C.~Fefferman in \cite{f:ma}
aimed at finding all scalar local invariants in CR geometry. This was motivated 
by the problem of understanding the  
local invariants which appear in the asymptotic expansions of the 
Bergman and Szeg\"o kernels of strictly pseudo-convex CR manifolds,
 in  similar way to which Riemannian invariants appear in the asymptotic expansion  
of the heat kernel; the study of the local invariants
in the singularities of these kernels led to important breakthroughs 
in \cite{beg:itccg} and more recently by Hirachi in \cite{hirachi1}.
 This program was later extended  to conformal geometry in \cite{fg:ci}. 
Both these geometries belong to a 
broader class of structures, the
{\it parabolic geometries}; these are structures which admit a principal bundle with 
structure group a parabolic subgroup $P$ of a semi-simple 
Lie group $G$, and a Cartan connection on that principle bundle 
(see the introduction in \cite{cg1}). 
An important question in the study of these structures 
is the problem of constructing all their local invariants, which 
can be thought of as the {\it natural, intrinsic} scalars of these structures.

  In the context of conformal geometry, the first (modern) landmark 
in understanding {\it local conformal invariants} was the work of Fefferman 
and Graham in 1985 \cite{fg:ci},
where they introduced the {\it ambient metric}. This allows one to 
construct local conformal invariants of any order in odd 
dimensions, and up to order $\frac{n}{2}$ in even dimensions. 
The question is then whether {\it all} invariants arise via this construction. 

The subsequent work of Bailey-Eastwood-Graham \cite{beg:itccg} proved that 
this is indeed true in odd dimensions; in even dimensions, 
they proved that the result holds  
when the weight (in absolute value) is bounded 
by the dimension. The ambient metric construction 
in even dimensions was recently extended by Graham-Hirachi, \cite{grhir}; this enables them to 
indentify in a satisfactory manner {\it all} local conformal invariants, 
even when the weight (in absolute value) exceeds the dimension.  

 An alternative 
construction of local conformal invariants can be obtained via the {\it tractor calculus} 
introduced by Bailey-Eastwood-Gover in \cite{bego}. This construction bears a strong 
resemblance to the Cartan conformal connection, and to 
the work of T.Y.~Thomas in 1934, \cite{thomas}. The tractor 
calculus has proven to be very universal; 
tractor buncles have been constructed \cite{cg1} for an entire class of parabolic geometries. 
The relation betweeen the conformal tractor calculus and the Fefferman-Graham 
ambient metric  has been elucidated in \cite{cg2}.

The present work, while pertaining to the question above
(given that it ultimately deals with the algebraic form of local 
{\it Riemannian} and {\it conformal} invariants), nonetheless addresses a different 
{\it type} of problem:  We here consider Riemannian invariants $P(g)$ for 
which the {\it integral} $\int_{M^n}P(g)dV_g$ remains invariant 
under conformal changes of the underlying metric; we then seek to understand 
the possible algebraic form of the {\it integrand} $P(g)$, 
ultimately proving that it can be de-composed 
in the way that Deser and Schwimmer asserted. 
It is thus not surprising that the prior work on 
 the construction and understanding of local {\it conformal} 
invariants plays a central role in this endeavor, in \cite{alexakis2,alexakis3}. 
We will explain in \cite{alexakis2} how some of the 
local conformal invariants that we {\it identify} in $P(g)$ 
would  be {\it expected} (given the properties of the ambient metric 
but also the insight obtained in \cite{beg:itccg}), while
 others are much less obvious.

On the other hand, our resolution of the Deser-Scwimmer conjecture will also rely 
heavily on a deeper understanding of the algebraic properties 
of the {\it classical} local Riemannian invariants.
 The fundamental theorem of invariant theory (see Theorem B.4 in 
\cite{beg:itccg} and also Theorem 2 in \cite{a:dgciI})
is used extensively throughout this series of papers. However, the most important 
algebraic tool on which our method relies are certain ``main algebraic Propositions'' 
presented in the present paper and \cite{alexakis2}.\footnote{A 
summary of these is provided in
subsection \ref{outline2} below.} 
These are purely algebraic propositions that
 deal with {\it local Riemannian invariants}. 
While the author was led to led to these Propositions 
out of the strategy that he felt was necessary to 
solve the Deser-Schwimmer conjecture, they can 
be thought of as results of independent interest. 
The {\it proof} of these Propositions, presented
 in \cite{alexakis4,alexakis5,alexakis6} is in fact 
not particularily intuitive. It is the author's 
sincere hope that deeper insight 
will be obtained in the future as to  {\it why} these algebraic 
 Propositions hold.

{\it Index Theory:} Questions similar to the Deser-Schwimmer
conjecture  arise naturaly 
in index theory; a good reference for such questions is   \cite{bgv}.
 For example, in the heat kernel proof of the index theorem (for Dirac operators) 
by Atiyah-Bott-Patodi \cite{abp}, 
the authors were led to consider integrals arising in the (integrated) expansion of 
the heat kernel over Riemannian manifolds of general Dirac operators, 
and sought to understand the 
local structure  of the integrand.\footnote{We note that
the geometric setting in \cite{abp} is more 
general than the one in the Deser-Scwimmer conjecture:
 In particular one considers vector bundles, equipped with an auxiliary connection,
  over compact Riemannian manifolds; the 
local invariants thus depend {\it both} on the 
curvature of the Riemannian metric {\it and} the curvature of the connection.} 
In that setting, however, the fact that one deals with 
a {\it specific}  integrand which arises in the heat kernel expansion 
plays a key role in the understanding of its local  structure. 
This is true both of the original proof of Patodi, Atiyah-Bott-Patodi \cite{abp} 
and of their subsequent simplifications and generalizations by Getzler, 
Berline-Getzler-Vergne, see \cite{bgv}. 

The closest analogous problem to the one considered here 
is the work of Gilkey and Branson-Gilkey-Pohjanpelto, \cite{gilkey,bgp}.  
In \cite{gilkey}, Gilkey considered Riemannian invariants $P(g)$ for 
which the {\it integral} $\int_{M^n}P(g)dV_g$ on any given 
(topological) manifold $M^n$ has a given
value, {\it independent of the metric $g$}. 
He proved that $P(g)$ must then be equal to a divergence, 
plus possibly a multiple of the Chern-Gauss-Bonnet integrand 
if the weight of $P(g)$ agrees with the dimension in absolute value. 
In \cite{bgp} the authors considered the problem of Deser-Schwimmer 
for localy conformally flat metics and derived the same decomposition 
(for {\it locally conformaly flat metrics}) as in \cite{gilkey}.
Although these two results can be considered  precursors of ours, the methods
there are entirely different from the ones here; it is highly unclear whether  
the methods of \cite{gilkey,bgp} could be applied to the problem at hand.

\subsection{Outline of the argument.}
\label{outline}

{\bf A one-page outline of the argument:}
 The Deser-Schwimmer conjecture is proven
by a multiple induction. At the roughest level, the induction
works as follows: Express $P(g)$ as a linear combination of complete
contractions:

\begin{equation}
\label{bo'az1}P(g)=\sum_{l\in L} a_l C^l(g),
\end{equation}
each $C^l(g)$ in the form (\ref{zaxaropoulos2}).  

 The different complete contractions $C^l(g)$ appearing above can be grouped up
 into ``categories'' according to certain algebraic features of the tensors
involved.\footnote{See the next page for more details.} Accordingly, we divide 
the index set $L$ into subsets $L^1,\dots,L^T$ so that 
the terms indexed in the same index set $L^t$ belong to the same
category (and vice versa), and $\bigcup_{t=1}^T L^t=L$; accordingly, we write:
\begin{equation}
 \label{diaolemenos}
P(g)=\sum_{t=1}^T\sum_{l\in L^t} a_l C^l(g).
\end{equation}
We will also introduce a {\it grading} among the set of categories:
 A given category of complete contractions
 will be ``better'' or ``worse'' than any other given category.
 For future reference, the ``best''
 category of complete contractions are the ones
 with $\frac{n}{2}$ factors.

Assume that in (\ref{diaolemenos}), for each pair $1\le \alpha<\beta\le T$
the category of complete contractions indexed in $L^\beta$
is  ``worse'' than the category of complete
contractions  indexed in $L^\alpha$. (Therefore,
in particular the ``worst'' category of complete contractions
in (\ref{diaolemenos}) is the category $\sum_{l\in L^T} a_l C^l(g)$).

The main step of our induction is to prove that {\it unless the
complete contractions $C^l(g), l\in L^T$ are in the ``best''
category},\footnote{I.e.~unless $P(g)$ is in the form (\ref{protrelos}).}
there exists a local conformal invariant $W(g)$
and a divergence of a vector field $div_iT^i(g)$
so that:

\begin{equation}
\label{bo'az2}\sum_{l\in L^T} a_l
C^l(g)-W(g)-div_iT^i(g)=\sum_{l\in L^{new}} a_l C^l(g),
\end{equation} where the complete contractions in the RHS of the
above belong to categories that are all ``better'' than the
category of $\sum_{l\in L^T} a_l C^l(g)$.

\par Observe that once this ``main step'' is proven, we can iteratively apply it to 
derive that there exists a local conformal invariant $\tilde{W}(g)$ and
a divergence $div_i\tilde{T}^i(g)$ so that:

\begin{equation}
\label{bo'az3}
 P(g)-\tilde{W}(g)-div_i\tilde{T}^i(g)=\tilde{P}(g),
 \end{equation} where
$\tilde{P}(g)$ is a linear combination of terms with $\frac{n}{2}$
factors. Furthermore, $\int_{M^n}\tilde{P}(g)dV_g$ {\it is also a global
conformal invariant.} Therefore, invoking the main theorem of
\cite{a:dgciII}, we derive that $\tilde{P}(g)$ can be written in
the form:

\begin{equation}
\label{bo'az4}\tilde{P}(g)=W'(g)+(Const)\cdot
\operatorname{Pfaff}(R_{ijkl}),
\end{equation}
where $W'(g)$ is a local conformal
invariant.\footnote{$\operatorname{Pfaff}(R_{ijkl})$ is the
Pfaffian of the curvature tensor (i.e. the Gauss-Bonnet
integrand).} Therefore, combining (\ref{bo'az3}) and
(\ref{bo'az4}) we derive the Deser-Schwimmer conjecture.

\subsection{A more detailed outline of the present paper.}
\label{outline2}

This series of papers (the present paper and 
\cite{alexakis2,alexakis3,alexakis4,alexakis5,alexakis6}) 
can naturally be divided into two parts: Part I (which
consists of the present paper together with \cite{alexakis2} and \cite{alexakis3})
 proves the Deser-Schwimmer
conjecture {\it subject to proving certain ``main algebraic propositions''}; 
these are Proposition \ref{pregiade} in the  present paper, 
and the two propositions 3.1, 3.2 in the section ``The important tools'' in \cite{alexakis2}. 
 Part II (which consists of the papers \cite{alexakis4},
\cite{alexakis5}, \cite{alexakis6})
 are devoted to proving these ``main algebraic propositions''.
Thus, the logical dependence of this work is that the present paper and \cite{alexakis2,alexakis3}
depend on \cite{alexakis4,alexakis5,alexakis6}. 
 Here we present a more detailed, yet broad, outline of this entire work, putting 
emphasis on the results proven in the present paper. In the subsequent papers
of this series, we will  provide further synopses of the 
other main ideas that appear in this work. 
\newline

{\bf ``Categories'' and the notion of ``better'' vs.~``worse'' categories:}
We now explain in more detail the notion of
``categories'' explained above, and how one category is ``better''
or ``worse'' than another category. Firstly,
recall that the curvature tensor $R_{ijkl}$ admits a natural decomposition 
 into its trace-free part
(the Weyl tensor) and its trace parts (consisting essentially of
the Ricci tensor),\footnote{See (\ref{Weyl}) below.} we will write out the global conformal
invariant $P(g)$ as a linear combination of complete contractions
involving covariant derivatives of the Weyl tensor $\nabla^{(m)}W$ and covariant
derivatives of the Schouten tensor  $\nabla^{(a)}P$:\footnote{The Schouten tensor,
defined in (\ref{WSch}), is a trace-adjustment of the Ricci tensor. For the
purpose of this brief introduction, the reader may wish to think
of the Schouten tensor as ``essentially'' the Ricci tensor.}

\begin{equation}
\label{bo'az5}P(g)=\sum_{l\in L} a_l
contr^l(\nabla^{(m_1)}W\otimes\dots\otimes\nabla^{(m_a)}W
\otimes\nabla^{(p_1)}P\otimes\dots\otimes\nabla^{(p_q)}P).
\end{equation}

\par Then, two complete contractions in the above form belong to
the same ``category'' if they have the same number $a+b$ of
factors (in total), and also if they have the same number $b$ of
factors $\nabla^{(p)}P$. Furthermore, if we consider two complete
contractions $C^1(g)$ and $C^2(g)$ in the
above form, then $C^1(g)$ is ``worse'' than $C^2(g)$ if it has
{\it fewer} factors in total. If $C^1(g)$, $C^2(g)$ have the same
 number of factors {\it in total},  then $C^1(g)$ is ``worse'' than
$C^2(g)$ if it has more factors $\nabla^{(a)}P$.

\par Thus: Let
$\sigma$ be the minimum total number of factors among all the
complete contractions indexed in $L$ in (\ref{bo'az1}).
Among the complete
contractions with $\sigma$ factors, let $s$ be the maximum number of
factors $\nabla^{(p)}P$. Then the ``worst'' complete contractions
in (\ref{bo'az5}) are the ones with $\sigma-s$ factors $\nabla^{(m)}W$ and $s$
factors $\nabla^{(p)}P$. Denote the index set of the ``worst'' complete
contractions by $L^s_\sigma\subset L$. We define
$P(g)_{worst-piece}:=\sum_{l\in L_\sigma^s} a_l C^l(g)$.

\par Our main claim is that if $\sigma<\frac{n}{2}$
then $P(g)_{worst-piece}$ can be expressed as follows:

\begin{equation}
\label{thegoal} P(g)_{worst-piece}=W(g)+div_iT^i(g)+\sum_{f\in F^1} a_g
C^f_g(\phi)+\sum_{f\in F^2} a_g C^f_g(\phi),
\end{equation}
where $W(g)$ is a local conformal invariant,\footnote{In fact, $W(g)=0$ {\it unless} $s=0$.} 
$div_iT^i(g)$ is the divergence of a Riemannian vector 
field and each of the complete contractions indexed in $F^1, F^2$
 are in the form:

\begin{equation}
\label{protheform}
contr(\nabla^{(m_1)}W\otimes\dots\otimes\nabla^{(m_a)}W\otimes
\nabla^{(p_1)}P\otimes\dots\otimes\nabla^{(p_b)}P),
\end{equation}
with the following additionnal properties: 
The terms indexed in $F^1$ have more than $\sigma$
factors in total, while the terms indexed in $F^2$ have $\sigma$
factors in total but strictly fewer than $s$ factors
$\nabla^{(p)}P$. (In other words, the terms indexed in $F^1,F^2$ are
``better'' than the terms in $P(g)_{worst-piece}$).
\newline

{\bf The main ideas in the derivation of (\ref{thegoal}), and a discussion of the
difficulties:} The starting point in deriving (\ref{thegoal}) is 
to pass from the {\it invariance under integration} enjoyed by $P(g)$ to 
a {\it local formula} for its conformal variation. 

The main tool we developed in \cite{a:dgciI} (in
order to address the Deser-Schwimmer conjecture) is the so-called
super divergence formula. In one sentence, this formula applies to
the {\it conformal variation $I_g(\phi)$} of
$P(g)$,\footnote{We recall that
$I_g(\phi):=e^{n\phi}P(e^{2\phi}g)-P(g)$; thus $I_g(\phi)$ is a
differential operator, depending on an auxiliary function
$\phi$.} and explicitly expresses $I_g(\phi)$ as a divergence of
a vector-valued differential operator $X^i_g(\phi)$:

\begin{equation}
\label{bo'az6}I_g(\phi)=div_iX^i_g(\phi).
\end{equation}

\par Then, the main task in proving (\ref{thegoal}) is two-fold: Firstly, to
identify a ``piece'' in $I_g(\phi)$ which is in one-to-one 
correspondence with the  the ``worst piece'' of $P(g)$. 
Secondly, to {\it use} the fact that $I_g(\phi)$ can be expressed
as a divergence\footnote{Via the ``super divergence formula'' from
\cite{a:dgciI}.} to derive (\ref{thegoal}).

\par We distinguish two main cases in order to
 derive (\ref{thegoal}): Either $s>0$ or $s=0$. 
We prove (\ref{thegoal}) when $s>0$ in the present paper. 
We prove (\ref{thegoal}) when $s=0$ in \cite{alexakis2}.
We outline very roughly the
proof of (\ref{thegoal}) in these two cases, without illustrating
the use of the ``main algebraic Propositions'' in this proof.
Then, only for the case $s>0$ we explain very briefly how the
``main algebraic Propositions'' are used  in the proof.
\newline

{\bf An outline of the proof of (\ref{thegoal}) when $s>0$:} In
this case it not hard to ``recognize'' the worst piece of $P(g)$ in
$I_g(\phi)$. We let $I^s_g(\phi)$ stand for the linear combination of terms in
$I_g(\phi)$ of homogeneity $s$ in the function 
$\phi$.\footnote{In other words, $I^s_g(\phi)=\frac{d^s}{dt^s}|_{t=0}e^{nt\phi}P(e^{2t\phi}g)$.} 
We have proven in \cite{a:dgciI}
 that the super divergence formula can also be applied to $I^s_g(\phi)$.
Thus in this case, we consider $I^s_g(\phi)$. By virtue of the
conformal invariance of the Weyl tensor,\footnote{See
(\ref{transweyl}) below.} and the transformation law of the Schouten
tensor and the Levi-Civita connection\footnote{See
(\ref{wschtrans}), (\ref{levicivita}) respectively below.} we
observe that if we write out
$$P(g)_{worst-piece}=\sum_{l\in L_\sigma^s} a_l contr^l(\nabla^{(m_1)}W\otimes
\dots\otimes\nabla^{(m_{\sigma-s})}W\otimes
\nabla^{(p_1)}P\otimes\dots\otimes\nabla^{(p_s)}P),$$ (recall that
by the definition of $P(g)_{worst-piece}$,
all complete contractions in the RHS of the above will have $s$ factors $\nabla^{(p)}P$
and $\sigma-s$ factors $\nabla^{(m)}W$), then:
\begin{equation}
\label{bo'az6,5} \begin{split} & I^s_g(\phi)=(-1)^s\sum_{l\in
L_\sigma^s} a_l contr^l(\nabla^{(m_1)}W\otimes
\dots\otimes\nabla^{(m_{\sigma-s})}W\otimes
\nabla^{(p_1+2)}\phi\otimes\dots\\&\otimes\nabla^{(p_s+2)}\phi)
+\sum_{j\in Junk-Terms} a_j C^j_g(\phi),
\end{split}
\end{equation}
 where each of the terms $C^j_g(\phi)$, $j\in Junk-Terms$ has at least $\sigma+1$
factors in total. {\it The complete contractions
$contr^l(\nabla^{(m_1)}W\otimes
\dots\otimes\nabla^{(m_{\sigma-s})}W\otimes
\nabla^{(p_1+2)}\phi\otimes\dots\otimes\nabla^{(p_s+2)}\phi)$
arise from the complete contractions
$contr^l(\nabla^{(m_1)}W\otimes
\dots\otimes\nabla^{(m_{\sigma-s})}W\otimes
\nabla^{(p_1)}P\otimes\dots\otimes\nabla^{(p_s)}P)$ by just
replacing each factor $\nabla^{(p)}_{r_1\dots
r_p}P_{r_{p+1}r_{p+2}}$ by a factor $\nabla^{(p+2)}_{r_1\dots
r_{p+2}}\phi$}. Thus, (\ref{bo'az6,5}) provides us with a way to
``recover'' the worst piece $P(g)_{worst-piece}$ from
$I^s_g(\phi)$.

\par Now, we recall that $\int_MI^s_g(\phi)dV_g=0$, thus we can
apply the main result of \cite{a:dgciI}\footnote{I.e. the ``super
divergence formula''.} to the above integral equation and derive a
{\it local equation} which expresses $I^s_g(\phi)$ as a divergence
of a vector-valued differential operator. In fact, the super
divergence formula gives much more precise information: It shows
 that there exists a linear combination
$X^i_g(\phi)$ of partial contractions, $X^i_g(\phi)=\sum_{r\in R}
a_r C^{r,i}_g(\phi)$,\footnote{``Partial contractions'' with one
free index, to be precise.} where each $C^{r,i}_g(\phi)$ is in
the form:
\begin{equation}
\label{bo'az7}pcontr(\nabla^{(m'_1)}R\otimes
\dots\otimes\nabla^{(m_{\sigma-s})}R\otimes
\nabla^{(b_1)}\phi\otimes\dots\otimes\nabla^{(b_s)}\phi)
\end{equation} (each $\nabla^{(m')}R$ is the $m^{th}$ iterated
covariant derivative of the curvature tensor, and each
$\nabla^{(b)}\phi$ is the $b^{th}$ iterated covariant derivative
of the function $\phi$), so that:

\begin{equation}
\label{bo'az6,75} \begin{split} &(-1)^s\sum_{l\in L_\sigma^s} a_l
contr^l(\nabla^{(m_1)}W\otimes
\dots\otimes\nabla^{(m_{\sigma-s})}W\otimes
\nabla^{(p_1+2)}\phi\otimes\dots\otimes\nabla^{(p_s+2)}\phi)=
\\&div_i\sum_{r\in R} a_r C^{r,i}_g(\phi)+\sum_{j\in Junk-Terms}
a_j C^j_g(\phi).
\end{split}
\end{equation}

{\it Furthermore, the super divergence formula also implies that
each $b_i\ge 2$} (apart from certain very special cases where we
may have $b_i=1$ for some of the vector fields $C^{r,i}_g(\phi)$--for
 the purposes of this introduction, we will assume that each
$b_i\ge 2$). We will now show how the main claim, (\ref{thegoal}), can be derived from 
(\ref{bo'az6,75}) when $s=\sigma$. We will then dicsuss why this 
direct approach fails when $0<s<\sigma$. 
\newline

{\it Proof of (\ref{thegoal}) in the case $s=\sigma$:}  Now, in
the case $\sigma=s$, we derive in subsection \ref{sofiavembo}
below that the vector field needed for (\ref{thegoal}) is the
vector field $X^i(g)$  that formally arises from
$X^i_g(\phi)=\sum_{r\in R} a_r C^{r,i}_g(\phi)$ in
(\ref{bo'az6,75}) by replacing each factor
$\nabla^{(b_i)}_{x_1\dots x_{b_i}}\phi$ by a factor
$-\nabla^{(b_i-2)}_{x_1\dots x_{b_i-2}}P_{x_{b_i-2}x_{b_i}}$
(observe that the condition $b_i\ge 2$ is necessary for this
operation to be well-defined).
\newline

\par On the other hand, in the case $s<\sigma$ one {\it cannot}
 derive (\ref{thegoal}) by directly applying the super divergence
formula to the integral equation $\int_{M^n} I^s_g(\phi)dV_g=0$
and then replacing the factors $\nabla^{(b)}\phi$ as above. We
next discuss {\it why} this direct approach will fail in this case:

{\it The difficulty in deriving (\ref{thegoal}) when $0<s<\sigma$:} If one were to directly
apply the super divergence formula to the integral equation
$\int_{M^n} I^s_g(\phi)dV_g=0$, one would derive a local equation
in the form (\ref{bo'az6,75}). Now, if one were to pick out the
terms with $\sigma$ factors in (\ref{bo'az6,75}) and then replace
the factors $\nabla^{(b)}_{r_1\dots r_b}\phi$ ($b\ge 2$) by
factors $\nabla^{(b-2)}_{r_1\dots r_b}P_{r_{b-1}r_b}$ (as in the case $s=\sigma$),
 one would derive an equation:

\begin{equation}
\label{bo'az7,5}
\begin{split} &(-1)^s\sum_{l\in L_\sigma^s} a_l
contr^l(\nabla^{(m_1)}W\otimes
\dots\otimes\nabla^{(m_{\sigma-s})}W\otimes
\nabla^{(p_1+2)}\phi\otimes\dots\otimes\nabla^{(p_s+2)}\phi)=
\\&div_i\sum_{r\in R} a_r C^{r,i}(g)+\sum_{k\in K} a_k C^k(g)+
\sum_{j\in Junk-Terms} a_j C^j_g(\phi),
\end{split}
\end{equation}
where the terms indexed in $Junk-Terms$ will have at least
$\sigma+1$ factors, {\it but the terms indexed in $K$ will be in
the form (\ref{protheform}) with  $\sigma$ factors in total, and
may have as many as $\sigma=1$ factors $\nabla^{(a)}P$}.
 In other words, the complete contractions
indexed in $K$ {\it do not necessarily} have
fewer than $s$ factors $\nabla^{(a)}P$. Therefore in the
 language of the ``one-page summary'',
the terms indexed in $K$ in the RHS of (\ref{bo'az7,5}) are not necessarily
 ``better'' than the terms in the LHS of (\ref{bo'az7,5}).

 \par Therefore in the case $0<s<\sigma$ we will use the super
 divergence formula applied to $I^s_g(\phi)$ {\it in a less straightforward way} to derive a
 stronger claim than (\ref{bo'az6,75}):
\newline

{\it The remedy when $0<s<\sigma$:} We will  prove that there
 exists a linear combination of vector fields, $\sum_{y\in Y} a_y
 C^{y,i}_g(\phi)$, where each $C^{y,i}_g(\phi)$ is in the form:

 $$pcontr(\nabla^{(m_1)}W\otimes\dots\otimes\nabla^{(m_{\sigma-s})}W
 \otimes\nabla^{(b_1)}\phi\otimes\dots\otimes\nabla^{(b_s)}\phi)$$
(with each $b_j\ge 2$) so that:
\begin{equation}
\label{bo'az8} (-1)^s\sum_{l\in L_\sigma^s} a_l
C^l_g(\phi)=div_i\sum_{y\in Y} a_y C^{y,i}_g(\phi)+\sum_{j\in
Junk} a_j C^j_g(\phi),
\end{equation}
where the terms indexed in $Junk$ have at least $\sigma+1$ factors
in total. (A brief discussion explaining the derivation of
(\ref{bo'az8}) is provided further down in this 10-page summary,
in ``A rough discussion of the ``main algerbraic   Proposition'').
\newline

\par Then, (\ref{bo'az8}) implies (\ref{thegoal}):
For each $y\in Y$ formally construct 
a vector field
 $C^{y,i}(g)$ in the form (\ref{protheform}) (with one free index) 
by replacing the factors $\nabla^{(b_j)}_{r_1\dots
r_{b_j}}\phi$ by factors
 $-\nabla^{(b_j-2)}_{r_1\dots r_{b_{j-2}}}P_{r_{b_{j-1}}r_{b_j}}$. We
 then derive (in section \ref{peresprez} below) that the divergence needed for
 (\ref{thegoal}) is precisely $\sum_{y\in Y} a_y C^{y,i}(g)$.

 {\it Note:} Observe that in this case $s>0$, (\ref{thegoal}) holds
{\it without} a
 local conformal invariant $W(g)$ in the RHS.
\newline

{\bf A rough description of the ``main algrebraic Proposition'' \ref{pregiade}
and of its use in proving equation (\ref{thegoal}) when $s>0$.}
\newline

{\bf The ``main algrebraic Proposition'' \ref{pregiade}:}

 First a little notation. We will be considering linear
 combinations of tensor fields, $\sum_{l\in L_\mu} a_l
 C^{l,i_1\dots i_\mu}_g(\Omega_1,\dots,\Omega_p)$, where each
 $C^{l,i_1\dots i_\mu}_g(\Omega_1,\dots,\Omega_p)$
is a partial contraction (with $\mu$ free indices) in the form:
\begin{equation}
\label{paparigas}
pcontr(\nabla^{(m_1)}R\otimes\dots\nabla^{(m_r)}R\otimes\nabla^{(a_1)}
\Omega_1\otimes\dots\otimes\nabla^{(a_p)}\Omega_p),
\end{equation}
with a given number (say $\tau(=r+p)$) of factors in total; among 
these a given number $p$ of factors are in the form 
$\nabla^{(a)}\Omega_x, 1\le x\le p$,\footnote{In other words, they are $a^{th}$ 
covariant derivatives of a scalar function $\Omega_x$.} and the remaining 
$\tau-p$ are in the form $\nabla^{(m)}R$.\footnote{The $m^{th}$ 
covariant derivatives of the curvature tensor $R$.}  Notice also
that there is a  given number $p$ of factors
$\nabla^{(a)}\Omega_x, 1\le x\le p$.
We furthermore require that
each $a_j\ge 2$ for each tensor field above,\footnote{Tensor fields
in the form with this property will be called ``acceptable''.} and that each tensor
 field has no {\it internal contractions}.\footnote{Recall from \cite{a:dgciI} that 
in a tensor field in the form (\ref{paparigas}),
an {\it internal contraction} is a pair of two indices that belong to the same 
factor and are contracting against each other.}

\par We also let $\sum_{l\in L'} a_l
 C^{l,i_1\dots i_{b_l}}_g(\Omega_1,\dots,\Omega_p)$ stand for a 
linear combination of (acceptable) tensor fields
in the form (\ref{paparigas}), each with rank $b_l \ge \mu+1$. 
Recall that for each free index ${}_{i_s}$ 
in $C^{l,i_1\dots i_{b_l}}_g$, the divergence 
$div_{i_s}C^{l,i_1\dots i_{b_l}}_g$ is a sum of $\tau$ 
partial contractions of rank $b_l-1$: the first
 summand arises when we hit the first factor $T_1$ in $C^{l,i_1\dots i_{b_l}}_g$
by a derivative $\nabla^{i_s}$ and contract the upper index ${}^{i_s}$ against 
the free index ${}_{i_s}$; the second summand 
arises when we hit the first factor $T_1$ in $C^{l,i_1\dots i_{b_l}}_g$
by a derivative $\nabla^{i_s}$ and contract the upper index ${}^{i_s}$ against 
the free index ${}_{i_s}$, etc.

For tensor fields in the form (\ref{paparigas}) we will let
$Xdiv_{i_h}C^{l,i_1\dots i_{a_l}}_g(\Omega_1,\dots,\Omega_p)$
stand for the sum of the $\tau-1$ terms in $div_{i_h}C^{l,i_1\dots
i_{b_l}}_g(\Omega_1,\dots,\Omega_p)$
 where the derivative $\nabla^{i_h}$ may hit any factor {\it other than the one} to
which the free index ${}_{i_h}$ belongs.\footnote{This rather
strange definition fits in with the conclusion of the super
divergence formula--see section \ref{mainconseq} below.} The
assumption of the ``main Proposition'' \ref{pregiade} is that:

\begin{equation}
\label{ekloges2007}
\begin{split}
& \sum_{l\in L_\mu} a_lXdiv_{i_1}\dots Xdiv_{i_\mu}
 C^{l,i_1\dots i_\mu}_g(\Omega_1,\dots,\Omega_p)=
\\&\sum_{l\in L'} a_lXdiv_{i_1}\dots Xdiv_{i_{b_l}}
 C^{l,i_1\dots i_{b_l}}_g(\Omega_1,\dots,\Omega_p)+(Junk-Terms),
\end{split}
\end{equation}
where $(Junk-Terms)$ here stands for a generic linear combination
of complete contractions with at least $\tau+1$ factors.\footnote{Whereas
the terms in the LHS of the above each have $\tau$ factors.}

\par The claim of the ``main algebraic Proposition'' is that there exists a linear
combination of acceptable $(\mu+1)$-tensor fields, say
$\sum_{h\in H} a_h C^{h,i_1\dots i_{\mu+1}}_g(\Omega_1,\dots,\Omega_p)$,  with
 each $C^{h,i_1\dots i_{\mu+1}}_g(\Omega_1,\dots,\Omega_p)$
in the form (\ref{paparigas}), so that:

\begin{equation}
\label{paparigasconcl}
\begin{split}
& \sum_{l\in L_\mu} a_l
 C^{l,(i_1\dots i_\mu)}_g(\Omega_1,\dots,\Omega_p)-Xdiv_{i_{\mu+1}}
\sum_{h\in H} a_h C^{h,(i_1\dots i_\mu)i_{\mu+1}}_g(\Omega_1,\dots,\Omega_p)
\\&=(Junk-Terms),
\end{split}
\end{equation}
where the $(Junk-Terms)$ in the above stand for a linear
combination of complete contractions with at least $\tau+1$
factors. Here the symbol ${}^{(i_1\dots i_\mu)}$ means that we are
{\it symmetrizing} over the indices ${}^{i_1},\dots,{}^{i_\mu}$.

\par We next highlight how (\ref{paparigasconcl}) can be used to prove
(\ref{thegoal}) when $s>0$:
\newline

{\bf The use of the ``main algebraic Proposition'' in deriving (\ref{thegoal}) (when $s>0$):}
We present here the argument from section \ref{prfkillthet} in brief:

Equation (\ref{thegoal}) is proven by a new induction. Write out:
\begin{equation}
\label{muvaz} P(g)_{worst-piece}=\sum_{l\in L} a_l
contr^l(\nabla^{(m_1)}W\otimes\dots\otimes\nabla^{(m_{\sigma-s})}W\otimes
\nabla^{(p_1)}P\otimes\dots\otimes\nabla^{(p_s)}P).
\end{equation}
 We assume
that among the complete contractions in
$P(g)_{worst-piece}$\footnote{Recall that the complete
contractions in $P(g)_{worst-piece}$ are all in the form
(\ref{protheform}) with $\sigma$ factors in total, of which $s\ge
1$ are in the form $\nabla^{(a)}P$.} the minimum number of
internal contractions is $\beta\ge 0$. 
We  denote by $L_\beta\subset L$ the index set of complete
contractions with $\beta$ internal contractions. (Thus the
complete contractions indexed in $L\setminus L_\beta$ will each
have at least $\beta+1$ internal contractions).
 We will then show in section \ref{prfkillthet} that there exists a divergence of a Riemannian vector field,
$div_iT^i(g)$, as allowed in the statement of
Conjecture \ref{desschwim} such that:

\begin{equation}
\label{vakifleri}
\begin{split}
& \sum_{l\in L_\beta} a_l
contr^l(\nabla^{(m_1)}W\otimes\dots\otimes\nabla^{(m_{\sigma-s})}W\otimes
\nabla^{(p_1)}P\otimes\dots\otimes\nabla^{(p_s)}P)=
\\&div_i T^i(g)+ \sum_{t\in T} a_t contr^t(g)+(Allowed),
\end{split}
\end{equation}
where the complete contractions indexed in
$T$ are in the form (\ref{protheform}) with $\sigma$ factors in
total, of which $s$ are in the form $\nabla^{(a)}P$, {\it and with
 $\beta+1$ internal contractions in total}. Furthermore
 $(Allowed)$ stands for a generic linear combination of
 complete contractions that are allowed in the right hand
 side of (\ref{thegoal}).\footnote{In other words, the complete
 contractions indexed in $(Allowed)$ are
{\it either} complete contractions with more than $\sigma+1$
factors in total, {\it or} they are complete contractions in the
form (\ref{protheform}) with $\sigma$ factors in total, but
strictly fewer than $s$ factors $\nabla^{(a)}P$.} Observe that if
we can show then will follow by iteratively repeating this step at
most $\frac{n}{2}$ times.\footnote{This is because a complete
contraction in the form (\ref{protheform}) with weight $-n$ can contain at most
$\frac{n}{2}$ internal contractions.}
\newline

{\it Mini-Outline of the proof of (\ref{vakifleri}):} We recall
that
\begin{equation}
\label{sabine1}
contr^l(\nabla^{(m_1)}W\otimes\dots\otimes\nabla^{(m_{\sigma-s})}W\otimes
\nabla^{(p_1+2)}\phi\otimes\dots\otimes\nabla^{(p_s+2)}\phi)
\end{equation}
stands for the complete contraction that arises from
\begin{equation}
\label{sabine2}contr^l(\nabla^{(m_1)}W\otimes\dots\otimes\nabla^{(m_{\sigma-s})}W\otimes
\nabla^{(p_1)}P\otimes\dots\otimes\nabla^{(p_s)}P)
\end{equation}
by replacing each factor $\nabla^{(a)}_{r_1\dots r_a}P_{ij}$ by a
factor $\nabla^{(a+2)}_{r_1\dots r_aij}\phi$. We will denote by
$contr^l(\phi)$  the complete contraction in the form
(\ref{sabine1}); we will also denote by $\overline{contr}^l(\phi)$
the complete contraction:
\begin{equation}
\label{sabine3}
contr^l(\nabla^{(m_1)}R\otimes\dots\otimes\nabla^{(m_{\sigma-s})}R\otimes
\nabla^{(p_1+2)}\phi\otimes\dots\otimes\nabla^{(p_s+2)}\phi)
\end{equation}
which arises from $contr^l(\phi)$ by formally replacing each
factor $\nabla^{(m)}_{r_1\dots r_m}W_{ijkl}$ by a factor
$\nabla^{(m)}_{r_1\dots r_m}R_{ijkl}$ (possibly times a
constant--but for the purposes of this introduction we will ignore
this fact). Observe that the resulting tensor fields still have
$\beta$ internal contractions in total, and also have each
function $\phi$ differentiated at least twice. We will then {\it
prove} in section \ref{reducindstat}
that the integral equation $\int_{M^n}
I^s_g(\phi)dV_g=0$ implies a {\it new} integral equation in the form:
\begin{equation}
\label{sabine4} \int_{M^n} (-1)^\beta\sum_{l\in L_\beta}
a_l\overline{contr}^l(\phi)+\sum_{v\in V} a_v
contr^v_g(\phi)+(Junk-Terms)dV_g=0.
\end{equation}
Here the complete contractions indexed in $V$ are in the form
(\ref{sabine3}) and
 have $\sigma$ factors and at least $\beta+1$ internal contractions in
total,\footnote{(And none of these internal contractions involve
two indices from among the indices ${}_i,{}_j,{}_k,{}_l$ in a
factor $\nabla^{(m)}R_{ijkl}$--this detail is only relevant for
the next sentence).} and each factor $\nabla^{(b_i)}\phi$ with
$b_i\ge 2$. Then, applying the super
divergence formula,\footnote{See section \ref{mainconseq} below.}
we derive a {\it local} equation:

\begin{equation}
\label{bastakardiamou1} \begin{split} & \sum_{u\in U_\beta}
a_uXdiv_{i_1}\dots Xdiv_{i_\beta}\overline{pcontr}^{u,i_1\dots
i_\beta}(\phi)+ \sum_{v\in V} a_v Xdiv_{i_1}\dots
Xdiv_{i_{b_v}}C^{v,i_1\dots i_{b_v}}_g(\phi)
\\&=(Junk-Terms),
\end{split}
\end{equation}
where the tensor fields $\overline{pcontr}^{u,i_1\dots
i_\beta}(\phi)$, $C^{v,i_1\dots i_{b_v}}_g(\phi)$ arise from the
complete contractions $\overline{contr}^u(\phi)$,
$contr^v_g(\phi)$ by formally replacing each internal
contraction\footnote{Which by hypothesis will consist of two
indices in the same factor contracting against each other--i.e.
two indices in the form $(\nabla^a,{}_a)$} by a free
index.\footnote{I.e. in the notation of the previous footnote we
{\it erase} the index $\nabla^a$ and we make the index ${}_a$ into
a free index.} Furthermore recall that $b_v>\beta$ for each $v\in V$.

\par Now, applying the ``main algebraic proposition''
\ref{pregiade}\footnote{We just set $\Omega_1=\dots=\Omega_p=\phi$.}
 to the above, we
derive that there exists a linear combination of
 $(\beta+1)$-tensor fields \\$\sum_{h\in H} a_h C^{h,(i_1\dots
i_\beta)i_{\beta+1}}_g(\phi)$ in the form (\ref{paparigas}) so that:
\begin{equation}
\label{bastakardiamou2} \sum_{u\in U_\beta} a_u
\overline{pcontr}^{u,(i_1\dots i_\beta)}(\phi)=Xdiv_{i_{\beta+1}} \sum_{h\in H} a_h
C^{h,(i_1\dots
i_\beta)i_{\beta+1}}_g(\phi)+(Junk-Terms).
\end{equation}
Finally,  we formally replace each factor
$\nabla^{(m)}_{r_1\dots r_m}R_{ijkl}$ by a  factor
$\nabla^{(m)}_{r_1\dots r_m}W_{ijkl}$,\footnote{(times a constant,
which we ignore for these purposes)} and then make all the free
indices into internal contractions.\footnote{By this we mean that
for each free index ${}_{i_h}$ which belongs to a factor $T_{a\dots c}$, we
formally add a derivative $\nabla^{i_h}$ onto the factor $T_{a\dots c}$ and
contract it against the index ${}_{i_h}$; thus 
we obtain a factor $\nabla^{i_h}T_{a\dots c}$.} Denote this 
formal operation by $Weylify[\dots]$. Then 
$Weylify[\sum_{h\in H} a_h C^{h,(i_1\dots
i_\beta)i_{\beta+1}}_g(\phi)]$ is the divergence $div_iT^i(g)$ that is needed for  
(\ref{vakifleri}). 
\newpage

\section{Conventions, Background,
and the Super divergence formula from \cite{a:dgciI}.}

\subsection{Conventions and Remarks.}

\par We introduce some conventions that will be used throughout
this series of papers. Firstly, we recall two notions introduced in \cite{a:dgciI}:
\begin{definition}
 \label{subl.comb}
Given any (formal) linear combination $\sum_{l\in L} a_l C^l$
and any subset $L'\subset L$, then the
linear combination $\sum_{l\in L'} a_l C^l$ will be called a
{\it sublinear combination} of $\sum_{l\in L} a_l C^l(g)$.\footnote{In the 
introduction we spoke of a ``piece''
of $\sum_{l\in L} a_l C^l$, for simplicity.}
\end{definition}
 We also recall
 the notion of an ``internal contraction''
 for any (complete or partial) contraction: 
\begin{definition}
\label{int.contr}
Consider any
 complete or partial contraction $C(T_1,\dots,T_a)$, involving the
tensors $T_1,\dots,T_a$.  Then an ``internal contraction''
is a pair of indices ${}_i,{}_j$ that belong to the same
factor $T_h$ and are contracting against each other in $C(T_1,\dots,T_a)$.
\end{definition}
 Finally, as in \cite{a:dgciI} and \cite{a:dgciII} we define the
``length'' of a complete contraction to stand for the number of
its factors.

Now, a few minor conventions:

\par Firstly, when we say that a local equation, say $\sum_{t\in T} a_t C^t=\sum_{y\in Y} a_y C^y$
holds modulo terms of length $\ge\tau+1$, we will mean that there exists a linear 
combination of complete contractions with at least $\tau+1$ factors, $\sum_{f\in F} a_f C^f$,
such that $\sum_{t\in T} a_t C^t=\sum_{y\in Y} a_y C^y+\sum_{f\in F} a_f C^f$. 

\par Secondly, when we write $\nabla^{(m)}$, $m$ will stand
for the number of differentiations. When we write $\nabla^m$,
${}^m$ will stand for a raised index.\footnote{The reader should
note that this convention was {\it not} adopted in \cite{a:dgciI}
and \cite{a:dgciII}.}

Thirdly, we will often be referring to factors
 $\nabla^{(m)}R_{ijkl}$, $\nabla^{(p)}Ric_{ij}$ and $R$ (the third factor
being the scalar curvature)
 in complete and partial contractions below. Whenever we write
$\nabla^{(m)}R_{ijkl}$, we will be assuming that no two of the indices
 ${}_i,{}_j,{}_k,{}_l$ are contracting against each other
(unless stated otherwise). Also,
 in $\nabla^{(p)}Ric_{ij}$, no two of the indices ${}_i,{}_j$ will be
 contracting against each other (unless stated otherwise). Moreover,
for brevity we will not be explicitly writing out all the indices
that belong to the different terms. For example, when we refer to 
factors $\nabla^{(m)}R_{ijkl}$ we have written out the four lower indices
of the curvature tensor but not the covariant derivative indices.

\par Furthermore, throughout this paper we will often
write out complete contractions with {\it two or more} factors
$\nabla^{(m)}R_{ijkl}$ or $\nabla^{(p)}Ric_{ij}$. When we do so,
and hence have the indices ${}_{ijkl}$ or ${}_{ij}$ appearing
repeatedly as lower indices, we will {\it not} be assuming that
these indices are contracting against each other. I.e.~if we have
a factor $T=\nabla^{(m)}R_{ijkl}$ and $T'=\nabla^{(m')}R_{ijkl}$
appearing in the same complete contraction, the indices
${}_{ijkl}$ in $T$ and ${}_{ijkl}$ in $T'$ are {\it not} assumed
to be contracting against each other. We only use this notation to
avoid writing $\nabla^{(m)}R_{ijkl}$, $\nabla^{(m')}R_{i'j'k'l'}$,
$\nabla^{(m'')}R_{i''j''k''l''}$ etc. For each factor $\nabla^{(m)}R_{ijkl}$
 the indices ${}_i,{}_j,{}_k,{}_l$ will be called {\it internal indices};
 for each factor $\nabla^{(p)}Ric_{ab}$ the indices ${}_a,{}_b$
will be called {\it internal indices}.
\newline

{\it ``Mini-Appendices'':} Throughout this series of papers, we will
sometimes postpone the proof of certain claims; the reader will be
referred to ``Appendices'' or ``Mini-Appendices'' further down in
the paper. These Appendices often refer to very special cases of
more general claims which require special proofs; the
reader who is interested only in the broad ideas in these papers may
wish to circumvent these sections.

\subsection{ Background: Some useful formulas}

{\it Standard formulas:}\footnote{Unless mentioned otherwise,
these formulas come from \cite{a:dgciI}.}
The curvature tensor $R_{ijkl}$ of a Riemannian manifold is given
 by the formula:

\begin{equation}
\label{curvature}
[\nabla_i\nabla_j-\nabla_j\nabla_i]X_l=R_{ijkl}X^k.
\end{equation}

\par Moreover, the Ricci tensor $Ric_{ik}$ arises from $R_{ijkl}$
by contracting the indices ${}_j,{}_l$:

\begin{equation}
\label{ricci}
Ric_{ik}=R_{ijkl}g^{jl}.
\end{equation}
The Schouten tensor is a trace-adjustement of Ricci
curvature:

\begin{equation}
\label{WSch}
P_{ij}=\frac{1}{n-2}[{Ric_{ij}}-\frac{R}{2(n-1)}
g_{ij}].
\end{equation}
(Here $Ric_{ij}$ stands for Ricci curvature and $R$
stands for scalar curvature $R_{ijkl}g^{ik}g^{jl}$).

We also recall the Weyl tensor:

\begin{equation}
\label{Weyl} W_{ijkl}=
R_{ijkl}-[P_{jk}g_{il}+P_{il}g_{jk}-P_{jl}g_{ik}-P_{ik}g_{jl}],
\end{equation}
which is conformally invariant, i.e:

\begin{equation}
\label{transweyl}
W_{ijkl}(e^{2\phi} g)=e^{2\phi}{W}_{ijkl}(g).
\end{equation}

\par Furthermore, we recall the Cotton tensor:

\begin{equation}
\label{cotton}
C_{ijk}={\nabla}_kP_{ij}-{\nabla}_jP_{ik},
\end{equation}
which is related to the Weyl curvature in the following way:

\begin{equation}
\label{cottweyl}
\nabla^iW_{ijkl}=(3-n)C_{jkl}.
\end{equation}

 The Ricci curvature transforms as follows:

\begin{equation}
\label{ricci}
{Ric}_{ab}(e^{2\phi}g)={Ric}_{ab}(g)+
(2-n){\nabla}^{(2)}_{ab}\phi - {\Delta}\phi
g_{ab}+(n-2)(\nabla_a{\phi}\nabla_b{\phi}-\nabla^k{\phi}\nabla_k{\phi}
g_{ab}),
\end{equation}

  While the Schouten tensor has the following transformation law:
\begin{equation}
\label{wschtrans}
P_{ab}(e^{2\phi}g)=P_{ab}(g)-\nabla^{(2)}_{ab}{\phi}+
\nabla_a{\phi}
\nabla_b{\phi}-\frac{1}{2}\nabla^k\phi\nabla_k\phi g_{ab}.
\end{equation}

\par The curvature tensor transforms:

\begin{equation}
\label{curvtrans}
\begin{split}
&R_{ijkl}(e^{2\phi}g)=e^{2\phi}[R_{ijkl}(g)+ \nabla^{(2)}_{il}\phi
g_{jk}+\nabla^{(2)}_{jk}\phi g_{il}-\nabla^{(2)}_{ik}\phi g_{jl}-\nabla^{(2)}_{jl}\phi g_{ik}
\\&+\nabla_i\phi\nabla_k{\phi}g_{jl}+\nabla_j{\phi}\nabla_l{\phi}g_{ik}-\nabla_i{\phi}\nabla_l {\phi}
g_{jk} -\nabla_j{\phi}\nabla_k{\phi}g_{il}
\\&+|\nabla\phi|^2g_{il}g_{jk}- |\nabla\phi|^2g_{ik}g_{lj}].
\end{split}
\end{equation}
We also recall following transformation law for the
Levi-Civita connection under general conformal transformations
${\hat{g}}_{ij}(x)=e^{2\phi} g_{ij}(x)$:

\begin{equation}
\label{levicivita} {\nabla}_k {\eta}_l(e^{2\phi}g)=
{\nabla}_k{\eta}_l(g) -\nabla_k{\phi} {\eta}_l -\nabla_l{\phi} {\eta}_k
+\nabla^s{\phi} {\eta}_s g_{kl}.
\end{equation}

\par Finally, on certain rare occasions we will be using the transformation law
of the curvature tensor $R_{ijkl}$ {\it under variations of the metric
$g_{ij}$ by a symmetric 2-tensor $v_{ij}$:}

\begin{equation}
\label{tasia} \frac{d}{dt}|_{t=0}[R_{ijkl}(g_{ab}+tv_{ab})]=
\frac{1}{2}[\nabla^{(2)}_{il}v_{jk}+\nabla^{(2)}_{jk}v_{il}
-\nabla^{(2)}_{ik}v_{jl}-\nabla^{(2)}_{jl}v_{ik}]+Q(R,v),
\end{equation}
where $Q(R,v)$ stands for a quadratic expression involving the
curvature tensor $R_{abcd}$ and the 2-tensor $v_{ef}$.

 \subsection{The main consequence of the super
divergence formula.} \label{mainconseq}

In this subsection we codify a consequence of the {\it super divergence formula},
which was the main result in \cite{a:dgciI};
(We recall from the 10-page outline that this formula 
considers Riemannian operators $L_g(\phi)$, depending
both on the metric and auxiliary functions, whose integral
over any closed manifold is always zero, and expresses them 
as divergences of explicitly constructed vector fields). Here 
we codify  into a Lemma
a main consequence of this formula, which
is what we will mostly be using in this series of papers.
\newline

We start with some notation:

\par We will be considering 
 complete contractions
$C^l_{g}(\psi_1\dots ,\psi_Z)$ in the {\it normalized} form:

\begin{equation}
\label{ens}
\begin{split}
&contr(\nabla^{a_1\dots
a_t}\nabla^{(m_1)}R_{ijkl}\otimes\dots\otimes \nabla^{b_1\dots
b_u}\nabla^{(m_s)}R_{i'j'k'l'}\otimes
\\& \nabla^{c_1\dots c_v}\nabla^{(p_1)}Ric_{ij}\otimes\dots\otimes
\nabla^{d_1\dots d_z}\nabla^{(p_q)}Ric_{i'j'}\otimes R\otimes \dots
\otimes R\otimes
\\& \nabla^{w_1\dots w_y}\nabla^{(a_1)}\psi_1\otimes
\nabla^{x_1\dots x_\alpha}\nabla^{(a_Z)}\psi_Z),
\end{split}
\end{equation}
where we are making the following notational conventions: In the
factors \\ $\nabla^{a_1\dots a_t}\nabla^{(m)}_{r_1\dots r_m}
R_{ijkl}$ (we are using this generic notation for the first $s$
factors), each of the indices ${}^{a_1},\dots ,{}^{a_t}$
 is contracting against one of the indices
${}_{r_1}\dots ,{}_{r_m},{}_{i},{}_{j},{}_{k},{}_l$.
 Moreover, none of the indices ${}_{r_1}\dots ,{}_{r_m},{}_i,{}_j,{}_k,{}_l$ are
contracting between themselves.

\par For the next $q$ factors (in the generic form
$\nabla^{c_1\dots c_v}\nabla^{(p)}_{r_1\dots r_p}Ric_{ij}$),
each of the indices
${}^{c_1},\dots ,{}^{c_v}$ is contracting against one of the
indices ${}_{r_1},\dots ,{}_{r_p},{}_i,{}_j$ and also none
 of the indices ${}_{r_1},\dots ,{}_{r_p},{}_i,{}_j$ are
contracting between themselves.
Each factor $R$ is a scalar curvature term.

\par Finally, for the last $Z$ factors (which we denote
by the generic notation $\nabla^{w_1\dots w_y}
\nabla^{(a)}_{u_1\dots u_a}\psi$) we assume that each of the
indices ${}^{w_1},\dots , {}^{w_y}$ is contracting against one of
the indices ${}_{u_1},\dots ,{}_{u_a}$ and none of the indices
${}_{u_1},\dots ,{}_{u_a}$ are contracting between themselves.
Note that any complete contraction in the form $contr(\nabla^{(m)}R\otimes\dots\otimes\nabla^{(m')}R
\otimes\nabla^{(p)}\psi_1\otimes\dots\otimes\nabla^{(p')}\psi_h)$ 
can be expressed as a linear combination of contractions in the form (\ref{ens}), 
by just repeatedly applying the curvature identity and the second Bianchi identity.  
\newline

Consider a set $\{C^l_g(\psi_1,\dots,\psi_s) \}_{l\in L}$ of 
normalized complete contractions, indexed in $L$.
Let $L_M\subset L$ stand for an index set of complete
contractions
 $C^l_{g}(\psi_1,\dots ,\psi_Z)$ with a total of
$q+\rho=M$ factors $\nabla^{c_1\dots c_v}\nabla^{(p)}_{r_1\dots
r_p}Ric_{ij}$ and $R$. We assume that for some $M\ge 0$, all index
sets $L_s$ with $s>M$ are empty. We then denote by $C^l_{g}
(\psi_1,\dots ,\psi_Z,\Omega^M)$ the complete contraction that
formally arises
 from each $C^l_{g}(\psi_1,\dots ,\psi_Z)$, $l\in L_M$, by replacing each factor
$\nabla^{c_1\dots c_v}\nabla^{(p)}_{r_1\dots r_p}Ric_{ij}$ by a
factor $-\nabla^{c_1\dots c_v}\nabla^{(p+2)}_{r_1\dots
r_pij}\Omega$ and each factor $R$ by a factor $-2\Delta\Omega$.
($\Omega$ will be a scalar function).

\par (Thus $C^l_{g}(\psi_1,\dots ,\psi_Z,\Omega^{M})$, $l\in L_M$,
is a complete contraction in the form:

\begin{equation}
\label{ens2}
\begin{split}
&(-1)^M 2^\rho contr(\nabla^{a_1\dots
a_t}\nabla^{(m_1)}R_{ijkl}\otimes \dots\otimes \nabla^{b_1\dots
b_u}\nabla^{(m_s)}R_{i'j'k'l'}\otimes
\\& \nabla^{c_1\dots c_v}\nabla^{(p_1+2)}\Omega\otimes\dots\otimes
\nabla^{d_1\dots d_z}\nabla^{(p_q+2)}\Omega\otimes
\Delta\Omega\otimes \dots \otimes \Delta\Omega\otimes
\\& \nabla^{w_1\dots w_y}\nabla^{(a_1)}\psi_1\otimes
\nabla^{x_1\dots x_\alpha}\nabla^{(a_Z)}\psi_Z).)
\end{split}
\end{equation}

\begin{definition}
\label{eraseindex}
Consider any
internal contraction in $C^l_{g} (\psi_1,\dots
,\psi_Z,\Omega^{M})$, say $\zeta=({}^a,{}_a)$ (notice that ${}^a$
 must necessarily be a derivative index).
We 
then say that {\it we replace the internal contraction
$\zeta$ by a free index} if we {\it erase} the index ${}^a$
and make the index ${}_a$ into a free index.
\end{definition}

\par We thus obtain a
1-tensor field $(C^l)^{i_1}_{g}(\psi_1,\dots ,\psi_Z,\Omega^{M})$
of weight $-n+1$ (the free index ${}_{i_1}$ is the index ${}_a$ above).
The same formal definition can also be applied to $k$ internal
contactions: If we pick out $k$ internal contactions, say
 $({}^{a_1},{}_{a_1}),\dots,({}_{a_s},{}_{a_s})$ and then erase
the indices ${}^{a_1},\dots ,{}^{a_k}$ and make the indices
 ${}_{a_1},\dots,{}_{a_k}$ into free indices ${}_{i_1},\dots,{}_{i_k}$
 we obtain a
$k$-tensor field \\$(C^l)^{i_1\dots i_k}_{g}(\psi_1,\dots ,
\psi_Z,\Omega^{M})$ of weight $-n+k$.

This language convention (of making an
internal contraction into a free index)
will be used throughout this series of papers.

\begin{definition}
\label{spaier} 
 Now for each $l\in L_M$, we denote by
$(C^l)^{i_1\dots i_{b_l}}_{g}(\psi_1,\dots ,\psi_Z, \Omega^{M})$
the tensor field that arises from $C^l_{g} (\psi_1,\dots
,\psi_Z,\Omega^{M})$ by making all the internal contractions into
free indices. We denote by $$Xdiv_{i_1}\dots
Xdiv_{i_{b_l}}(C^l)^{i_1\dots i_{b_l}}_{g}(\psi_1, \dots
,\psi_Z,\Omega^{M})$$ the sublinear combination in $div_{i_1}\dots
div_{i_{b_l}}(C^l)^{i_1\dots i_{b_l}}_{g}(\psi_1, \dots
,\psi_Z,\Omega^{M})$ that arises when each $\nabla_{i_u}$ is
allowed to hit any factor other than the one to which ${}_{i_u}$
belongs.
\end{definition}

{\bf The main consequence of the super divergence formula:}

\begin{lemma}
\label{labrome} 
Assume an integral equation:

\begin{equation}
\label{rome} \int_{M^n} \Sum_{l\in L} a_l C^l_{g} (\psi_1,\dots
,\psi_Z)+ \Sum_{h\in H} a_h C^h_{g}(\psi_1,\dots ,\psi_Z)dV_{g}=0,
\end{equation}
which is assumed to hold for every compact $(M^n,g)$, and every
 $\psi_1,\dots,\psi_Z\in C^\infty (M^n)$. Here the
complete contractions indexed in $L$ have length $\sigma$ and are
in the {\it normalized form} (\ref{ens}), and the complete contractions indexed
in $H$ have length $>\sigma$. We let $M$ stand for the maximum
 number of factors $\nabla^{(p)}Ric$ and $R$ (in total) 
among the complete 
contractions $C^l_g(\dots), l\in L$; let $L_M\subset L$ be the index 
set of complete contractions with $M$ factors $\nabla^{(p)}Ric, R$ (in total).

We claim:
\begin{equation}
\label{szabo} \Sum_{l\in L_M} a_l (-1)^{b_l} Xdiv_{i_1}\dots
Xdiv_{i_{b_l}}(C^l)^{i_1\dots i_{b_l}}_{g}(\psi_1, \dots
,\psi_Z,\Omega^{M})=0,
\end{equation}
modulo complete contractions of length $\ge\sigma +1$.
\end{lemma}

{\it Proof:} We will show this claim in two steps. Initially, we
show
 that for some linear combination $\Sum_{h\in H'} a_h C^h_{g}
(\psi_1,\dots ,\psi_Z,\Omega^{M})$ of complete contractions with
length $\ge\sigma +1$ we have:

\begin{equation}
\label{siouda} \int_{M^n} \Sum_{l\in L_M} a_l C^l_{g}
(\psi_1,\dots ,\psi_Z,\Omega^{M})+ \Sum_{h\in H'} a_h
C^h_{g}(\psi_1,\dots ,\psi_Z,\Omega^{M})dV_{g}=0.
\end{equation}

{\it Proof of (\ref{siouda}):} Let us denote the integrand of (\ref{rome})
 by $L_{g}(\psi_1,\dots ,\psi_Z)$.

\par  We consider any dimension $N\ge n$ and denote by
$L_{g^N}(\psi_1,\dots ,\psi_Z)$ the re-writing of
$L_{g}(\psi_1,\dots ,\psi_Z)$ in dimension $N$.

Then, as shown in \cite{a:dgciI} (using the silly divergence formula) we
 derive that for any $N\ge n$,
 any $(M^N,g)$ and any $\psi_1,\dots ,\psi_Z\in
C^\infty(M^N)$:

\begin{equation}
\label{giannhs} \int_{M^N} L_{g^N}(\psi_1,\dots ,\psi_Z)dV_{g^N}=0.
\end{equation}

Now, let $L^{M}_{g^N}(\psi_1,\dots ,\psi_s,\Omega^{M}):=
\frac{\partial^{M}}{\partial\lambda^{M}}|_{\lambda
=0}[e^{(N-n)\lambda\Omega(x)}
L_{e^{2\lambda\Omega(x)}g^N}(\psi_1,\dots ,\psi_s)]$. It follows from (\ref{giannhs}) that:

\begin{equation}
\label{intNpsi} \int_{M^N} L^{M}_{g}(\psi_1,\dots
,\psi_Z,\Omega^{M}) dV_{g^N}=0.
\end{equation}

 Now, using (\ref{curvtrans}) and (\ref{levicivita})
 and the transformation law for
the volume form: $dV_{e^{2\lambda\Omega(x)}g}=
e^{N\lambda\Omega(x)}dV_{g}$, it follows that we can re-express
(\ref{intNpsi}) as follows:

\begin{equation}
\label{intNpsiana}
\begin{split}
& N^{M}\int_{M^N} \Sum_{l\in L_M} a_l C^l_{g^N}
(\psi_1,\dots ,\psi_Z,\Omega^{M})dV_{g^N}
\\&+ {\Sigma}_{x=0}^{M}
N^x\int_{M^N} \Sum_{u\in U^x} a_u C^u_{g}(\psi_1,\dots
,\psi_Z,\Omega^{M})dV_{g^N}=0,
\end{split}
\end{equation}
where the summands $C^u_{g^N}(\psi_1,\dots ,\psi_Z,\Omega^{M}),
u\in U^x$ are independent of the dimension $N$. Also, each
$C^u_{g^N}(\psi_1,\dots ,\psi_s,\Omega^{M}), u\in U^{M}$
 has at least $\sigma +1$ factors (possibly with factors
$\Omega$ without derivatives). Picking $M^N=M^n\times S^1\dots \times S^1$
with the product metric $g^N=g^n+(dt^1)^2+\dots +(dt^{N-n})^2$
we derive an integral equation in dimension $n$, where $N$
is just a free variable:

\begin{equation}
\label{intNpsianan}
\begin{split}
& N^{M}\int_{M^n} \Sum_{l\in L_M} a_l C^l_{g}
(\psi_1,\dots ,\psi_Z,\Omega^{M})dV_{g}
\\&+ {\Sigma}_{x=0}^{M}
N^x\int_{M^n} \Sum_{u\in U^x} a_u C^u_{g}(\psi_1,\dots
,\psi_Z,\Omega^{M})dV_{g}=0.
\end{split}
\end{equation}
Therefore, viewing the above as a polynomial in $N$  and
 restricting attention to the coefficient of $N^{M}$ we derive:

\begin{equation}
\label{concluNton} \int_{M^n} \Sum_{l\in L_M} a_l C^l_{g}
(\psi_1,\dots ,\psi_Z,\Omega^{M})dV_{g}+
 \int_{M^n} \Sum_{u\in U^{M}} a_u
C^u_{g}(\psi_1,\dots ,\psi_Z,\Omega^{M})dV_{g}=0,
\end{equation}
where each $C^u_{g}(\psi_1,\dots ,\psi_Z,\Omega^{M})$ has length
$\ge \sigma +1$. This is exactly (\ref{siouda}).
\newline

\par Now, we denote the integrand in (\ref{concluNton}) by
$L_{g} (\psi_1,\dots ,\psi_Z,\Omega^{M})$ and we apply the super
divergence
 formula to $L_{g}(\psi_1,\dots ,\psi_Z,\Omega^{M})$. We focus on the
sublinear combination $supdiv_{+}[L_{g} (\psi_1,\dots
,\psi_Z,\Omega^{M})]$ in $supdiv[L_{g} (\psi_1,\dots
,\psi_Z,\Omega^{M})]$ that consists of complete contractions of
length $\sigma$ and with no internal contractions. By virtue of
the super divergence formula and Lemma 8 in \cite{a:dgciI}, we
derive:

\begin{equation}
\label{katsis} supdiv_{+}[L_{g} (\psi_1,\dots
,\psi_Z,\Omega^{M})]=0,
\end{equation}
modulo complete contractions of length $\ge\sigma +1$.

\par On the other hand, by the algorithm for the super divergence
 formula in \cite{a:dgciI}, we derive:

\begin{equation}
\label{katsis} 
\begin{split}
&supdiv_{+}[L_{g} (\psi_1,\dots
,\psi_Z,\Omega^{M})]
\\&= \Sum_{l\in L_M} a_l (-1)^{b_l}
Xdiv_{i_1}\dots Xdiv_{i_{b_l}}(C^l)^{i_1\dots i_{b_l}}_{g}(\psi_1,
\dots ,\psi_Z,\Omega^{M}).
\end{split}
\end{equation}

Combining the two above equations we derive (\ref{szabo}). $\Box$

\section{From the super divergence formula for $I_{g}(\phi)$ back to
$P(g)$: The two main claims of this series of papers.}

\par Throughout this section, $P(g)$ will be a Riemannian scalar of weight $-n$
with the feature that $\int_{M^n}P(g)dV_{g}$ is a global conformal
invariant (see Definition \ref{globconfinv}).

\par Let us begin by writing $P(g)$ as a linear
 combination:

\begin{equation}
\label{pnorm} P(g)={\sum}_{l\in L} a_l C^l(g),
\end{equation}
where each complete contraction $C^l(g)$ is in the form:

\begin{equation}
\label{normcontrx}
\begin{split}
&contr({\nabla}^{(m_1)} W\otimes
\dots\otimes {\nabla}^{(m_s)}W
 \otimes {\nabla}^{(p_1)}P\otimes\dots \otimes
{\nabla}^{(p_q)}P).
\end{split}
\end{equation}

\par Our next two Propositions flesh out the claims
made in the first page of our ``10-page outline''. We will define 
 the ``worst piece'' in $P(g)$ and claim that
by subtracting a divergence
and a local conformal invariant we can cancel it out modulo
introducing ``better'' correction terms. The ``worst piece''
will consist of terms with a given number $\sigma$ of factors in
total and a given
number $s$ of factors $\nabla^{(a)}P$ (see the next paragraph).
The two propositions correspond to the cases $s>0$ and $s=0$.
\newline

\par Consider $P(g)$ as in (\ref{pnorm}). Denote by $\sigma$ the minimum number of
factors among the complete contractions indexed in $L$
in (\ref{pnorm}). Denote by $L_\sigma\subset L$ the index
set of those complete contractions. Now, denote by
$\Theta_r\subset L_\sigma$ the index set of complete
contractions with $r$ factors $\nabla^{(a)}P$ (and hence
$\sigma-r$ factors $\nabla^{(m)}W$). We note that some of these sets
may apriori be empty.

\begin{proposition}
\label{killtheta}
 Suppose that $P(g)=\Sum_{l\in L} a_l C^l(g)$ is a linear combination of
 contractions in the form (\ref{normcontrx}), and the minimum
  number of factors among the contractions
$C^l(g)$ is $\sigma<\frac{n}{2}$. We assume that (for $P(g)$) the
sets ${\Theta}_{\sigma},\dots ,{\Theta}_{s+1}$ are empty, where
$1\le s\le \sigma$, but $\Theta_s$ is not empty. We claim
that there
 is a Riemannian vector field $T^i(g)$ so that

\begin{equation}
\label{killthetat} {\Sigma}_{l\in {\Theta}_s} a_l C^l(g)- div_i
T^i(g)={\Sigma}_{r\in R} C^r(g),
\end{equation}
where each $C^r(g)$ is either in the form (\ref{normcontrx}) with
length $\sigma$ and fewer than $s$ factors $\nabla^{(p)}P$, or it
has length $>\sigma$.
\end{proposition}

\par Clearly, if we can show the above Proposition, then by iterative
 repetition we can derive that there is a vector field
$T^i(g)$ so that $P(g)-div_iT^i(g)=\Sum_{l\in L'} a_l
C^l(g)+\Sum_{j\in J} a_j C^j(g)$, with each $C^l(g), l\in L'$  in
the form:
\begin{equation}
\label{karami}
contr(\nabla^{(m_1)}W\otimes\dots\otimes\nabla^{(m_\sigma)}W),
\end{equation}
while each $C^j(g)$ will have at least $\sigma+1$ factors. Thus,
if we can show
 Proposition \ref{killtheta}, we will be reduced to proving Theorem \ref{thetheorem} 
in the case where all complete
 contractions in $P(g)$ with $\sigma$ factors are in the form
 (\ref{karami}).
\newline

\par Our next ``main claim'' applies precisely to that setting:

\begin{proposition}
\label{bigprop}
\par Consider any $P(g)$,
$P(g)=\Sum_{l\in L} a_l C^l(g)$ where each $C^l(g)$ has length
$\ge\sigma$, and each $C^l(g)$ of length $\sigma$ 
 is in the form (\ref{karami}). Denote by $L_\sigma\subset L$ 
the index set of terms with length $\sigma$.

We claim that there is a local conformal invariant $W(g)$
 of weight $-n$ and also a vector field $T^i(g)$ as in the statement
of Theorem \ref{thetheorem}, so that:

\begin{equation}
\label{karami2} \Sum_{l\in L^\sigma} a_l C^l(g)-W(g) -div_i
T^i(g)=0
\end{equation}
modulo complete contractions of length $\ge\sigma+1$.
\end{proposition}

\par We observe that if we can show the above two propositions,
then by iterative repetition our Theorem will follow, in
view of \cite{a:dgciII}.
\newline

\par Now, in the remainder of the present paper  we will
explain how to derive Proposition \ref{killtheta} in the case $\sigma \ge 3$ 
{\it assuming} the  ``main algebraic Proposition''  \ref{pregiade} below.
Proposition \ref{bigprop} in the case $\sigma \ge 3$ 
will be proven in \cite{alexakis2, alexakis3} {\it assuming} 
another two  ``main algebraic Propositions'' which are formulated in \cite{alexakis2}. 
The cases $\sigma <3$ of 
Propositions \ref{killtheta} and \ref{bigprop} will 
be proven in \cite{alexakis3}.
The three ``main algebraic propositions are then proven 
in \cite{alexakis4, alexakis5, alexakis6}.

\section{Proof of Proposition \ref{killtheta} in the easy case $s=\sigma$.}
\label{prfkillthet}

\par  We will
distinguish two cases: Either $s=\sigma$ or
$s<\sigma$.\footnote{Observe that in the case $s=\sigma$, all
complete contractions in $P(g)|_{\Theta_\sigma}$ contain only
factors $\nabla^{(p)}P$.} We will firstly show the claim when
$s=\sigma$. This proof is much easier than the case $s<\sigma$,
but it will contain simple forms of certain arguments that will be
used throughout this series of papers. It also is instructive,
in the sense that it can illustrate how the super divergence
formula applied to $I_g(\phi)$ can be used to
understand the algebraic structure of $P(g)$.

\begin{definition}
\label{sublinear} If $P(g)$ is in the form
$P(g)=\Sum_{l\in L} a_l C^l(g)$ then for any subset
$A\subset L$, we
 will denote by $P(g)|_A= {\sum}_{l\in A} a_l C^l(g)$.

\par Finally, for complete contractions $C(g)$, $C_g(\phi)$ of weight
$-n$,  we define the operation $Image_{\psi}^d$ as follows:

$$Image_{\psi}^d[C(g)]={\partial}^d_{\lambda}|_{\lambda =0}
\{ e^{n\lambda\psi(x)}C(e^{2\lambda\psi(x)})\},$$
and
$$Image_{\psi}^d[C_{g}(\phi)]={\partial}^d_{\lambda}|_{\lambda =0}\{ e^{n\lambda\psi(x)} C_{e^{2\lambda\psi(x)}
g(x)}(\phi)\}.$$
\end{definition}

\subsection{Proof of Proposition \ref{killtheta} when $s=\sigma$.}
\label{sofiavembo}

Our main tool will be to use the super divergence formula
applied to $I^\sigma_g(\phi)$ in order to show that the
sublinear combination $\sum_{s\in \Theta_\sigma} a_s C^s(g)$ in $P(g)$
is equal to a divergence modulo ``better'' correction terms.

\par Recall that $I_g(\phi):=e^{n\phi}P(e^{2\phi}g)-P(g)$ is the ``image''
of $P(g)$ under conformal variations of the metric $g$;
recall that $I^\sigma_g(\phi)$
consists of the terms in $I_g(\phi)$ which have
homogeneity $\sigma$ in the function $\phi$.

\par We have two tools at our disposal: Firstly,
 the ``super divergence formula'' for $I_g(\phi)$. Secondly, we
will momentarily show how the ``worst piece'' in $P(g)$ (i.e. 
the sublinear combination $\sum_{s\in \Theta_\sigma} a_s C^s(g)$--see
the discussion above Proposition \ref{killtheta})) is in
almost one-to-one corresponence with a
particular sublinear combination in $I^\sigma_g(\phi)$.

\par Let us flesh out the second remark: Observe that given the
formula (\ref{ens}) for $P(g)$, $I^\sigma_g(\phi)$ can be explicitly computed
by applying the identities (\ref{transweyl}), (\ref{wschtrans}) and
 (\ref{levicivita}). With a simple
observation we can derive much more
precise information:

 Since $I^\sigma_g(\phi)$ consists of terms of homogeneity
$\sigma$ in $\phi$ and the minimum number of factors in $P(g)$
is $\sigma$, we observe that the only complete contrqactions in $P(g)$
which can give rise to a term with $\sigma$ factors in $I_g(\phi$) are the ones
indexed in $\Theta_\sigma$. In fact, we can derive more: 
For each $s\in \Theta_\sigma$ we define $C^{s,\iota}_g(\phi)$
to stand for the complete contraction that arises from $C^s(g)$
by replacing each of the factors $\nabla^{(a)}_{r_1\dots r_a}P_{ij}$ by
$-\nabla^{(a+2)}_{r_1\dots r_a ij}\phi$. Then formulas
(\ref{wschtrans}), (\ref{transweyl}) and (\ref{levicivita}) imply that:

\begin{equation}
\label{iana} I^{\sigma}_{g}(\phi)= {\Sigma}_{s\in
{\Theta}_{\sigma}}
 a_s C^{s,\iota}_{g}(\phi)+ {\Sigma}_{k\in K} a_k
C^k_{g}(\phi),
\end{equation}
where each $C^k_{g}(\psi)$ is a complete contraction in the
 form:

\begin{equation}
\label{linisym}
\begin{split}
& contr(\nabla^{(m_1)}R_{ijkl}\otimes\dots\otimes\nabla^{(m)}R_{i'j'k'l'}\otimes
\nabla^{(p_1)}Ric_{ab}\otimes\dots\otimes\nabla^{(p_q)}Ric_{a'b'}\otimes
\\& R^\alpha\otimes\nabla^{(\nu_1)}\phi\otimes\dots\otimes\nabla^{(\nu_s)}\phi),
\end{split}
\end{equation}
with  length $\ge \sigma +1$.

\par Now, we are ready to prove our
 Proposition \ref{killtheta} in this case $s=\sigma$.
\newline

{\it Mini-outline of the proof of Proposition \ref{killtheta} when $s=\sigma$:} 
The proof relies strongly  on the super divergence formula.
We will show that this formula implies that
modulo terms with length $\ge\sigma+1$:

$$I^\sigma_g(\phi)=div_i \sum_{r\in R} a_r C^{r,i}_g(\phi),$$
where each vector field $C^{r,i}_g(\phi)$ (${}^i$ is the free index)
if a partial contraction in the form:
 $$pcontr(\nabla^{(a_1)}\phi\otimes\dots\otimes\nabla^{(a_\sigma)}\phi),$$
where each $a_i\ge 2$. We will
then show that the vector field $T^i(g)$ which formally
arises from $\sum_{r\in R} a_r C^{r,i}_g(\phi)$ by replacing
each factor $\nabla^{(a)}_{t_1\dots t_a}\phi$  by
$\nabla^{(a-2)}_{r_1\dots r_{a-2}}P_{r_{a-1}r_a}$ 
satisfies the claim of Proposition \ref{killtheta}.
\newline

{\it Proof of Proposition \ref{killtheta} when $s=\sigma$:} 
 We first consider the case where there is no complete
 contraction $C^s(g)$ with $s\in {\Theta}_{\sigma}$ which
contains a
 factor $P^i_i$. That implies that there is no complete
 contraction $C^{s,\iota}_{g}(\phi)$ in (\ref{iana}) with a
factor ${\nabla}^{(p)}_{r_1\dots r_p}\phi =\Delta \phi$. We will
 refer to this as the simplifying assumption.
\newline

{\it Proof of Proposition \ref{killtheta} (when $s=\sigma$), under
the simplifying
 assumption.}
\newline

\par We claim that there is a linear combination of
 vector-valued Riemannian differential operators in
 $\phi$,
$\{ C^{j_i}_{g}(\phi)\}_{i\in I}$, each in the form:

\begin{equation}
\label{tens1} pcontr({\nabla}^{({\nu}_1)}\phi\otimes\dots\otimes
{\nabla}^{({\nu}_{\sigma})}\phi),
\end{equation}
 with one free index and
${\nu}_1,\dots ,{\nu}_{\sigma}\ge 2$ so that
modulo complete contractions of length $\ge \sigma +1$:

\begin{equation}
\label{imagecanc1} {\Sigma}_{s\in {\Theta}_{\sigma}}
 a_s C^{s,\iota}_{g}(\phi)- {\Sigma}_{i\in I} a_i
div_{j_i} C^{j_i}_{g}(\phi)=0.
\end{equation}

{\it Proof of (\ref{imagecanc1}):} Recall the algorithm for the
super divergence formula from \cite{a:dgciI}. By Lemma 20 in
\cite{a:dgciI}, we only need to
 restrict our attention to the good, hard and undecided
 descendants of each $C^{s,\iota}_{g}(\phi)$, $s\in
{\Theta}_{\sigma}$. By Lemma 16 in \cite{a:dgciI},
these will all
 be $\vec{\xi}$-contractions $C^t_{g}(\phi,\vec{\xi})$ in
 the form:

\begin{equation}
\label{descphixi}
contr({\nabla}^{({\nu}_1-a_1)}\phi\otimes\dots\otimes
{\nabla}^{({\nu}_{\sigma}-a_{\sigma})}\phi\otimes\vec{\xi}\otimes
\dots\otimes\vec{\xi}),
\end{equation}
where each factor $\vec{\xi}$ contracts against a factor
${\nabla}^{({\nu}-a)}\phi$. Furthermore, since we have no
 factors $\Delta\phi$ in any $C^{s,\iota}_{g}(\phi)$,
$s\in {\Theta}_{\sigma}$, it follows that each $\nu -a\ge 2$. If a
$\vec{\xi}$-contraction above has $M$ factors $\vec{\xi}$, we
perform $M-1$ integrations by parts. The correction terms that we
introduce have length $\ge \sigma +1$.  So, we indeed derive
(\ref{imagecanc1}). $\Box$
\newline

\par We then construct Riemannian vector fields
$C^{i,j_i}(g)$ out of each
 Riemannian vector-valued differential operator
$C^{i,j_i}_{g}(\phi)$
 by substituting each factor ${\nabla}^{(\nu)}_{a_1\dots
a_{\nu}}{\phi}$ by a factor $-{\nabla}^{(\nu -2)}_{a_1\dots a_{\nu
-2}}P_{a_{\nu -1} a_{\nu}}$. We see that each
$div_{j_i}C^{j_i}(g)$ is a linear combination:

$$div_{j_i}C^{j_i}(g)=(-1)^{\sigma}
{\Sigma}_{s\in S^i}a_s C^s(g)$$ of complete contractions in the
form:

$$contr({\nabla}^{(m'_1)}P \otimes\dots\otimes {\nabla}^{(m'_{\sigma})} P).$$

{\it Derivation of Proposition \ref{killtheta} (with $s=\sigma$)
from (\ref{imagecanc1}):} We use the fact that (\ref{imagecanc1})
holds formally (see \cite{a:dgciI} for a definition of this notion)).
 We then {\it repeat} the sequence of permutations
of indices by which we make the linearization of the left hand
side of (\ref{imagecanc1}) {\it formally zero} to the linear
combination:

$$P(g)|_{\Theta_\sigma}-div_i\sum_{i\in I}C^{i,j_i}(g).$$
It follows that we can also make the above formally equal also, modulo
introducing correction terms by virtue of the identities
$\nabla_i\nabla_jX_l-\nabla_j\nabla_iX_l=R_{ijkl}X^k$ and
$\nabla_aP_{bc}-\nabla_bP_{ac}=\frac{1}{n-3}\nabla^dW_{abdc}$.

Observe that the correction terms that we obtain by virtue
of the above identities are precisely in the form allowed by our
Proposition \ref{killtheta}. This concludes the proof of our claim
in this case. $\Box$
\newline

{\it Proof of Proposition \ref{killtheta} (when $s=\sigma$) in the
general case (without the simplifying assumption).}
\newline

\par We now consider the case where the complete contractions
 $C^l(g)$, $l\in {\Theta}_{\sigma}$
 are allowed to contain factors $P^i_i$.

\par In this case we observe that if $C^l(g)$ contains $A$ factors $P^i_i$, then
$C^{l,\iota}_{g}(\phi)$ will contain $A$
 factors $\Delta \phi$. Recall the super divergence
 formula for $I^{\sigma}_{g}(\phi)$:

\begin{equation}
\label{ianaRpro} {\Sigma}_{l\in {\Theta}_{\sigma}} a_l
C^{l,\iota}_{g}(\phi)= {\Sigma}_{k\in K} a_k div_{i_k}
C^{i_k}_{g}(\phi),
\end{equation}
modulo complete contractions in the form (\ref{linisym}) with
 length $\ge \sigma +1$.

\par The problem is, now, that there might be vector fields
$C^{k,i_k}_{g}(\phi)$ which are in the form (\ref{tens1})
 with one free index and with a factor ${\nabla}_i\phi$.
Hence the procedure carried out for the previous simple case
 cannot be carried over to this case (because we can not replace the factors $\nabla\phi$, with only
 one derivative).
\newline

  So, in this case, we claim the following:

\begin{lemma}
\label{treatR} Consider (\ref{ianaRpro}).
There is a subset of the vector fields
$\{C^{k,i_k}_{g}(\phi)\}_{k\in K}$, indexed in $K^{\sharp}$,
$\{C^{k,i_k}_{g}(\phi)\}_{k\in K^\sharp}$, in the form
(\ref{tens1}), with the property that each $C^{k,i_k}_{g}(\phi)$,
$k\in K^\sharp$ contains
 factors ${\nabla}^{(l)}\phi$ with $l\ge 2$, so that:

\begin{equation}
\label{cancextpro} {\Sum}_{l\in {\Theta}_{\sigma}} a_l
C^{l,\iota}_{g} (\phi) - {\Sum}_{k\in K^{\sharp}} a_k div_{i_k}
C^{i_k}_{g}(\phi)=0
\end{equation}Thomas Watson 

modulo complete contractions of length $\ge \sigma +1$.
\end{lemma}

\par Let us notice that if we can prove the above,
we can then repeat the argument from the previous case, by using
(\ref{cancextpro}). Hence, we will have proven  our Proposition
\ref{killtheta} (when $s=\sigma$) in full generality.
\newline

{\it Proof of Lemma \ref{treatR}:} We will construct the set
$K^{\sharp}\subset K$.

\par We consider the  set of good, hard or undecided
 descendants (see the last definition in subsection 5.1 in 
\cite{a:dgciI} for a description of these notions) of the complete
 contractions $C^{l,\iota}_{g}(\phi), l\in
{\Theta}_{\sigma}$ with $\vec{\xi}$-length $\sigma$,
 and proceed to integrate by parts as explained in the
algorithm for the super divergence formula in \cite{a:dgciI}.
 We impose the restriction that any
factor $\vec{\xi}_i$ which contracts against a factor $\nabla\phi$
will not be integrated by parts, provided there
 is another factor $\vec{\xi}_i$ which does not contract against a factor $\nabla\phi$.
Furthermore, whenever along the iterative
 integration by parts we obtain a
$\vec{\xi}$-contraction  of $\vec{\xi}$-length $\sigma$ whose only
factors $\vec{\xi}$ contract against a factor ${\nabla}\phi$, we
cross it out and index it in the set $H$. The
$\vec{\xi}$-contractions that are not crossed out give rise to the
divergences $\{ a_k div_{i_k} C^{i_k}_{g}(\phi)\}_{k\in
K^{\sharp}}$.

\par Thus, we derive the equation:

\begin{equation}
\label{symplh} {\Sum}_{l\in {\Theta}_{\sigma}} a_l C^{l,\iota}_{g}
(\phi)+ {\Sum}_{k\in K^{\sharp}} a_k div_{i_k} C^{i_k}_{g}(\phi)+
PO[{\Sum}_{h\in H} a_h C^h_{g}(\phi,\vec{\xi})]=0,
\end{equation}
which holds modulo complete contractions of length $\ge \sigma
+1$. By construction, each vector field $C^{i_k}_{g}(\phi)$ has
 length $\sigma$ and is in the form (\ref{tens1})
with each $\nu_i\ge 2$.

Therefore, it suffices to show:

\begin{equation}
\label{symplh2} PO[{\Sum}_{h\in H} a_h C^h_{g}(\phi,\vec{\xi})]=0,
\end{equation}
modulo complete contractions of length $\ge \sigma +1$. Hence it
would suffice to show:

\begin{equation}
\label{symplheasy} {\Sum}_{h\in H} a_h C^h_{g}(\phi,\vec{\xi})=0,
\end{equation}
modulo $\vec{\xi}$-contractions of $\vec{\xi}$-length $\ge\sigma
+1$.

\par We do this as follows: Notice that in any
$\vec{\xi}$-contraction $C^h_g(\phi,\vec{\xi})$, $h\in H$, the function $\phi$
appears only in expressions $\nabla_i\phi\vec{\xi}^i$, or in factors
$\nabla^{(a)}\phi$ with $a\ge 2$. Let us consider the
$\vec{\xi}$-contraction with the maximum number $M$ of factors
${\nabla}_i\phi\vec{\xi}^i$. Suppose they are
 indexed in  $H^M\subset H$.
Notice that $M<\sigma$, otherwise we would have
$\sigma=\frac{n}{2}$ (we see this by considering the weight). If
we can show that:

\begin{equation}
\label{deixto} {\Sum}_{h\in H^M} a_h C^h_{g}(\phi,\vec{\xi})=0
\end{equation}
modulo $\vec{\xi}$-contractions of $\vec{\xi}$-length $\ge \sigma
+1$, then (\ref{symplheasy}) will follow by
 induction. We write each $C^h_{g}(\phi,\vec{\xi})$ with
$h\in H^M$ as follows:

$$C^h_{g}(\phi,\vec{\xi})={C'}^h_{g}(\phi)\cdot ({\nabla}_k\phi
\vec{\xi}^k)^M.$$

\par For any $h\in H^{M}$, we then define
$PO^{*}[C^h_{g}(\phi,\vec{\xi})]$ to stand for the sublinear
combination in $PO[C^h_{g}(\phi,\vec{\xi})]$ which arises as
follows: We integrate by parts with respect to each factor
$\vec{\xi}^k$ and then force each
  derivative ${\nabla}^k$ to hit a factor
${\nabla}^{(a)}_{r_1\dots r_a}\phi$ ($a\ge 2$) in
${C'}^h_{g}(\phi)$. We define

$$PO^{-}[C^h_{g}(\phi,\vec{\xi})]:=
PO[C^h_{g}(\phi,\vec{\xi})]- PO^{*}[C^h_{g}(\phi,\vec{\xi})].$$
Notice that by definition, each complete contraction of length
$\sigma$ in $PO^{-}[C^h_{g}(\phi,\vec{\xi})]$ will have strictly
fewer than $M$ factors $\nabla\phi$.

\par We write out the super divergence formula as follows:

\begin{equation}
\label{manylappro}\begin{split} &I^{\sigma}_{g}(\phi)+
{\Sigma}_{k\in  K^{\sharp}} a_k div_{i_k} C^{k,i_k}_{g}(\phi)+
{\Sigma}_{h\in H\setminus H^{M}}
 a_h PO[C^h_{g}(\phi,\vec{\xi})] +
\\& {\Sigma}_{h\in H^{M}} a_h  \{ PO^{*}
[C^h_{g}(\phi,\vec{\xi})] + PO^{-} [C^h_{g}(\phi,\vec{\xi})]\} =0,
\end{split}
\end{equation}
modulo complete contractions of length $\ge \sigma +1$.

\par Now, let us observe: Each complete contraction in
(\ref{manylappro}) that does not belong to ${\Sigma}_{h\in H^M}
a_h \{  PO^{*} [C^h_{g}(\phi,\vec{\xi})]\}$  will
 have fewer than $M$ factors $\nabla\phi$.
 This follows from the fact that $M$ is the maximum number of factors
${\nabla}^k \phi\vec{\xi}_k$ among the $\vec{\xi}$-contractions
$C^h_{g}(\phi,\vec{\xi})$, and since each complete
 contraction $C^{l,\iota}_{g}(\phi)$, $l\in
{\Theta}_{\sigma}$ and each vector field $C^{k,i_k}_{g}(\phi)$
have only factors ${\nabla}^{(a)}\phi$, $a\ge 2$, by construction.

\par We now claim that

\begin{equation}
\label{smax} {\Sigma}_{h\in H^{M}} a_h   PO^{*}
[C^h_{g}(\phi,\vec{\xi})] =0,
\end{equation}
modulo complete contractions of length $\ge \sigma +1$. This holds
because (\ref{manylappro}) holds formally, and since (\ref{smax}) is the sublinear
combination in (\ref{manylappro}) of complete contractions of
length $\sigma$ with $M$ factors $\nabla\phi$.

\par  Now, (\ref{smax}) also holds formally.
Write out:
$${\Sigma}_{h\in H^{M}} a_h   PO^{*}
[C^h_{g}(\phi,\vec{\xi})]={\Sigma}_{t\in T} a_t C^t_{g}(\phi),$$
where each complete contraction $C^t_{g}(\phi)$ is in the form:

\begin{equation}
\label{skata} contr({\nabla}^{(m_1)}_{r_1\dots r_{m_1}}\phi\otimes
\dots\otimes {\nabla}^{(m_{\sigma -M})}_{t_1\dots t_{m_{\sigma
-M}}}\phi\otimes {\nabla}_{y_1}\phi\otimes\dots\otimes
{\nabla}_{y_M}\phi).
\end{equation}
We observe that the linear combination ${\Sigma}_{t\in T} a_t
C^t_{g}(\phi)$ arises from the linear combination ${\Sigma}_{h\in
H^M} a_t {C'}^h_{g}(\phi)\cdot ({\nabla}_k\phi\vec{\xi}^k)^M$ by
making each factor $\vec{\xi}^k$ into a derivative ${\nabla}^k$,
then allowing the derivative $\nabla^k$ to hit any of the 
factors $\nabla^{(A)}\phi$ in $C'^h_g(\phi)$ and adding all the complete
 contractions we thus obtain.

In particular, each factor ${\nabla}_k\phi$ in any $C^t_{g}(\phi)$
contracts against a factor ${\nabla}^{(a)}\phi$, $a\ge 3$.

\par Now, for each $C^t_{g}(\phi)$
let $\tilde{C}^t_{g}(\phi)$ stand for the complete contraction of
weight $-n+2M$ which arises from $C^t_{g}(\phi)$ by erasing each
factor ${\nabla}_i\phi$ and also erasing the index against which
${}_i$ contracts. Since (\ref{smax}) holds formally, it follows
that:\footnote{A rigorous proof of this fact can be found 
in the Appendix below---see the operation $Erase$.}

\begin{equation}
\label{smax3ana2}{\Sigma}_{t\in T} a_t \tilde{C}^t_{g}(\phi)=0.
\end{equation}

 But (\ref{smax3ana2}) just tells us that:

$$\Sum_{h\in H^{M}} a_h (\sigma -M)^{M}{C'}^h_{g}(\phi)=0.$$

Therefore, we have shown (\ref{deixto}). $\Box$
\newline

We have fully proven the Proposition \ref{killtheta} when
$s=\sigma$. $\Box$

\section{Proposition \ref{killtheta} in the hard case (where $s<\sigma$).}
\label{reducindstat}

\subsection{Technical Tools:}

{\bf Useful identities:} Now, we will put down a few identities
that will prove useful later on.

{\it Decomposition of the Weyl tensor:}
 Recall the Weyl tensor $W_{ijkl}$, see (\ref{Weyl}). Consider the tensor
$T=\nabla^{r_{a_1}\dots r_{a_x}}
\nabla^{(m)}_{r_1\dots r_m}W_{ijkl}$
where each index ${}^{r_{a_s}}$ is contracting against the
(derivative) index ${}_{r_{a_s}}$, and all the other indices are
free. We have then introduced the language
convention that the tensor $T$ {\it has $x$ internal
contractions}.

\par We will decompose the tensor $T$ into a linear
combination of tensors in the form $\nabla^{(m)}R_{ijkl}$.
By just applying formula (\ref{Weyl}) we find:

\begin{equation}
\label{decompo1}
\begin{split} 
 &\nabla^{r_{a_1}\dots
r_{a_x}}\nabla^{(m)}_{r_1\dots r_m}W_{ijkl}=\nabla^{r_{a_1}\dots
r_{a_x}}\nabla^{(m)}_{r_1\dots r_m}R_{ijkl}
\\&+\Sum_{z\in
Z^{\delta=x+1}} a_z T^z(g)+\Sum_{z\in Z^{\delta=x+2}} a_z T^z(g),
\end{split}
\end{equation}
where $\Sum_{z\in Z^{\delta= x+1}} a_z T^z(g)$ stands for a linear
combination of tensor products of the form $\nabla^{r_{a_1}\dots
r_{a_x}}\nabla^{(m)}_{r_1\dots r_m}Ric_{sq}\otimes g_{vb}$ in the
same free indices as $T$, with the feature that there are a total
of $x+1$ internal contractions in the tensor
$\nabla^{(m)}Ric_{sq}$ (including the one in the tensor $Ric_{ab}={R^i}_{aib}$
itself). $\Sum_{z\in Z^{\delta= x+2}} a_z T^z(g)$ stands for a
linear combination of tensor products of the form
$\nabla^{r_{a_1}\dots r_{a_x}}\nabla^{(m)}_{r_1\dots r_m}R\otimes
g_{vb}\otimes g_{hj}$ ($R$ stands for the scalar curvature) in the
same free indices as $T$, with the feature that there are a total
of $x+2$ internal contractions in the tensor $\nabla^{(m)}R$
(including the two in the factor $R={R^{st}}_{st}$ itself). If $m>0$ we will use
the contracted second Bianchi identity to think of
$\nabla^{r_{a_1}\dots r_{a_x}}\nabla^{(m)}_{r_1\dots r_m}R$ as a
factor $2\nabla^{r_{a_1}\dots
r_{a_x}r_{a_{x+1}}}\nabla^{(m-1)}_{r_1\dots
r_{m-1}}Ric_{r_mr_{a_{x+1}}}$, modulo introducing quadratic
correction terms.

\par Next useful identity: We consider a factor $T$ in the form\\
$T=\nabla^{r_{a_1}\dots r_{a_x}}\nabla^{(m)}_{r_1\dots
r_m}W_{r_{m+1}r_{m+2}r_{m+3}r_{m+4}}$ where again each of the
indices ${}^{r_{a_v}}$ is contracting against the index
${}_{r_{a_v}}$, and moreover now at least one of the indices
${}^{r_{a_v}}$ is contracting against one of the internal indices
${}_{r_{m+1}},\dots ,{}_{r_{m+4}}$. We then calculate:

\begin{equation}
\label{decompo2} \begin{split} &\nabla^{r_{a_1}\dots
r_{a_x}}\nabla^{(m)}_{r_1\dots
r_m}W_{r_{m+1}r_{m+2}r_{m+3}r_{m+4}}=
\frac{n-3}{n-2}\nabla^{r_{a_1}\dots r_{a_x}}\nabla^{(m)}_{r_1\dots
r_m} R_{r_{m+1}r_{m+2}r_{m+3}r_{m+4}} \\&+\Sum_{z\in Z^{\delta=x}}
a_z T^z(g)+\Sum_{z\in Z^{\delta=x+1}} a_z T^z(g),
\end{split}
\end{equation} where $\Sum_{z\in Z^{\delta=x+1}} a_z T^z(g)$
stands for the same generic linear combination as before. 
$\Sum_{z\in
Z^{\delta=x}} a_z T^z(g)$ only appears in the case where there are
two indices ${}^{r_{a_b}}, {}^{r_{a_c}}$ contracting against two
internal indices in $W_{ijkl}$ (and moreover the indices
${}_{r_{a_b}},{}_{r_{a_c}}$ do not belong to the same block
${}_{[ij]},{}_{[kl]}$). $\Sum_{z\in
Z^{\delta=x}} a_z T^z(g)$ stands for a linear combination of tensors
$\nabla^{r_{a_1}\dots r_{a_{x-1}}}\nabla^{(m)}Ric_{ab}$ with $x$
internal contractions (also counting the internal contraction in
the factor $Ric_{ab}$ itself) , and with the extra feature that one of the
indices ${}^{r_{a_1}},\dots ,{}^{r_{a_{x-1}}}$ is contracting
against one of the internal indices ${}_a,{}_b$ in $Ric_{ab}$.
\newline

{\bf The ``fake'' second Bianchi identities for the
derivatives of the Weyl tensor:}
We recall that the Weyl tensor $W_{ijkl}$ is antisymmetric in the
indices ${}_i,{}_j$ and ${}_k,{}_l$, and also $W_{ijkl}=W_{klij}$. It also satisfies the
first Bianchi identity. Nevertheless, it does not satisfy the
second Bianchi identity. We now present certain substitutes
 for the second Bianchi identity:

Firstly, if the indices ${}_r,{}_i,{}_j,{}_k,{}_l$ are all free then:

\begin{equation} \label{in&out}
\nabla_rW_{ijkl}+\nabla_jW_{rikl}+\nabla_iW_{jrkl}=\sum
(\nabla^sW_{srty}\otimes g),
\end{equation}
where the symbol $\sum(\nabla^sW_{srty}\otimes g)$ stands for a
linear combination of a tensor product of the three-tensor
$\nabla^sW_{sqty}$ (i.e. essentially the Cotton tensor) with an
un-contracted metric tensor. The exact form of $\sum
(\nabla^sW_{sqty}\otimes g)$ is not important for our study so we do not write
it down.

\par On the other hand, if the indices ${}_i,{}_j,{}_k,{}_l$ are free we then
have:

\begin{equation}
\label{in&out2}
\nabla^s_sW_{ijkl}+\frac{n-2}{n-3}\nabla^s_jW_{sikl}+\frac{n-2}{n-3}\nabla^s_iW_{jskl}=
\Sum (W,g)+\sum Q(R),
\end{equation}
where the symbol $\Sum (W,g)$ stands for a linear combination of
tensor products: $\nabla^{ik}W_{iakb}\otimes g_{cd}$ ($g_{cd}$ is
an un-contracted metric tensor-note that there are {\it two}
internal contractions in the factor $\nabla^{ik}W_{ijkl}$) and the
symbol $\sum Q(R)$ stands for some linear combination of
quadratic expressions in the curvature tensor. Again the exact
form of these expressions is not important so we do not write them
down.

\par On the other hand, if the indices ${}_r,{}_i,{}_j,{}_l$ are free  then:

\begin{equation}
\label{in&out2c}
\nabla^k\nabla_rW_{ijkl}+\nabla^k\nabla_jW_{rikl}+\nabla^k\nabla_iW_{jrkl}=
\Sum Q(R).
\end{equation}

\par Furthermore, we have that the analogue of the second Bianchi identity clearly holds if both
the index ${}_r$ and one of the indices ${}_i,{}_j$ are involved in an
internal contraction:

\begin{equation}
\label{in&out2d}
\nabla^{ri}\nabla_rW_{ijkl}+\nabla^{ri}\nabla_jW_{rikl}+\nabla^{ri}\nabla_iW_{jrkl}=
\Sum Q(R).
\end{equation}

\par Lastly, we also note the identity:

\begin{equation}
\label{in&out3}
\Delta\nabla^kW_{ijkl}+\nabla^{rk}_iW_{irkl}+\nabla^{rk}_i
W_{rjkl}= \Sum Q(R).
\end{equation}

\par Let us also recall the identity:

\begin{equation}
\label{cotton}
\nabla_aP_{bc}-\nabla_bP_{ac}=\frac{1}{n-3}\nabla^dW_{abcd}.
\end{equation}

\par These identities will be useful in the context of the next
formal constructions. 
\newline

{\bf The operations ``$Weylify$'' and ``$Riccify$'':}
These two operations are formal operations that act 
on complete contractions in the forms (\ref{constr1}), (\ref{constr2}) and produce
complete contractions in the forms (\ref{liniPearly}), 
(\ref{liniPearly2}), respectively. We will show two important
technical Lemmas concerning these two formal operations,
Lemmas \ref{paparidis} and \ref{paparidis2}.

For the first construction, we will be considering complete
and partial contractions (with no internal
contractions) in the form:

\begin{equation}
\label{constr1}
contr(\nabla^{(m_1)}R_{ijkl}\otimes\dots\otimes\nabla^{(m_s)}R_{ijkl}\otimes
\nabla^{(p_1)}\psi\otimes\dots\otimes\nabla^{(p_q)}\psi\otimes\nabla\upsilon\otimes\dots
\otimes\nabla\upsilon)
\end{equation}
($\upsilon$ is a scalar) with the following restrictions: In each
complete contraction and vector field there are $a\ge 0$ factors
$\nabla\upsilon$ ($a$ is fixed) and $q$ factors $\nabla^{(p)}\psi$
($q$ also
 fixed). We require that none of the factors $\nabla\upsilon$ are
 contracting between themselves
and none of them contains a free index. Furthermore, we require that for any 
factor $\nabla^{(p)}\psi$ which is not contracting against a
 factor $\nabla\upsilon$, $p\ge 2$.

\begin{definition}
\label{weylify}
 We consider a collection of
 complete contractions, $\{C^l_{g}(\psi^q,\upsilon^a)\}_{l\in
L}$ and a collection of such vector fields
$\{C^{r,i}_{g}(\psi^q,\upsilon^a)\}_{r\in R}$ in the form
(\ref{constr1}). Assume that the complete
contractions and vector fields above all have a given length $\tau
+a$.

 We define an operation $Weylify[\dots ]$ that acts on the
contractions and vector fields above by performing the following
operations: Each factor $\nabla^{(m)}_{r_1\dots r_m}R_{ijkl}$ that
is not contracting against a factor $\nabla\upsilon$ is replaced
by a factor $\nabla^{(m)}_{r_1\dots r_m}W_{ijkl}$. Also, each
factor $\nabla^{(p)}_{r_1\dots r_p}\psi$ that is not contracting
against against any factor $\nabla\upsilon$ is replaced by a
factor $\nabla^{(p-2)}_{r_1\dots r_{p-2}}P_{r_{p-1}r_p}$.

\par Now, any factor $T=\nabla^{(m)}_{r_1\dots
r_m}R_{r_{m+1}r_{m+2}r_{m+3}r_{m+4}}$ that is contracting against
$s>0$ factors $\nabla\upsilon$, {\it with the restriction that
all these $s$ factors $\nabla\upsilon$ are contracting against
derivative indices} will be replaced as follows: Suppose it is the
indices ${}_{r_{a_1}},\dots ,{}_{r_{a_x}}$ that are contracting
against the factors $\nabla\upsilon$. Then, we replace $T$ by
$\nabla^{r_{a_1}\dots r_{a_x}}\nabla^{(m)}_{r_1\dots
r_m}W_{r_{m+1}r_{m+2}r_{m+3}r_{m+4}}$.

\par On the other hand, if there are
internal indices also contracting against factors
$\nabla\upsilon$, we replace $T$ by
$\frac{n-2}{n-3}\nabla^{r_{a_1}\dots
r_{a_x}}\nabla^{(m)}W_{ijkl}$. Now, each factor
$\nabla^{(p)}_{r_1\dots r_p}\psi$ with $p\ge 2$ that is
contracting against $w$ factors $\nabla\upsilon$ (say the indices
$r_{a_1}\dots r_{a_w}$) is replaced by $\nabla^{r_{a_1}\dots
r_{a_w}}\nabla^{(p-2)}_{r_1\dots r_{p-2}}P_{r_{p-1}r_p}$. Finally,
every expression $\nabla_i\psi\nabla^i\upsilon$ is replaced by a
factor $P^a_a$. In the end, we also erase all the factors
$\nabla\upsilon$ (they were left uncontracted).
\end{definition}

 Thus, by acting on the complete
 contractions and vector fields in the form (\ref{constr1})
 with the operation $Weylify[\dots ]$, we obtain complete
 contractions and vector fields of length $\tau$ in the form:

\begin{equation}
\label{liniPearly}
\begin{split}
&contr(\nabla^{f_1\dots f_y}\nabla^{(m_1)}W_{ijkl}\otimes\dots
\otimes \nabla^{u_1\dots u_p}\nabla^{(m_s)}W_{i'j'k'l'}
\\& \nabla^{a_1\dots
a_{t_1}}\nabla^{(u_1)}P_{ij}\otimes\dots\otimes \nabla^{c_1\dots
c_{t_s}}\nabla^{(u_z)}P_{i'j'}\otimes P^a_a\otimes\dots\otimes P^b_b),
\end{split}
\end{equation}
where we are making the following conventions: In each factor
$\nabla^{f_1\dots f_y}\nabla^{(m)}_{r_1\dots r_m}W_{ijkl}$
 each of the the indices ${}^{f_1},\dots ,{}^{f_y}$ contracts
 against one of the indices ${}_{r_1},\dots ,{}_l$, while no two of
the indices ${}_{r_1},\dots ,{}_l$ contract between themselves. On
the
 other hand, for each factor $\nabla^{y_1\dots y_t}
\nabla^{(u)}_{a_1\dots a_u}P_{ij}$, each of the
 upper indices ${}^{y_1},\dots ,{}^{y_t}$ contracts
against one of the indices ${}_{a_1},\dots ,{}_{a_u},{}_i,{}_,{}_j$.
Moreover, none of the indices ${}_{a_1},\dots ,{}_{a_u},{}_i,{}_j$
contract between themselves.

\begin{definition}
\label{deltasetc} Consider any complete contraction (or tensor
field) of
 the form (\ref{liniPearly}), with the properties described above. We
will let $\delta_W$ stand for the number of internal
contractions in all the factors $\nabla^{(m)}W_{ijkl}$, and the
number $\delta_P$ to stand for the number of internal contractions
among all the factors $\nabla^{(p)}P_{ij}$.
\end{definition}

 We see that for a
contraction or vector field $C_{g}(\psi^q,\upsilon^a)$ in the form
(\ref{constr1}), the complete contraction or vector field
$Weylify[C_{g}(\psi^q,\upsilon^a)]$ will have length $\tau$ and a
total of $q$ factors in the form $\nabla^{(p)}P$, and
$\delta_W+\delta_P=a$. This operation extends to linear
combinations of contractions.
One last definition prior to stating our Lemma:

\begin{definition}
 \label{Xdiv}
 For any vector field $C^i_g$ in the form
(\ref{constr1}), $Xdiv_iC^i_g$ will stand for the sublinear
combination in $div_i C^i_g$ where $\nabla_i$ is not allowed to
hit the factor to which the free index ${}_i$ belongs, {\it nor}
any of the factors $\nabla\upsilon$.
\end{definition}

\par Now, our claim regarding the operation $Weylify$ is
 the following:

\begin{lemma}
\label{paparidis} Assume an equation:

\begin{equation}
\label{upo9eto} \Sum_{l\in L} a_l
C^l_{g}(\psi^q,\upsilon^a)-Xdiv_i \Sum_{r\in R} a_r
C^{r,i}_{g}(\psi^q,\upsilon^a)=0,
\end{equation}
that holds modulo complete contractions of length $\ge\tau+a+1$. Here the contractions and
 tensor fields are in the form (\ref{constr1}) with length $\tau +a$.
 We claim:

\begin{equation}
\label{kataka9ia} 
\begin{split}
&\Sum_{l\in L} a_l
Weylify[C^l_{g}(\psi^q,\upsilon^a)]-div_i \Sum_{r\in R} a_r
Weylify[C^{r,i}_{g}(\psi^q,\upsilon^a)]
\\&= \Sum_{d\in D^1} a_d
C^d(g)+\Sum_{d\in D^2} a_d C^d(g),
\end{split}
\end{equation}
where each $C^d(g)$ is in the form (\ref{liniPearly}) (with length
$\tau$) and moreover if $d\in D^1$ then $C^d(g)$ has
less than $q$ factors $\nabla^{(p)}P$, while if $d\in D^2$ then
$C^d(g)$ has $q$ factors $\nabla^{(p)}P$ but
 also
$\delta_W+\delta_P\ge a+1$. This equation holds modulo complete
contractions of length $\ge\tau +1$.
\end{lemma}

{\it Proof:} We will use the fact that (\ref{upo9eto})
holds formally to {\it repeat} the formal applications
of identities that make the LHS  of (\ref{upo9eto}) formally zero to the LHS of
(\ref{kataka9ia}); the RHS of (\ref{kataka9ia}) will
then arise as correction terms in this process.
 Now, we first observe that it would be sufficient
to show that

\begin{equation}
\label{billy}Weylify\{\Sum_{l\in L} a_l
C^l_{g}(\psi^q,\upsilon^a)-Xdiv_i \Sum_{r\in R} a_r
C^{r,i}_{g}(\psi^q,\upsilon^a)\}
\end{equation}
 is equal to the right hand side
of (\ref{kataka9ia}). That this is sufficient is clear because the
contraction that arises in each
$$div_i
Weylify[C^{r,i}_{g}(\psi^q,\upsilon^a)]$$ when $\nabla^i$ hits the
factor to which ${}_i$ belongs is clearly in the form $C^d(g),
d\in D^2$, and moreover because for each $r\in R$
$$Weylify\{Xdiv_i
C^{r,i}_{g}(\psi^q,\upsilon^a)\}=Xdiv_i
Weylify\{C^{r,i}_{g}(\psi^q,\upsilon^a)\},$$ modulo contractions of
length $\ge\tau +1$.

\par Next, we use the fact that (\ref{upo9eto}) holds exactly
(with no correction terms) at the linearized level (i.e.~if we
replace each complete contraction $C_{g}(\psi^q,\upsilon^a)$ by
$linC_{g}(\psi^q,\upsilon^a)$).\footnote{See the introduction of 
 \cite{a:dgciI} for a definition 
of linearization.} We ``memorize'' the sequence of
permutations of indices (and applications of the distributive
rule) by which we can make the linearization of (\ref{upo9eto})
formally zero. We may then repeat the same sequence of
permutations to the left hand side of (\ref{billy}), to make it
vanish, modulo introducing correction
terms, as follows:

\begin{enumerate}
\item{ We introduce correction terms of length $\ge\tau +1$ by
virtue of (\ref{curvature}) when we permute derivative indices in
a factor $\nabla^{(m)}W_{ijkl}$ or when we permute the first $p-2$
derivative indices in a factor $\nabla^{(p-2)}P_{ij}$.}

\item{We introduce correction terms of the form $\Sum_{d\in D^2}
a_d C^d(g)$ by virtue of (\ref{in&out}) when we apply the ``fake''
second Bianchi identity to the indices ${}_{r_m},{}_i,{}_j$ in a
factor $\nabla^{(m)}_{r_1\dots r_m}W_{ijkl}$ with
 no internal contractions involving internal indices.}

\item{We introduce correction terms of length $>\tau$ or of the
form $\Sum_{d\in D^2} a_d C^d(g)$, by virtue of the identities
(\ref{in&out2}), (\ref{in&out3}) when we  apply the ``fake'' second
Bianchi identity to the indices ${}_{r_m},{}_i,{}_j$ in a factor
$\nabla^{(m)}_{r_1\dots r_m}W_{ijkl}$ with one or two internal
contractions respectively.}

\item{We introduce correction terms of the form $\Sum_{d\in D^1}
a_d C^d(g)$ from the right hand side of (\ref{cotton}) when we
want to switch the indices ${}_{r_{m-2}},{}_{r_{m-1}}$ in a factor
$\nabla^{(p)}_{r_1\dots r_{m-2}}P_{r_{m-1}r_m}$.}
\end{enumerate}

 That completes the proof of our claim. $\Box$
\newline

{\it The operation $Riccify$:} We now define the operation $Riccify$ that acts on
complete contractions $C(\Omega^q,\psi^s,\upsilon^a)$ and vector
fields $C^i(\Omega^q,\psi^s,\upsilon^a)$ in the form:

\begin{equation}
\label{constr2}
\begin{split}
&contr(\nabla^{(m_1)}R_{ijkl}\otimes\dots\otimes\nabla^{(m_t)}
R_{ijkl}\otimes\nabla^{(a_1)}\psi_1\otimes\nabla^{(a_s)}\psi_s
\otimes
\\&\nabla^{(p_1)}\Omega\otimes\dots\otimes\nabla^{(p_q)}\Omega\otimes\nabla\upsilon\otimes\dots
\otimes\nabla\upsilon)
\end{split}
\end{equation}
with length $\tau +a$ (and with $a$ factors $\nabla\upsilon$),
 where both the $s$ factors $\nabla^{(u)}\psi_h$ and 
the $q$ factors $\nabla^{(p)}\Omega$
 are subject
 to the same restrictions as for the factors $\nabla^{(p)}\psi$
in the contractions in the form (\ref{constr1}). In particular:
 In each
complete contraction and vector field in the above form there are
$a\ge 0$ factors $\nabla\upsilon$ ($a$ is fixed) and $q$ factors
$\nabla^{(p)}\psi$ ($q$ also
 fixed). Also, none of the factors $\nabla\upsilon$ are
 contracting between themselves
and none of them contains a free index. Also, we require that any
factor $\nabla^{(p)}\Omega$ or $\nabla^{(p)}\psi$ which is not
contracting against a  factor $\nabla\upsilon$ must have $p\ge 2$.
 Moreover, we assume that the complete
contractions and vector fields above all have a fixed length $\tau
+a$.

\begin{definition}
\label{riccify}
 We define an operation $Riccify[\dots]$ that acts on complete
and partial contractions
  in the form (\ref{constr2}) as follows:
 We replace each factor \\$\nabla^{(m)}_{r_1\dots
r_m}R_{r_{m+1}r_{m+2}r_{m+3}r_{m+4}}$ for which the indices
${}_{r_{a_1}},\dots ,{}_{r_{a_d}}$ are contracting against factors
$\nabla\upsilon$ by a factor $\nabla^{r_{a_1}\dots
r_{a_d}}\nabla^{(m)}_{r_1\dots
r_m}R_{r_{m+1}r_{m+2}r_{m+3}r_{m+4}}$ We also replace each factor
$\nabla^p_{r_1\dots r_p}\Omega$ $(p\ge 2$) for which the indices
${}_{r_{a_1}},\dots ,{}_{r_{a_d}}$ are contracting against factors
$\nabla\upsilon$ by a factor $\nabla^{r_{a_1}\dots
r_{a_d}}\nabla^{(p)}_{r_1\dots r_{p-2}}Ric_{r_{p-1}r_p}$. Then, we
replace all expressions $\nabla_i\Omega\nabla^i\upsilon$ by a
factor $\frac{1}{2}R$. Finally, we replace each factor
$\nabla^{(a)}_{r_1\dots r_a}\psi_h$ (for which the indices
$r_{v_1},\dots r_{v_c}$ are contracting against factors
$\nabla\upsilon$) by an expression $\nabla^{r_{v_1}\dots
r_{v_c}}\nabla^{(a)}_{r_1\dots r_a}\psi_h$.
 In the end we also erase all the factors
$\nabla\upsilon$ (they have been left uncontracted).
\end{definition}

Thus acting by the
operation $Riccify$ on complete and partial contractions in
the form (\ref{constr2}) we obtain complete and partial contractions in
 the form:

\begin{equation}
\label{liniPearly2}
\begin{split}
&contr(\nabla^{f_1\dots f_y}\nabla^{(m_1)}R_{ijkl}\otimes\dots
\otimes \nabla^{u_1\dots u_p}\nabla^{(m_t)}R_{ijkl}\otimes\nabla^{v_1\dots v_z}\nabla^{(a_1)}\psi_1
\\& \otimes
\dots\otimes \nabla^{q_1\dots q_w}\nabla^{(a_s)}\psi_s\otimes
\nabla^{a_1\dots a_{t_1}}\nabla^{(u_1)}Ric_{ij}\otimes\dots\otimes
\nabla^{c_1\dots c_{t_s}}\nabla^{(u_z)}Ric_{ij}).
\end{split}
\end{equation}

\begin{definition}
\label{eisenbud} For contractions in the form (\ref{liniPearly2})
we define $\delta_R$ to stand for the total number of internal
contractions in the factors $\nabla^{(m)}R_{ijkl}$ and
$\delta_{Ric}$ to stand for the total number of internal
contractions in the factors $\nabla^{(p)}Ric$ (including the one
in the factors $Ric$ themselves) and also $\delta_\psi$ to stand
for the total number of internal contractions in the factors
$\nabla^{(a)}\psi$.
\end{definition}
 (Note: In the
future we will sometimes denote this operation $Riccify$  by
$\Omega to Ric$).

{\it Note:} In (\ref{constr2}), we may have $s=0$.
Furthermore, we recall from Definition \ref{Xdiv} that if $C^i_g$ is a vector field in the form
(\ref{constr2}) then $Xdiv_i$ will stand for the sublinear
combination in $div_i C^i_g$ where $\nabla_i$ is not allowed to
hit the factor to which the free index ${}_i$ belongs, {\it nor}
any of the factors $\nabla\upsilon$.

 Our Lemma is then the following:

\begin{lemma}
\label{paparidis2} Assume an equation:

\begin{equation}
\label{upo9eto2} \Sum_{l\in L} a_l
C^l_{g}(\Omega^q,\psi^s,\upsilon^a)-Xdiv_i \Sum_{r\in R} a_r
C^{r,i}_{g}(\Omega^q,\psi^s,\upsilon^a)=0,
\end{equation}
which holds modulo complete contractions of length
$\ge\tau+a+1$.

\par We claim:

\begin{equation}
\label{kataka9ia'}
\begin{split}
 & \Sum_{l\in L} a_l
Riccify[C^l_{g}(\Omega^q,\phi^s,\upsilon^a)]-div_i \Sum_{r\in R}
a_r Riccify[C^{r,i}_{g}(\Omega^q,\psi^s,\upsilon^a)]
\\&= \Sum_{d\in
D^1} a_d C^d_{g}(\psi^s)+\Sum_{d\in D^2} a_d C^d_{g}(\psi^s),
\end{split}
\end{equation}
where each $C^d_{g}(\phi^s)$ is in the form (\ref{liniPearly2})
(with length $\tau$) and moreover if $d\in D^1$ we will have
 that
$C^d_{g}(\phi^s)$ has less than $q$ factors $\nabla^{(p)}Ric$
 but will also have $\delta_R+\delta_{Ric}+\delta_\psi\ge a$,
while if $d\in D^2$ then $C^d_{g}(\phi^s)$ has $q$ factors
$\nabla^{(p)}Ric$ but also $\delta_R+\delta_{Ric}+ \delta_\psi\ge
a+1$. This equation holds modulo complete contractions of length
$\ge\tau+1$.
\end{lemma}

{\it Proof:} The proof is an easier version of the proof of the
previous Lemma. We use the fact that (\ref{upo9eto2}) holds
formally and we {\it repeat} the applications of the formal
identities and the distributive rule that make (\ref{upo9eto2})
formally zero to the LHS of (\ref{kataka9ia'}).

\par Now, in (\ref{kataka9ia'}), we use the identity
$\nabla_iR=2\nabla^kRic_{ik}$ ($R$ here is the scalar curvature)
once if needed, and we may assume that all the complete
contractions in the LHS of (\ref{kataka9ia'}) have any factors
$\nabla^{(p)}Ric$ (i.e. factors involving the Ricci curvature)
being in the form $\nabla^{f_1\dots f_b}\nabla^{(p)}_{r_1\dots
r_p}Ric_{ab}$ (where each of the indices ${}^{f_1},\dots,{}^{f_b}$
is contracting against one of the indices
${}_{r_1},\dots,{}_{r_p},{}_a,{}_b$, and none of the lower indices
are contracting between themselves), or in the form $R$ (scalar
curvature).

 Furthermore, when we {\it repeat} the permutations by which the
 LHS of (\ref{upo9eto2}) is made formally zero to the LHS of
(\ref{kataka9ia'}), we may assume wlog that the upper indices in
each factor  $\nabla^{f_1\dots f_b}\nabla^{(p)}_{r_1\dots
r_p}Ric_{ab}$ are {\it not} permuted (since they correspond to
factors $\nabla\upsilon$ in the LHS of (\ref{upo9eto2})).

\par Therefore, the RHS in (\ref{kataka9ia'}) can arise either when the
divergence index $\nabla_i$ in $div_i
Riccify[C^{r,i}_{g}(\Omega^q,\psi^s,\upsilon^a)]$ hits the factor
to which ${}_i$ belongs, {\it or} by virtue of the identity:
$$\nabla_aRic_{bc}-\nabla_bRic_{ac}=\nabla^dR_{abcd}$$
(where by the observation above the indices ${}_a,{}_b,{}_c$ will
not be contracting against each other).  $\Box$

\subsection{Proof of Proposition \ref{killtheta} when $s<\sigma$: Reduction to
an inductive statement.}

\par In the rest of this section we will explain how to derive
 Proposition \ref{killtheta} in the case where $\sigma\ge 3$.
 The cases $\sigma=1,\sigma =2$ will be covered in the
 paper \cite{alexakis3} in this series.

Recall (see the discussion above Definition \ref{deltasetc}) 
that we are assuming that for $P(g)$, $\Theta_s\ne
\emptyset$ and $\Theta_h=\emptyset$ for each $h>s$. We write
$P(g)|_{\Theta_s}$ as a linear combination:

\begin{equation}
\label{avni} P(g)|_{\Theta_s}=\Sum_{l\in L} a_l C^l(g)
\end{equation}
(modulo longer complete contractions), where each $C^l(g)$ is of
the form:

\begin{equation}
\label{liniP}
\begin{split}
&contr(\nabla^{f_1\dots f_y}\nabla^{(m_1)}_{r_1\dots
r_{m_1}}W_{ijkl}\otimes\dots \otimes \nabla^{s_1\dots
s_z}\nabla^{(m_f)}_{t_1\dots t_{m_f}}W_{i'j'k'l'}
\\&\otimes
\nabla^{y_1\dots y_t}\nabla^{(u_1)}_{a_1\dots
a_u}P_{ij}\otimes\dots\otimes \nabla^{y'_1\dots y'_o}
\nabla^{(u_r)}_{b_1\dots b_o}P_{i'j'}\otimes (P^a_a)^K),
\end{split}
\end{equation}
with the usual conventions: In each factor
$\nabla^{f_1\dots f_y}\nabla^{(m)}_{r_1\dots r_m}W_{ijkl}$  each of the the
raised indices ${}^{f_1},\dots ,{}^{f_y}$
contracts against one of the indices ${}_{r_1},\dots ,{}_l$, while
no two of the indices ${}_{r_1},\dots ,{}_l$ contract between
themselves. On the other hand,
 in each factor $\nabla^{y_1\dots y_t}\nabla^{(u)}_{a_1\dots a_u}
P_{ij}$,  each of the raised indices
${}^{y_1},\dots ,{}^{y_t}$ contracts
against one of the indices ${}_{a_1},\dots ,{}_{a_u},{}_i,{}_j$. Moreover, none of
the indices ${}_{a_1},\dots ,{}_{a_u},{}_i,{}_j$ contract between themselves. We
call such complete contractions $W$-{\it normalized}.

\par By virtue of the curvature identity it is clear
that modulo introducing correction terms of length
$\ge\sigma+1$, we can write $P(g)|_{\Theta_s}$ as a linear
 combination of $W$-normalized complete contractions
$C^l(g)$.
\newline

\begin{definition}
\label{tourkika}
\par Now, for each complete contraction $C^l(g)$ in the form (\ref{liniP}),
we defined $\delta_W$ to stand for the number of internal
contractions among the factors $\nabla^{(m)}W_{ijkl}$. We defined
$\delta_P$ to stand for the number of internal contractions among
the factors $\nabla^{(p)}P_{ij}$ {\it plus} the number $K$ of
factors $P^a_a$. In order to distinguish these numbers among the
various complete contractions $C^l(g)$, $l\in L$, we will write
$\delta_W(l),\delta_P(l)$. We also define
$\delta(l)=\delta_W(l)+\delta_P(l)$ (sometimes we will write
$\delta$ instead of $\delta(l)$). This notation trivially extends
to vector fields in the form (\ref{liniP}) with one free index.
\end{definition}

\par Furthermore, in the cases $s=\sigma-1$ and $\sigma-2$ we will
introduce an extra piece of notation purely for technical reasons:

{\it Special definition:}   If
 $s=\sigma-2$ then $P(g)|_{\Theta_s}$ is ``good'' if the only complete
 contraction in $P(g)|_{\Theta_s}$ with $\sigma-2$ factors $P^a_a$
 is of the form $(Const)\cdot\Delta^{\frac{n}{2}-\sigma-2}\nabla^{il}
 W_{ijkl}\otimes\nabla^{i'l'}{{W_{i'}}^{jk}}_{l'}\otimes(P^a_a)^{\sigma-2})$
 (when $\sigma<\frac{n}{2}-1$)
 or $(Const)\cdot contr(\nabla^lW_{ijkl}\otimes\nabla_{l'}W^{ijkl'}\otimes
 (P^a_a)^{\sigma-2})$ when $\sigma=\frac{n}{2}-1$. If $s=\sigma-1$, then
$P(g)|_{\Theta_s}$ is ``good'' if all complete contractions in $P(g)|_{\Theta_s}$
have $\delta_W+\delta_P=\frac{n}{2}-1$.\footnote{
In other words, if there are complete contractions in
$P(g)|_{\Theta_s}$ with $\delta_W+\delta_P<\frac{n}{2}-1$
 then $P(g)|_{\Theta_s}$ is ``good'' if no complete contractions
in $P(g)|_{\Theta_s}$
have $\sigma-2$ factors $P^a_a$.}

\par We will prove in the paper \cite{alexakis3} in this series the following Lemma: 
\begin{lemma}
\label{tistexnes} 
There exists a divergence $div_iT^i(g)$ so that
$$P(g)|_{\Theta_s}-div_iT^i(g)=\sum_{l\in \Theta'_s} a_l
C^l(g)+\sum_{t\in T} a_t C^t(g).$$ Here each $C^t(g)$ is in the
form (\ref{liniP}) and has fewer than $s$ factors $\nabla^{(p)}P$.
The complete contractions indexed in $\Theta'_s$ are in the form
(\ref{liniP}) with $s$ factors $\nabla^{(p)}P$ and moreover this
linear combination is good. 
\end{lemma}
Lemma \ref{tistexnes} will be proven in \cite{alexakis3}, by {\it explicitly}
 constructing the divergence $div_iT^i(g)$. (There is no recourse to 
the ``main algebraic Proposition'').   
Therefore, for the rest of this
section when $s=\sigma-1$ or $s=\sigma-2$ we will be assuming that
$P(g)|_{\Theta_s}$ is good.

\par We consider $\mu=min_{l\in L} \delta(l)$
 (recall that $L$ is the index set on the right
hand side of (\ref{avni}). We denote by
$L_\mu\subset L$ to be the set for which $l\in L_\mu$
 if and only if $\delta(l)=\mu$.
 We claim the following:

\begin{proposition}
\label{burj} Under the assumptions of Proposition \ref{killtheta} (and assuming the Lemma \ref{tistexnes},
\footnote{Recall in particular the definition of the index set
$\Theta_s$, and that we have written out
$P(g)|_{\Theta_s}=\sum_{l\in L} C^l(g)$ (modulo longer complete
contractions); recall also that if $s=\sigma-1$ or $s=\sigma-2$
then $P(g)|_{\Theta_s}$ is assumed to be good.} we claim that there
is a linear combination $T^i(g)=\Sum_{r\in R} a_r C^{r,i}(g)$,
where each $C^{r,i}(g)$ is in the form (\ref{liniP}) with length
$\sigma$, weight $-n+1$  and $\delta=\mu$, so that modulo complete
contractions of length $\ge\sigma +1$:

\begin{equation}
\label{poutses} \Sum_{l\in L_\mu} a_l C^l(g)-div_i \Sum_{r\in R}
a_r C^{r,i}(g)=\Sum_{u\in U} a_u C^u(g) +\Sum_{x\in X} a_x C^x(g),
\end{equation}
where each $C^u(g)$ is in the form (\ref{liniP}) with $s$ factors
$\nabla^{(p)}P$ and $\delta=\mu +1$.  Each $C^x(g)$ is in the form
 (\ref{liniP}) with $s-1$ factors $\nabla^{(p)}P$.
\end{proposition}

The remainder of this paper is devoted to proving the above
(subject to the ``main algebraic Proposition'' \ref{pregiade}).
 For now, we note that  Proposition \ref{pregiade} 
implies Proposition \ref{killtheta}, by iterative
repetition: After a finite number of applications of the above, we
will be left with correction terms that are of the form
$\Sum_{x\in X} a_x\dots$. This is because we are dealing with
complete contractions of a {\it fixed} weight $-n$, thus there can
be at most $\frac{n}{2}$ internal contractions in any such
complete contraction.
\newline

{\bf Proof of Proposition \ref{burj}:}
\newline

We firstly wish to understand explicitly how the terms of length $\sigma$ in $I^s_g(\psi)$
arise from $P(g)|_{\Theta_s}$. Then, we reduce Proposition \ref{burj}
to the Lemmas  \ref{corola}, \ref{cancwant}.

\par We consider $I^s_{g}(\psi)(:=\frac{d^s}{dt^s}|_{t=0}[e^{t\cdot n\psi}P(e^{t\cdot 2\psi}g)-P(g)])$. 
It follows straightorwardly from the transformation law of the Schouten tensor
that:

\begin{equation}
\label{souha} I^s_{g}(\psi)=(-1)^s\Sum_{l\in L} a_l C^l_{g}(\psi)+(Junk),
\end{equation}
where each $C^l_{g}(\psi)$ arises from $C^l(g)$ (which 
is in the form (\ref{liniP})) by replacing each
factor $\nabla^{a_1\dots a_t}\nabla^{(p)}_{t_1\dots t_p}P_{ij}$ 
by $\nabla^{a_1\dots a_t}\nabla^{(p_2)}_{t_1\dots t_pij}\psi$. 
Explicitly, it will be in the form:

\begin{equation}
\label{onlywpsi}
\begin{split}
 & contr(\nabla^{a_1\dots
a_t}\nabla^{(m_1)}W_{ijkl}\otimes\dots\otimes\nabla^{b_1\dots
b_u}\nabla^{(m_{\sigma-s})}W_{i'j'k'l'}
\\&\otimes\nabla^{v_1\dots
v_x}\nabla^{(p_1+2)}\psi\otimes\dots\otimes \nabla^{y_1\dots
y_w}\nabla^{(p_s+2)}\psi),
\end{split}
\end{equation}
and will have $\delta_W+\delta_\psi\ge \mu$ ($\delta_\psi$ here
stands for the total number of internal contractions among the
 factors $\nabla^{(k)}\psi$).
$(Junk)$ stands for a generic linear combination of terms 
with at least $\sigma+1$ factors in the form $\nabla^{(m)}R$, 
$\nabla^{(a)}\psi$.  Furthermore, 
 $\int_{M^n}I^s_{g}(\psi)dV_{g}=0$; hence we may apply the
 super divergence formula to this integral equation.
 Now, for convenience,
 we polarize the function $\psi$ and thus we will be
considering $I^s_{g}(\psi_1,\dots ,\psi_s)$.
\newline

We will now re-write $I^s_g(\psi)$ as a linear combination of complete contractions
involving curvature, rather than Weyl, tensors:
\par By decomposing the Weyl tensor as in (\ref{Weyl}) and applying 
the curvature and Bianchi identities, we re-write $I^s_{g}(\psi_1,\dots ,\psi_s)$
as a linear combination:

\begin{equation}
\label{souha2} I^s_{g}(\psi_1,\dots,\psi_s)=\Sum_{b\in B} a_b
C^b_{g}(\psi_1,\dots ,\psi_s),
\end{equation}
where each $C^b_{g}(\psi_1,\dots ,\psi_s)$ is in the form:

\begin{equation}
\label{linisymric}
\begin{split}
&contr(\nabla^{f_1\dots f_y}\nabla^{(m_1)}R_{ijkl}\otimes\dots
\otimes \nabla^{g_1\dots g_p}\nabla^{(m_s)}R_{ijkl}
\\&\otimes \nabla^{y_1\dots
y_w}\nabla^{(d_1)}Ric_{ij}\otimes\dots\otimes\nabla^{x_1\dots
x_p}\nabla^{(d_q)}Ric_{ij} \otimes R^\alpha \otimes
\\& \nabla^{a_1\dots a_{t_1}}\nabla^{(u_1)}\psi_1\otimes\dots\otimes \nabla^{c_1\dots
c_{t_s}}\nabla^{(u_s)}\psi_s),
\end{split}
\end{equation}
with the usual conventions: In each factor
$\nabla^{f_1\dots f_y}\nabla^{(m)}_{r_1\dots r_m}R_{ijkl}$,
 each of the the indices ${}^{f_1},\dots ,{}^{f_y}$
contracts against one of the indices ${}_{r_1},\dots ,{}_l$, while
no two of the indices ${}_{r_1},\dots ,{}_l$ contract between
themselves. On the other hand, for each factor $\nabla^{y_1\dots
y_t} \nabla^{(u)}_{a_1\dots a_u}\psi_h$,  each of the
upper indices ${}^{y_1},\dots ,{}^{y_t}$ contracts against one of
the indices ${}_{a_1},\dots ,{}_{a_u}$. Moreover, none of the
indices ${}_{a_1},\dots ,{}_{a_u}$ contract between themselves.
For the factors $\nabla^{x_1\dots x_p}\nabla^{(u)}_{t_1\dots
t_u}Ric_{ij}$, we impose the condition that each of the indices
${}^{x_1},\dots ,{}^{x_p}$ must contract against one of the
indices ${}_{t_1},\dots ,{}_{t_u},{}_i,{}_j$. Moreover, we impose
the restriction that none of the indices ${}_{t_1},\dots
,{}_{t_u},{}_i,{}_j$ contract between themselves (this assumption
can be made by virtue of the contracted second Bianchi identity).

\begin{definition}
\label{normalcontr} A contraction in the form (\ref{linisymric})
with all the features described above, and with the additional
requirement that each factor $\nabla^{a_1\dots
a_{t}}\nabla^{(u)}\psi_h$ has $t+u\ge 2$ (i.e.~$\psi_h$ is
differentiated at least twice) will be called normal.

\par For any complete contraction in the form (\ref{linisymric}),
 $\delta_R$ will stand for the number of
internal contractions in factors $\nabla^{(m)}R_{ijkl}$.
$\delta_{Ric}$ will stand for the number of internal contractions
in factors $\nabla^{(p)}Ric_{ij}$, {\it where we also count the
internal contraction in $Ric_{ij}={R^k}_{ikj}$}, plus $2\alpha$,
where $\alpha$ stands for the number of factors $R$ (scalar
curvature). Lastly, $\delta_\psi$ will stand for the total number
of internal contractions in the factors of the form
$\nabla^{a_1\dots a_{t}}\nabla^{(u)}\psi_h$.
\end{definition}

\par By the formula (\ref{onlywpsi}), we see
 that the sublinear combination of length $\sigma$
in $I^s_{g}(\psi_1,\dots ,\psi_s)$ consists of complete
contractions with at least two derivatives on each function
$\psi_h$.
\newline

Let us now understand more concretely how a given term in the 
form (\ref{onlywpsi}) gives rise to terms of the form (\ref{linisymric}).
We first introduce some definitions:

\begin{definition}
\label{skeletos}
  For each complete contraction $C^l_{g}(\psi)$,
$l\in L_\mu$, let us denote by $C^{l,\iota}_{g}(\psi)$ the
complete contraction (times a constant) that arises from
$C^l_{g}(\psi)$ by replacing the factors $\nabla^{(m)}W_{ijkl}$
according to the following rule: If $\nabla^{(m)}W_{ijkl}$ does
not have an internal contraction involving one of the indices
${}_i,{}_j,{}_k,{}_l$, we replace it by $\nabla^{(m)}R_{ijkl}$. If it has at least one
internal contraction involving one of the indices ${}_i,{}_j,{}_k,{}_l$,
 we replace it by $\frac{n-3}{n-2}\nabla^{(m)}R_{ijkl}$.
\end{definition}

\par Observe that by construction, if $C^l(g)$ has $\delta_W+\delta_P=b$, then
$C^{l,\iota}_{g}(\psi_1,\dots ,\psi_s)$ has
$\delta_R+\delta_{\psi}=b$, and no factors $\nabla^{(p)}Ric$ or $R$.

\par In particular, $C^{l,\iota}_{g}(\psi_1,\dots ,\psi_s)$
 will be in the form:

\begin{equation}
\label{linisym2}
\begin{split}
 & (Const)\cdot contr(\nabla^{a_1\dots
a_t}\nabla^{(m_1)}R_{ijkl}\otimes\dots\otimes\nabla^{b_1\dots
b_u}\nabla^{(m_{\sigma-s})}R_{i'j'k'l'}
\\&\otimes\nabla^{v_1\dots
v_x}\nabla^{(p_1+2)}\psi_1\otimes\dots\otimes\nabla^{y_1\dots
y_w}\nabla^{(p_s+2)}\psi_s).
\end{split}
\end{equation}

\begin{definition}
\label{proka}
Consider any $C^l_g(\psi_1,\dots,\psi_a)$ in the form (\ref{linisymric}) 
with $\sigma$ factors. 
If $C^l_{g}(\psi_1,\dots ,\psi_s)$ has $q=0$ and
$\delta =\mu$ it will be called a target. If $C^l_{g}(\psi_1,\dots
,\psi_s)$ has $q=0$ and $\delta >\mu$, it will be called a
contributor.

\par If $C^l_{g}(\psi_1,\dots ,\psi_s)$ has $q>0$ and
$\delta > \mu$ we call it 1-cumbersome. We call
$C^l_{g}(\psi_1,\dots ,\psi_s)$ 2-cumbersome if it has $q>0$ and
$\delta = \mu$ and the feature that each
 factor $\nabla^{a_1\dots a_t}\nabla^{(p)}_{r_1\dots r_p}Ric_{ij}$ has $t>0$ and
the index ${}_j$ is contracting against one of the indices
${}^{a_1},\dots ,{}^{a_t}$.

\par Finally, when we say $C^l_g(\psi_1,\dots,\psi_s)$ is  
``cumbersome'', we will mean it is either 1-cumbersome or
2-cumbersome.
\end{definition}

\par We make the convention that when $\Sum_{j\in J} a_j C^j_{g}(\psi_1,\dots
, \psi_s)$, \\$\Sum_{f\in F} a_f C^f_{g}(\psi_1,\dots , \psi_s)$
appear {\it on the right hand sides} of  equations below, they
will stand for {\it generic} linear combinations of contributors
and cumbersome complete contractions, respectively.

 Then using the
 decomposition of the Weyl curvature (\ref{Weyl}), we explicitly write each
 $C^l_{g}(\psi_1,\dots, \psi_s)$ as a linear combination 
of terms in the above forms:

\par For each $l\in L_\mu$, it follows that:

\begin{equation}
\label{lj}
\begin{split}
&C^l_{g}(\psi_1,\dots ,\psi_s)= C^{l,\iota}_{g}(\psi_1,\dots
,\psi_s)+ \Sum_{j\in J} a_j C^j_{g}(\psi_1,\dots ,\psi_s) +
\\& \Sum_{f\in F} a_f C^f_{g}(\psi_1,\dots ,\psi_s),
\end{split}
\end{equation}
while for each $l\in L\setminus L_\mu$:

\begin{equation}
\label{l'j} C^l_{g}(\psi_1,\dots ,\psi_s)= \Sum_{j\in J} a_j
C^j_{g}(\psi_1,\dots ,\psi_s) + \Sum_{f\in F} a_f
C^f_{g}(\psi_1,\dots ,\psi_s),
\end{equation}
where  each $C^f_{g}(\psi_1,\dots
,\psi_s)$ has $\delta_R +\delta_\psi +\delta_{Ric}\ge \mu+1$ (and
hence is 1-cumbersome). This follows since
$C^l_{g}(\psi_1,\dots,\psi_s)$, $l\in L\setminus L_\mu$  has
$\delta_W+\delta_P\ge \mu+1$. 

{\it Remark:} We observe that for each complete contration 
$C^f_g(\psi_1,\dots,\psi_s)$ in the RHSs of (\ref{lj}), (\ref{l'j}) with 
$\alpha>0$ factors $R$ (of the scalar curvature) will respectively satisfy 
$\delta\ge \mu+2\alpha$, $\delta\ge\mu+1+2\alpha$. This is because a factor $R$ 
in the RHS can only arise from an (undifferentiated)  factor $W_{ijkl}$
in the LHS of (\ref{lj}), (\ref{l'j}); thus a factor with no 
internal contractions in the LHS gives rise to a factor $R={R^{ab}}_{ab}$ 
with two internal contractions. (This remark will be useful in \cite{alexakis3}).
\newline

\par In view of the form (\ref{onlywpsi}) where each complete 
contraction has $\delta_W+\delta_P\ge \mu$, we derive that:

 \begin{equation}
 \label{polis}
\begin{split}
&I^s_{g}(\psi_1,\dots ,\psi_s)=\Sum_{l\in L_\mu} a_l
C^{l,\iota}_{g}(\psi_1,\dots ,\psi_s)+ \Sum_{j\in J} a_j
C^j_{g}(\psi_1,\dots ,\psi_s) 
\\&+ \Sum_{f\in F} a_f
C^f_{g}(\psi_1,\dots ,\psi_s)+(Junk);
\end{split}
 \end{equation}
$L_\mu$ here is the same index set as in Proposition \ref{burj}. The linear combinations
 $\sum_{j\in J}\dots,\sum_{f\in F}\dots$ are generic linear combinations
  of contributors and cumbersome complete contractions (see definition \ref{proka}).
$(Junk)$ stands for a generic linear combination of terms 
with at least $\sigma+1$ factors in total.

\par Now, for the next Lemma, we will let $Z_g(\psi_1,\dots ,\psi_s)$
stand for {\it any} linear combination in the form above, where
$\sum_{l\in L_\mu}\dots$ is the {\it same} linear combination as in
$I^s_{g}(\psi_1,\dots ,\psi_s)$, while $\sum_{f\in
F}\dots,\sum_{j\in J}\dots$, $(Junk)$  are allowed to be {\it generic} linear
combinations of the forms described above.  In these generic
linear combinations $Z_g(\psi_1,\dots ,\psi_s)$ we will still be
assuming that $Z_g(\psi_1,\dots ,\psi_s)$ is symmetric in the
functions $\psi_1,\dots,\psi_s$.

\par We partition the index set $F$ into subsets: We let $f\in F^{q,z}$
if and only if \\$C^f_{g}(\psi_1,\dots ,\psi_s)$ has $q$ factors of
the form $\nabla^{a_1\dots a_t}\nabla^{(p)}Ric$ or $R$ and also
has $\delta_R+\delta_{Ric}+\delta_\psi=z$. We also define
$F^q=\bigcup_{z\ge \mu}F^{q,z}$. One last language convention
before stating our claims: We will say that the index set $F^{q}$
(or more generally $F$) is {\it bad} if there are complete
contractions $C^f_g(\psi_1,\dots ,\psi_s)$, $f\in F^q$ with at
least $\sigma-2$ factors in the form $\Delta\psi_h$ {\it or} $R$
(scalar curvature).
\newline

{\it The main Claims:}

\begin{lemma}
\label{corola} Consider any  $Z_{g}(\psi_1,\dots ,\psi_s)$,
written out as a linear combination in the form (\ref{polis}).
Assume that $\int_{M^n}Z_g(\psi_1,\dots,\psi_\sigma)dV_g=0$ for every $(M^n,g)$ and every
 function $\psi_1=\dots=\psi_\sigma=\psi\in C^\infty (M^n)$.\footnote{This is just a
re-statement of the fact that $Z_g(\psi_1,\dots,\psi_\sigma)$ 
is symmetric in the functions $\psi_1,\dots,\psi_\sigma$.}
 Assume also that for a given $q_1>0$, $F^q=\emptyset$ for every $q>q_1$.
 Moreover, we assume that for
 a given  $z_1\ge \mu$
$F^{q_1,z}=\emptyset$ for every $z<z_1$. We make different claims for the two cases
$z_1>\mu$ and $z_1=\mu$.

\par If $z_1> \mu$ and $F^{q_1}$ is not bad,\footnote{See the language convention above.}
 we claim that there
 is a linear combination of vector fields,
$\Sum_{h\in H^{q_1,z_1}}a_h C^{h,i}_{g}(\psi_1,\dots , \psi_s)$
where each $C^{h,i}_{g}(\psi_1,\dots ,\psi_s)$ is in the form
(\ref{linisymric}) with length $\sigma$, $q+\alpha =q_1$,
$\delta_R+\delta_{Ric}+\delta_\psi= z_1$ and with one free index,
so that:

\begin{equation}
\label{cancq1}
\begin{split}
& \Sum_{f\in F^{q_1,z_1}} a_f C^f_{g}(\psi_1,\dots ,\psi_s) -
div_i \Sum_{h\in H^{q_1,z_1}} a_h C^{h,i}_{g}(\psi_1,\dots
,\psi_s)=
\\& \Sum_{s\in S} a_s C^s_{g}(\psi_1,\dots ,\psi_s) +
\Sum_{t\in T}
 a_t C^t_{g}(\psi_1,\dots ,\psi_s),
\end{split}
\end{equation}
 where each $C^s_{g}(\psi_1,\dots ,\psi_s)$ is in the form
(\ref{linisymric}), has length $\sigma$ and is not bad, and has
$q+\alpha=q_1-1$, $\delta_R+\delta_{Ric}+\alpha =z_1$. On the
other
 hand, each $C^t_{g}(\psi_1,\dots ,\psi_s)$ is
 of the form (\ref{linisymric}) with length $\sigma$ and
$q+\alpha =q_1$ factors $\nabla^{(p)}Ric$ or $R$,
$\delta_R+\delta_{Ric}+\delta_\psi=z_1+1$. The above holds modulo
complete contractions of length $>\sigma$.

\par In the case where $z_1=\mu$ and $F^{q_1}$ is not bad, we claim that there is a
linear combination of vector fields $\Sum_{h\in H^{q_1,z_1}} a_h
C^{h,i}_{g}(\psi_1,\dots ,\psi_s)$, where each
\\$C^{h,i}_{g}(\psi_1,\dots ,\psi_s)$ is in the form
(\ref{linisymric}) with length $\sigma$, $q+\alpha =q_1$,
$\delta_R+\delta_{Ric}+\delta_\psi= \mu$ and with one free index,
so that:

\begin{equation}
\label{cancq1'}
\begin{split}
& \Sum_{f\in F_{q_1,\delta_1}} a_f C^f_{g}(\psi_1,\dots ,\psi_s) -
div_i \Sum_{h\in H_{q_1,z_1}} a_h C^{h,i}_{g}(\psi_1,\dots
,\psi_s)=
\\& \Sum_{t\in T} a_t C^t_{g}(\psi_1,\dots ,\psi_s)
\end{split}
\end{equation}
 where each $C^t_{g}(\psi_1,\dots ,\psi_s)$ is
 in the form (\ref{linisymric}) (not bad) with length $\sigma$ and
$q+\alpha =q_1$ factors $\nabla^{(p)}Ric$ or $R$, $\delta=\mu+1$.
The above holds modulo complete contractions of length $>\sigma$.
\newline

\par Claim 2: Consider $I^s_g(\psi_1,\dots,\psi_s)$,
in the form (\ref{polis}), and suppose $F$ is bad.
We claim that there is a
linear combination of vector fields, \\$\sum_{h\in H^{q_1}} a_h
C^{h,i}_{g}(\psi_1,\dots ,\psi_s)$ where each
$C^{h,i}_{g}(\psi_1,\dots ,\psi_s)$ is in the form
(\ref{linisymric}) with length $\sigma$, $q+\alpha =q_1$,
 and with one free index, so that:

\begin{equation}
\label{antidiaploki}
\begin{split}
& \Sum_{f\in F_{q_1}} a_f C^f_{g}(\psi_1,\dots ,\psi_s) - div_i
\Sum_{h\in H_{q_1}} a_h C^{h,i}_{g}(\psi_1,\dots ,\psi_s)=
\\& \Sum_{x\in X} a_x C^x_{g}(\psi_1,\dots ,\psi_s),
\end{split}
\end{equation}
where each of the complete contractions $C^x_g(\psi_1,\dots
,\psi_s)$ is of the form (\ref{linisymric}) with length $\sigma$
and $q+\alpha \le q_1$ and $\delta\ge\mu+1$, and is not bad. The above holds
modulo complete contractions of length $>\sigma$.
\end{lemma}

{\it Note:} Claim 2 will be proven in \cite{alexakis3}. 
\newline

\par Observe that the  Lemma \ref{corola} implies that there
 is a linear combination of vector fields
$\Sum_{h\in H} a_h C^{h,i}_{g}(\psi_1,\dots ,\psi_s)$, where each
$C^{h,i}_{g}(\psi_1,\dots ,\psi_s)$ is a partial contraction of
length $\sigma$ in the form (\ref{linisymric})
 and with one free index, so that:

\begin{equation}
\label{killspec} \Sum_{f\in F} a_f C^f_{g}(\psi_1,\dots ,\psi_s)
-div_i \Sum_{h\in H} a_h C^{h,i}_{g}(\psi_1,\dots ,\psi_s)=
\Sum_{j\in J} a_j C^j_{g}(\psi_1,\dots ,\psi_s).
\end{equation}
Here the first sublinear combination is {\it not} generic, but
stands for the sublinear combination in (\ref{polis}). The above
holds modulo complete contractions of length $>\sigma$. Therefore, assuming we can prove
Lemma \ref{corola} we can then apply it to the integral equation
$\int_{M^n}I^s_g(\psi_1,\dots,\psi_s)dV_g=0$ (recall that $I^s_g(\psi_1,\dots,\psi_s)$
is in the form (\ref{polis})), to derive a new integral equation:

\begin{equation}
 \label{new.int}
\int_{M^n} \Sum_{l\in L_\mu} a_l
C^{l,\iota}_{g}(\psi_1,\dots , \psi_s)+\Sum_{j\in J} a_j
C^j_{g}(\psi_1,\dots ,\psi_s)
\\+\sum_{\zeta\in Z} a_\zeta
C^\zeta_g(\psi_1,\dots ,\psi_s) dV_{g}=0,
\end{equation}
which holds for every $(M^n,g)$ and every $\psi_1=\dots=\psi_\sigma=\psi\in C^\infty (M^n)$
(recall that our complete contractions are assumed to be symmetric in
the functions $\psi_1,\dots,\psi_s$), where the 
complete contractions $C^\zeta_g$ have length $\ge\sigma+1$
and the complete contractions $C^l_g,C^j_g$ are as described
below equation (\ref{polis}).\footnote{In particular, the
linear combination $\sum_{l\in L_\mu} a_l C^{l,\iota}_{g}(\psi_1,\dots , \psi_s)$ is {\it the same}
linear combination that appears in Proposition \ref{burj}, while the
linear combination $\sum_{j\in J}\dots$ is a generic linear combination of
complete contractions as explained below equation (\ref{polis}).}
\newline

\par Our next Lemma will then apply to the {\it new} integral equation 
(\ref{new.int}). In order to state it, we will need one extra 
definition:

\begin{definition}
\label{pongeconnes} For each complete contraction $C_{g}(\phi)$ or
vector field $C^i_{g}(\phi)$ in the form (\ref{linisymric}), with
no factors $\nabla^{(p)}Ric$ or $R$ and with $\delta=\mu$ (in
other words there are $\mu$ internal contractions and all of them
involve a derivative index), we denote by
$C^{|i_1\dots i_\mu}_{g}(\phi)\nabla_{i_1}\upsilon\dots\nabla_{i_\mu}\upsilon$,
$C^{i|i_1\dots i_\mu}_{g}(\phi)\nabla_{i_1}\upsilon\dots\nabla_{i_\mu}\upsilon$,
 the complete contraction or
vector field that arises from it by replacing each internal
contraction $(\nabla^a,{}_a)$ by an expression
$(\nabla^a\upsilon,{}_a)$.\footnote{We thus obtain complete
contractions and vector fields of length $\sigma +\mu$.}
\end{definition}

\begin{lemma}
\label{cancwant}
\par Assume an equation:
\begin{equation}
\label{bigdaddy} \begin{split} &\int_{M^n} \Sum_{l\in L_\mu} a_l
C^{l,\iota}_{g}(\psi_1,\dots , \psi_s)+\Sum_{j\in J} a_j
C^j_{g}(\psi_1,\dots ,\psi_s) \\&+\sum_{\zeta\in Z} a_\zeta
C^\zeta_g(\psi_1,\dots ,\psi_s) dV_{g}=0,
\end{split}
\end{equation}
which holds for every compact $(M^n,g)$ and every $\psi_1,\dots
,\psi_s\in C^\infty (M^n)$, and where each $C^\zeta_g$ has length
$\ge\sigma+1$. We then claim that there is a linear combination of
normalized vector fields $\Sum_{d\in D} a_d C^{d,i}_{g}
(\psi_1,\dots ,\psi_s)$, where each $C^{d,i}_{g} (\psi_1,\dots ,
\psi_s)$ is in the form (\ref{linisymric}) with no factors
$\nabla^{(p)}Ric$ or
 $R$ and with $\delta=\mu$, so that:

\begin{equation}
\label{soazig} 
\begin{split}
 &\Sum_{l\in L_\mu} a_l 
C^{l,\iota|i_1\dots i_\mu}_{g}(\psi_1,\dots , \psi_s)
\nabla_{i_1}\upsilon\dots\nabla_{i_\mu}\upsilon
 \\&-\Sum_{d\in D} a_d Xdiv_i
C^{d,i|i_1\dots i_\mu}_{g}(\psi_1,\dots , \psi_s)\nabla_{i_1}\upsilon\dots\nabla_{i_\mu}\upsilon=0,
\end{split}
\end{equation}
modulo complete contractions of length $\ge\sigma +\mu+1$.
\end{lemma}

\par Let us check how  Proposition
\ref{burj} follows from Lemma \ref{cancwant}.

\par Just observe that Lemma \ref{paparidis} and  (\ref{soazig})
imply that 
 the vector field $T^i:=Weylify[\Sum_{d\in D} a_d
\Pi\rho^\mu_\upsilon C^{d,i}_{g}(\psi_1,\dots ,\psi_s)]$ 
fulfils the requirements 
of Proposition \ref{burj}.
\newline
\par Thus, if we can show Lemmas \ref{corola}, \ref{cancwant},
Proposition \ref{killtheta} will follow.

\subsection{The main algebraic Proposition.}
\label{fundprop}

The Proposition  that we state in
this section will imply Lemmas \ref{corola} and \ref{cancwant}.
\newline

\par In order to state and prove the main algebraic proposition we will
 need some more terminology. We will be considering tensor fields
$C^{i_1\dots i_\alpha}_{g}(\Omega_1,\dots , \Omega_p)$ 
of length $\sigma$ (with no internal contractions) in
the  form:

\begin{equation}
\label{preform1}
\begin{split}
&pcontr(\nabla^{(m_1)}R_{ijkl}\otimes\dots\otimes\nabla^{(m_s)}R_{ijkl}
\otimes \nabla^{(b_1)}\Omega_1\otimes\dots\otimes \nabla^{(b_p)}\Omega_p);
\end{split}
\end{equation}
here $\sigma=s+p$ and ${}_{i_1},\dots ,{}_{i_\alpha}$ are the free indices.
Such a complete contraction will be called {\it acceptable} if each $b_i\ge 2$.
Recall the operation $Xdiv$ from Definition \ref{Xdiv}. 

\begin{proposition}
\label{pregiade}

Consider two linear combinations of acceptable tensor fields in
the form (\ref{preform1}),

$$\Sum_{l\in L_1} a_l
C^{l,i_1\dots i_{\alpha}}_{g} (\Omega_1,\dots
,\Omega_p), \Sum_{l\in L_2} a_l
C^{l,i_1\dots i_{\beta_l}}_{g} (\Omega_1,\dots
,\Omega_p),$$
 where each $C^{l,i_1\dots i_{\alpha}}_{g}$ above has length $\sigma\ge 3$ and a given
number $\sigma_1= \sigma-p$ of factors in the form $\nabla^{(m)}R_{ijkl}$.
 Assume  that for each
$l\in L_2$,  $\beta_l\ge \alpha+1$.  Assume that modulo complete contractions
of length $\ge\sigma+1$:

\begin{equation}
\label{hypothese2}
\begin{split}
&\Sum_{l\in L_1} a_l Xdiv_{i_1}\dots Xdiv_{i_{\alpha}}
C^{l,i_1\dots i_{\alpha}}_{g} (\Omega_1,\dots
,\Omega_p)+
\\&\Sum_{l\in L_2} a_l Xdiv_{i_1}\dots Xdiv_{i_{\beta_l}}
C^{l,i_1\dots i_{\beta_l}}_{g} (\Omega_1,\dots
,\Omega_p)=0.
\end{split}
\end{equation}

We claim that there is a linear
combination of acceptable $(\alpha+1)$-tensor fields
in the form (\ref{preform1}), $\sum_{r\in R} a_r C^{r,i_1\dots i_{\alpha+1}}_g(\Omega_1,\dots,\Omega_p)$,
with length $\sigma$ so that:\footnote{Recall that 
given a $\beta$-tensor field $T^{i_1,\dots i_\alpha\dots i_\beta}$,
$T^{(i_1\dots i_\alpha)\dots i_\beta}$ stands 
for a new tensor field that arises from
$T^{i_1,\dots i_\alpha\dots i_\beta}$ by symmetrizing over the indices ${}^{i_1},\dots,{}^{i_\alpha}$.}

\begin{equation}
\label{concln}
\begin{split}
&\Sum_{l\in L_1} a_l C^{l,(i_1\dots i_{\alpha})}_{g} (\Omega_1,\dots
,\Omega_p)=\sum_{r\in R} a_r Xdiv_{i_{\alpha+1}}C^{r,(i_1\dots i_\alpha)i_{\alpha+1}}_g
(\Omega_1,\dots,\Omega_p),
\end{split}
\end{equation}
modulo terms of length $\ge\sigma+1$.
\end{proposition}

{\it Note:} Observe that the conclusion (\ref{concln}) of 
this Proposition is equivalent to the equation:

\begin{equation}
\label{concln'}
\begin{split}
&\Sum_{l\in L_1} a_l C^{l,i_1\dots i_{\alpha}}_{g} (\Omega_1,\dots
,\Omega_p)\nabla_{i_1}\upsilon\dots\nabla_{i_\alpha}\upsilon=
\\&\sum_{r\in R} a_r Xdiv_{i_{\alpha+1}}
[C^{r,i_1\dots i_\alpha i_{\alpha+1}}_g(\Omega_1,\dots,\Omega_p)
\nabla_{i_1}\upsilon\dots\nabla_{i_\alpha}\upsilon],
\end{split}
\end{equation}
which holds modulo complete contractions of 
length $\ge\sigma+\alpha+1$. (Recall that 
$Xdiv_{i_{\alpha+1}}[\dots]$ in the RHS of the above 
stands for the sublinear combination of terms in $div_{i_{\alpha+1}}[\dots]$ 
where the derivative $\nabla^{i_{\alpha+1}}$ is not allowed to 
hit the factor to which the free index ${}_{i_{\alpha+1}}$ {\it nor}
any of the factors $\nabla\upsilon$). 
\newline
In the next subsection we show how the main algebraic proposition \ref{pregiade}
implies Lemmas \ref{corola}, \ref{cancwant}, and hence also Proposition 
\ref{killtheta}.

\subsection{Lemmas \ref{corola} and \ref{cancwant} follow from
 Proposition \ref{pregiade}.}
\label{peresprez}

\par We first check that Lemma \ref{cancwant} indeed follows from
Proposition \ref{pregiade}.
\newline

\par Our starting point will be to apply the super divergence formula to the integral equation
(\ref{bigdaddy}). 

\begin{definition}
\label{axrhst} For each $l\in L_\mu$ and each $j\in J$, we denote
by $C^{l,\iota|i_1\dots i_\mu}_{g}(\psi_1,\dots ,\psi_s)$,
$C^{j|i_1\dots i_{m_j}}_{g}(\psi_1,\dots ,\psi_s)$ the tensor
fields that arise from $C^{l,\iota}_{g}(\psi_1,\dots ,\psi_s)$,
\\$C^j_{g}(\psi_1,\dots ,\psi_s)$, respectively, by making all the
internal contractions into free indices (recall the 
definition \ref{eraseindex}).
\end{definition}

The super divergence formula applied to (\ref{bigdaddy}) gives the local equation:

\begin{equation}
\label{bigdaddy2}
\begin{split}
&(-1)^\mu\Sum_{l\in L_\mu} a_l Xdiv_{i_1}\dots
Xdiv_{i_\mu}C^{l,\iota|i_1\dots i_\mu}_{g}(\psi_1,\dots ,\psi_s)+
\\&\Sum_{j\in J} a_j (-1)^{m_j-1}
Xdiv_{i_1}\dots Xdiv_{i_{m_j}} C^{j|i_1\dots
i_{m_j}}_{g}(\psi_1,\dots ,\psi_s)=0,
\end{split}
\end{equation}
which holds modulo complete contractions of length $\ge\sigma +1$.

\par Clearly, each of the complete contractions $C^{l,\iota}$,
$C^j$ has factors $\nabla^{(b)}\psi_h$ with $b\ge 2$. Therefore,
each of the tensor fields in (\ref{bigdaddy2}) has factors
$\nabla^{(c)}\psi_h$ with $c\ge 1$, and moreover the factors
$\nabla\psi_h$ can only arise from factors $\Delta\psi_h$ by
replacing $\Delta\psi_h$ by $\nabla_a\psi_h$ (${}_a$ is a free
index).

\par For each $C^{l,\iota}$, $C^j$ appearing in (\ref{cancwant})
let $|\Delta|(l),|\Delta|(j)$  stand for the number of factors
$\Delta\psi_h$.  We define $|\Delta|_{Max}$ to stand for
$max_{f\in L_\mu, f\in J}|\Delta|(f)$.
\newline

\par We observe that if $|\Delta|_{Max}=0$, then 
Lemma \ref{cancwant} follows by just applying Proposition \ref{pregiade} to
(\ref{bigdaddy2}). In the case $|\Delta|_{Max}>0$ we {\it cannot} directly 
apply Proposition \ref{pregiade} to (\ref{bigdaddy2}) due to the presence of factors $\nabla\psi_h$
among certain tensor fields in (\ref{bigdaddy2}).
 We will treat the case $|\Delta|_{Max}>0$
further down.
\newline

{\it Lemma \ref{corola} follows from Proposition \ref{pregiade}
(general discussion):} (Refer to the  notation of Lemma \ref{labrome}). We
apply Lemma \ref{labrome} to the equation $\int_{M^n}Z_g(\psi_1,\dots,\psi_\sigma)dV_g$
(see the hypothesis of Lemma \ref{corola}) and
deduce that modulo complete contractions of length $\ge\sigma +1$:

\begin{equation}
\label{rhodesia}
\begin{split}
&\Sum_{f\in F_{q_1,z_1}} a_f Xdiv_{i_1}\dots Xdiv_{i_{z_1}}
C^{f,i_1\dots i_{z_1}}_{g}(\psi_1,\dots ,\psi_s,\Omega^{q_1})+
\\&\Sum_{z>z_1}\Sum_{f\in F_{q_1,z}}a_f
Xdiv_{i_1} \dots Xdiv_{i_z} C^{f,i_1\dots i_z}_{g}(\psi_1,\dots
,\psi_s,\Omega^{q_1})=0.
\end{split}
\end{equation}

\par We now define $|\Delta|(f)$ 
for each complete contraction
 $C^f$ in Lemma \ref{corola} to stand for the number of factors
$\Delta\psi_h$ {\it or} $\Delta\Omega$. (Observe that by
construction a factor $\Delta\Omega$ can only arise in
$C^{f,i_1\dots i_{z_1}}_{g}(\psi_1,\dots ,\psi_s,\Omega^{q_1})$ by
replacing some factor $R$ in $C^{f,i_1\dots
i_{z_1}}_{g}(\psi_1,\dots ,\psi_s)$ by $-2\Delta\Omega$).  We
define $|\Delta|_{Max}$ to stand for $max_{f\in F^{q_1,z}, z\ge
z_1}|\Delta|(f)$. We write $|\Delta|_{Max}=M$, for short.

\par  Lemma \ref{corola} in the case where $|\Delta|_{Max}=0$
can be shown by  applying Proposition \ref{pregiade} and
the operation $Riccify$ to (\ref{rhodesia}). The details of this
will be provided below, in the cases where $M:=|\Delta|_{Max}>0$. That
proof, if we set $M=0$ applies to show how Lemma \ref{corola}
follows from (\ref{rhodesia}), in the both the case $z_1=\mu$ and 
$z_1>\mu$.
\newline

{\bf Proof of Lemmas \ref{corola} and \ref{cancwant}:} We now consider equations (\ref{bigdaddy2})
and (\ref{rhodesia}) where $|\Delta|_{Max}>0$. Our strategy will
then be to reduce ourselves to the case where $|\Delta|_{Max}=0$ 
by a downward induction on $|\Delta|_{Max}$ (see below).
In this general situation, we will not show Lemmas \ref{corola}
and \ref{cancwant} all in one piece, but rather we will
distinguish cases. We distinguish three cases: Either
$|\Delta|_{Max}=\sigma -1$ or it is $\sigma -2$ or it is
$\le\sigma -3$. Here we consider only the case
$|\Delta|_{Max}\le\sigma-3$. The  cases 
$|\Delta|_{Max}=\sigma-1,|\Delta|_{Max}=\sigma-2$ will be
treated in \cite{alexakis3}. (For reference purposes, we codify the 
claim of Lemmas \ref{corola}, \ref{cancwant}
when $|\Delta|_{Max}=\sigma-1$, $|\Delta|_{Max}=\sigma-2$ in the end of this subsection.)
\newline

{\it Proof of Lemmas \ref{corola} and \ref{cancwant} in the case
$M=|\Delta|_{Max}\le\sigma -3$.}
\newline

{\it Outline:} We will claim the equations (\ref{vasko}),
(\ref{vasko2}), (\ref{timgow}), (\ref{timgow2}), (\ref{timgow3})
below, and will show how Lemmas \ref{corola} and \ref{cancwant}
will follow from these four equations. We then prove these four
equations (using Proposition \ref{pregiade}).
\newline

{\it Lemma \ref{cancwant}:} In (\ref{bigdaddy2}), we let
$L_\mu^K$, $K=1,\dots ,M$ to  the index set of the complete
contractions $C^{l,\iota}_{g}$ with $K$ factors $\Delta\psi_h$.
Accordingly, we let $J^K$, $K=1\dots ,M$ be the index set of
complete contractions $C^j_{g}$ with $K$ factors $\Delta\psi_h$.

Consider (\ref{bigdaddy2}). We claim  that there
exists a linear combination of vector fields, $\Sum_{v\in V^{M}} a_v
C^{v,i}_{g}(\psi_1,\dots ,\psi_s)$, with each $C^{v,i}_{g}$ in
the form (\ref{linisymric}) with length $\sigma$, $\delta =\mu$
and with no factors $\nabla^{(p)}Ric$ or $R$, but with $M$ factors
$\Delta\psi_h$, so that:\footnote{Note: We will be writing
$\Pi\rho^\mu[\dots]$ instead of $\Pi\rho^\mu_\upsilon[\dots]$ to
avoid confusion}

\begin{equation}
\label{vasko}
\begin{split} 
 &\Sum_{l\in L_\mu^{M}} a_l C^{l,\iota|i_1\dots i_\mu}_{g}(\psi_1,\dots ,\psi_s)
\nabla_{i_1}\upsilon\dots\nabla_{i_\mu}\upsilon=[Xdiv_i \Sum_{v\in V^{M}}
a_v C^{v,i|i_1\dots i_\mu}_{g}(\psi_1,\dots ,\psi_s)
\\&+\Sum_{r\in R^{M-1}} a_r C^{r|i_1\dots i_\mu}_{g}(\psi_1,\dots ,\psi_s)]\nabla_{i_1}\upsilon\dots\nabla_{i_\mu}\upsilon,
\end{split}
\end{equation}
where each 
$C^{r|i_1\dots i_\mu}_{g}(\psi_1,\dots ,\psi_s)$ 
on the RHS is a partial contraction in the form (\ref{linisymric}) 
(with $\mu$ free indices) 
with no
factors $\nabla^{(p)}Ric,R$, with $\delta=\mu$ but $M-1$
 factors $\Delta\psi_h$.

\par If we can prove the above, then
we will be reduced to proving Lemma \ref{cancwant} under the extra
assumption that for every $C^{l,\iota}_g$ in (\ref{bigdaddy}) will
have at most $M-1$ factors $\Delta\psi_h$.

\par  In this setting,  we define
$\delta_{min}(M)$ to stand for the minimum number of internal
contractions among the complete contractions $C^j, j\in J_M$ in (\ref{bigdaddy}). By
definition,  $\delta_{min}(M)\ge \mu+1$. We then claim
that there exists a linear combination of vector fields, $\Sum_{h\in
H^{M}_{\delta_{min}(M)}} a_h C^{h,i}_{g}(\psi_1,\dots ,\psi_s)$,
where each $C^{h,i}_{g}$ is in the form (\ref{linisymric}) with
length $\sigma$, $\delta=\delta_{min}(M)$ and with no factors
$\nabla^{(p)}Ric$ or $R$ but with $M$ factors $\Delta\psi_h$, so
that:

\begin{equation}
\label{vasko2} \begin{split} & \Sum_{j\in J^{M}_{\delta_{min}(M)}}
a_j C^{j|i_1\dots i_{\delta_{min}(M)}}_{g}(\psi_1,\dots
,\psi_s)\nabla_{i_1}\upsilon\dots \nabla_{i_{\delta_{min}(M)}}\upsilon
\\&-Xdiv_i \Sum_{h\in H^{M}_{\delta_{min}(M)}} a_h
C^{h,i||i_1\dots i_{\delta_{min}(M)}}_{g}(\psi_1,\dots ,\psi_s)
\nabla_{i_1}\upsilon\dots \nabla_{i_{\delta_{min}(M)}}\upsilon
\\&=\Sum_{r\in R^{M}} a_r 
C^{r|i_1\dots i_{\delta_{min}(M)}}_{g}(\psi_1,\dots ,\psi_s)\nabla_{i_1}\upsilon\dots \nabla_{i_{\delta_{min}(M)}}\upsilon,
\end{split}
\end{equation}
where each $C^{r|i_1\dots i_{\delta_{min}(M)}}_{g}(\psi_1,\dots ,\psi_s)$
 is a partial
contraction in the form (\ref{linisymric}) (with $\mu$ free indices) and with
$\delta_R=\delta_{min}(M)$ and $M-1$ factors $\Delta\psi_h$.

\par Observe that  (\ref{vasko}), (\ref{vasko2}),
imply Lemma \ref{cancwant}:
Iteratively applying them we reduce ourselves to proving Lemma
\ref{cancwant} under the additional assumption that each
$C^{l,\iota}$ has no factors $\Delta\psi_h$, and also each $C^j$
has no factors $\Delta\psi_h$. In that case we have already shown
how Lemma \ref{cancwant} directly follows from Proposition \ref{pregiade}. 
\newline

{\it Lemma \ref{corola}:} We make  analogous claims regarding Lemma
\ref{corola}. Consider (\ref{rhodesia}). We denote by $F^{q_1,z}_K$,
$K=1,\dots ,M$ the index set of complete contractions with K
factors $\Delta\psi_h$ or $\Delta\Omega$.
We initially consider the sublinear combination
indexed in $F^{q_1, z_1}_M$.  We then make two different claims,
for the two cases $z_1=\mu$ and $z_1>\mu$. If $z_1>\mu$, then
for some complete contraction $C^f_g$, $f\in F^{q_1,z_1}_M$, we may have factors $R$
(of the scalar curvature); if $z_1=\mu$ there can be no such factors (by definition). We
further subdivide $F^{q_1,z_1}_M$ into subsets
$F^{q_1,z_1}_{M,\alpha}$, $\alpha =0,\dots ,M$, where $f\in
F^{q_1,z_1}_{M,\alpha}$ if and only if $C^f$ has $\alpha$ factors
$R$ (and hence $M-\alpha$ factors $\Delta\psi$).

 We claim that for each of the index sets
$F^{q_1,z_1}_{M,\alpha}$ there is a linear combination of vector
fields, $\Sum_{r\in R^{q_1,z_1}_{M,\alpha}} a_r
C^{r,i}_{g}(\psi_1,\dots ,\psi_s)$ where each $C^{r,i}$ is in the
form (\ref{linisymric}) and has $q_1-\alpha$ factors
$\nabla^{(p)}Ric$, $\delta= z_1$ and $\alpha$ factors $R$ and
$M-\alpha$ factors $\Delta\psi_h$, so that modulo complete
contractions of length $\ge\sigma +1$:

\begin{equation}
\label{timgow}
\begin{split}
 &\Sum_{f\in F^{q_1, z_1}_{M,\alpha}} a_f
C^f_{g}(\psi_1,\dots ,\psi_s) -div_i \Sum_{r\in
R^{q_1,z_1}_{M,\alpha}} a_r C^{r,i}_{g}(\psi_1,\dots ,\psi_s)=
\\&\Sum_{d\in D_1} a_d C^d_{g}(\psi_1,\dots ,\psi_s)+
\Sum_{d\in D_2} a_d C^d_{g}(\psi_1,\dots ,\psi_s)+ \Sum_{d\in D_3}
a_d C^d_{g}(\psi_1,\dots ,\psi_s),
\end{split}
\end{equation}
where each $C^d$, $d\in D^1$ has $q_1$ factors $\nabla^{(p)}Ric$
and $M$ factors $\Delta\psi_h$ and $\delta=z_1+1$. Each $C^d, d\in
D_2$ has $q_1-1$ factors $\nabla^{(p)}Ric$ and $M$ factors
$\Delta\psi_1$ and $\delta=z_1$. Finally, each $C^d, d\in D_3$ has
$q_1$ factors $\nabla^{(p)}Ric$ and $M-1$ factors $\Delta\psi_h$
and $\delta=z_1$.
\newline

\par In the case where $z_1=\mu$,  we have noted that
no $C^f$, $f\in F^{q_1,\mu}_M$ has a factor $R$. We
then claim that there is a linear combination of vector fields,
$\Sum_{t\in T} a_t C^{t,i}_{g}(\psi_1,\dots ,\psi_s)$
 in the form (\ref{linisymric}) with $q_1$ factors
$\nabla^{(p)}Ric$
 and with $\delta =z_1=\mu$ and with $M$ factors
$\Delta\psi_h$, so
 that modulo complete contractions of length $\ge\sigma +1$:

\begin{equation}
\label{timgow2}
\begin{split}
 &\Sum_{f\in F^{q_1, z_1}_M} a_f
C^f_{g}(\psi_1,\dots ,\psi_s) -div_i \Sum_{t\in T} a_t
C^{t,i}_{g}(\psi_1,\dots ,\psi_s)=
\\&\Sum_{d\in D_1} a_d C^d_{g}(\psi_1,\dots ,\psi_s)+
\Sum_{d\in D_3} a_d C^d_{g}(\psi_1,\dots ,\psi_s),
\end{split}
\end{equation}
where the complete contractions on the right hand side are as
in the notation under (\ref{timgow}).
\newline

\par Assuming (for a moment) (\ref{timgow}), (\ref{timgow2}), we are reduced  to proving Lemma \ref{corola}
under the additional assumption that $F_M^{q_1,z_1}=\emptyset$. 
In that setting, 
we define $z_{min}(M)$ to stand for the minimum $z$ for
which $F^{q_1,z}_M\ne \emptyset$. By our hypothesis,
$z_{min}(M)>z_1= \mu$. On the other hand, 
some contractions $C^f_{g}$, $f\in F^{q_1,z_{min}(M)}_M$,  might
have factors $R$ (of the scalar curvature). We further
subdivide $F^{q_1,z_{min}(M)}_M$ into subsets
$F^{q_1,z}_{M,\alpha}$, $\alpha =0,\dots ,M$, where $f\in
F^{q_1,z}_{M,\alpha}$ if and only if $C^f$ has $\alpha$ factors
$R$ and  $M-\alpha$ factors $\Delta\psi$.

\par We claim that for each of the index sets
$F^{q_1,z_{min}(M)}_{M,\alpha}$, there is a linear combination of
vector fields, $\Sum_{r\in R^{q_1,z_{min}(M)}_{M,\alpha}} a_r
C^{r,i}_{g}(\psi_1,\dots ,\psi_s)$ where each $C^{r,i}$ is in the
form (\ref{linisymric}) and has $\delta= z_{min}(M)$ and $\alpha$
factors $R$ and $M-\alpha$ factors $\Delta\psi_h$, so that modulo
complete contractions of length $\ge\sigma +1$:

\begin{equation}
\label{timgow3}
\begin{split}
 &\Sum_{f\in F^{q_1, z_{min}(M)}_{M,\alpha}} a_f
C^f_{g}(\psi_1,\dots ,\psi_s) -div_i \Sum_{r\in
R^{q_1,z_{min}(M)}_{M,\alpha}} a_r C^{r,i}_{g}(\psi_1,\dots
,\psi_s)=
\\&\Sum_{d\in D_1} a_d C^d_{g}(\psi_1,\dots ,\psi_s)+
\Sum_{d\in D_2} a_d C^d_{g}(\psi_1,\dots ,\psi_s)+ \Sum_{d\in D_3}
a_d C^d_{g}(\psi_1,\dots ,\psi_s).
\end{split}
\end{equation}
In the above, each $C^d$, $d\in D_1$ is a complete
contraction with $\delta= z_{min}(M)+1$ and all the other features
being the same as the contractions $C^f$ indexed in $F^{q_1}_M$ 
(in particular they have $\delta=z_{min}$). Each $C^d,
d\in D_2$ is a complete contraction with $q=q_1-1$ and all the
other features being the same as the contractions $C^f$ indexed in $F^{q_1}_M$
(in particular they have $\delta=z_{min}$). 
Finally, each $C^d, d\in D_3$ is a complete contraction
with a total of $M-1$  factors $R$ or $\Delta\psi_h$, and
all the other features being the same as the contractions $C^f$ indexed in $F^{q_1}_M$
(in particular they have $\delta=z_{min}$).
\newline

\par We remark that (\ref{timgow3})
implies that modulo complete contractions of length $\ge\sigma +1$:
\begin{equation}
\label{timgow4}
\begin{split}
 &\Sum_{f\in F^{q_1, z_{min}(M)}_M} a_f
C^f_{g}(\psi_1,\dots ,\psi_s) -div_i \Sum_{\alpha=0}^{M}\Sum_{r\in
R^{q_1,z_{min}(M)}_{M,\alpha}} a_r C^{r,i}_{g}(\psi_1,\dots
,\psi_s)=
\\&\Sum_{d\in D_1} a_d C^d_{g}(\psi_1,\dots ,\psi_s)+
\Sum_{d\in D_2} a_d C^d_{g}(\psi_1,\dots ,\psi_s)+ \Sum_{d\in D_3}
a_d C^d_{g}(\psi_1,\dots ,\psi_s).
\end{split}
\end{equation}
Th terms in the RHS of the above have the same 
properties as the terms in the RHS of (\ref{timgow3}).
\newline

\par Thus, in order to derive Lemma \ref{cancwant} we need to show 
(\ref{vasko}), (\ref{vasko2}), and to derve Lemma \ref{corola} we need 
to show (\ref{timgow}), (\ref{timgow2}), (\ref{timgow3}). 
\newline

{\it Proof of equations (\ref{vasko}) and (\ref{vasko2}).}
\newline

\par Our aim is to apply Proposition \ref{pregiade} to
equation (\ref{bigdaddy2}). Since (\ref{axrhst}) is symmetric in the functions $\psi_1,\dots
,\psi_s$, we can just set $\psi_1=\dots
\psi_s=\psi$ and we lose no information. For notational
convenience, we will still write $\psi_1,\dots , \psi_s$ but the
functions $\psi_1,\dots ,\psi_s$ will in fact all be equal to
$\psi$. Now, by factoring out the factors $\Delta\psi$ we write:

\begin{equation}
\label{xarmosyno}
\begin{split}
& \Sum_{l\in L_M} C^{l,\iota}_{g}(\psi_1,\dots ,\psi_s)+\Sum_{j\in
J_M}C^j_{g}(\psi_1,\dots ,\psi_s) =
\\&\Sum_{l\in L_M} C^{l,\iota}_{g}(\psi_{M+1},\dots
,\psi_s)\Delta\psi_1\dots \Delta\psi_M+\Sum_{j\in J_M}
C^j_{g}(\psi_{M+1},\dots ,\psi_s)\Delta\psi_1\dots \Delta\psi_M.
\end{split}
\end{equation}

\par In view of (\ref{bigdaddy2}), we claim:

\begin{equation}
\label{bigdaddy3}
\begin{split}
&\Sum_{l\in L_\mu^{M}} a_l Xdiv_{i_1}\dots
Xdiv_{i_{\mu-M}}C^{\iota,l|i_1\dots
i_{\mu-M}}_{g}(\psi_{M+1},\dots ,\psi_s)+
\\&\Sum_{j\in J_M} a_j Xdiv_{i_1}\dots Xdiv_{i_{m_j-M}}
C^{j|i_1\dots i_{m_j-M}}_{g}(\psi_{M+1},\dots ,\psi_s)=0,
\end{split}
\end{equation}
modulo complete contractions of length $\ge\sigma -M+1$. (\ref{bigdaddy3}) 
follows by focusing on the sublinear combination in
(\ref{bigdaddy2}) that has $M$ factors $\nabla\psi_1,\dots
,\nabla\psi_M$ ({\it notice that this sublinear combination vanishes separately 
and all $\nabla\psi_h$'s are contracting against
derivative indices}), and the formally erasing the factors $\nabla\psi_h$ 
{\it and} the (derivative) indices against which they contract. 
This produces a new true equation,\footnote{This fact can be rigorously checked by applying 
the operation $Erase$--see the Appendix below.} which is precisely (\ref{bigdaddy3}). 

\par We  now directly apply Proposition \ref{pregiade} to
(\ref{bigdaddy3})\footnote{After first re-writing
(\ref{bigdaddy3}) in dimension $n-2M$.} (since by the hypothesis
that $|\Delta|_{Max}\le\sigma-3$ the real length
of the tensor fields in (\ref{bigdaddy3}) is $\ge 3$).
  In the case where $L_\mu^{M}\ne \emptyset$, we
 deduce that there is a linear combination of acceptable
$\mu -M+1$-tensor fields,\\ $\Sum_{h\in H} a_h C^{h,i_1\dots i_{\mu
-M+1}}_{g}(\psi_{M+1},\dots ,\psi_s)$, so that:

\begin{equation}
\label{yasser}
\begin{split}
&\Sum_{l\in L^{M}_\mu} a_l C^{\iota,l|i_1\dots
i_{\mu-M}}_{g}(\psi_{M+1},\dots ,\psi_s)\nabla_{i_1}\upsilon\dots
\nabla_{i_{\mu -M}}\upsilon=
\\&\Sum_{h\in H} a_h Xdiv_{i_{\mu-M+1}}C^{h,i_1\dots i_{\mu
-M+1}}_{g}(\psi_{M+1},\dots ,\psi_s)\nabla_{i_1}\upsilon\dots
\nabla_{i_{\mu -M}}\upsilon.
\end{split}
\end{equation}

\par Therefore, since the above holds formally, we observe that
the linear combination of vector fields needed for (\ref{vasko})
is precisely
\begin{equation}
\label{katseli}
\begin{split}
& [C^{h,i_1\dots i_{\mu
-M+1}}_{g}(\psi_{M+1},\dots ,\psi_s)\nabla_{i_1}\upsilon\dots
\nabla_{i_{\mu
-M}}\upsilon]
\\&\cdot \nabla^{k_1}\upsilon\nabla_{k_1}\psi_1\dots\nabla^{k_M}
\upsilon\nabla_{k_M}\psi_M.
\end{split}
\end{equation}

\par We prove (\ref{vasko2}) by a very similar argument. We
again use the notation (\ref{xarmosyno}), only now
$L_\mu^{M}=\emptyset$. We claim:

\begin{equation}
\label{bigdaddy4}
\begin{split}
&\Sum_{j\in J^{M}_{\delta_{min}(M)}} a_j Xdiv_{i_1}\dots
Xdiv_{i_{\delta_{min}(M)-M}} C^{j|i_1\dots
i_{\delta_{min}(M)-M}}_{g}(\psi_{M+1},\dots,\psi_s)+
\\&\Sum_{\delta>\delta_{min}(M)}\Sum_{j\in J^{M}_\delta} a_j Xdiv_{i_1}\dots
Xdiv_{i_{m_j-M}} C^{j|i_1\dots i_{m_{\delta_{min}-M}}}_{g}
(\psi_{M+1},\dots ,\psi_s)=0,
\end{split}
\end{equation}
modulo complete contractions of length $\ge\sigma -M+1$. This
again follows by focusing on the sublinear combination in
(\ref{bigdaddy2}) that has $M$ factors $\nabla\psi_1,\dots
,\nabla\psi_M$ (notice this sublinear combination vanishes separately) and  then applying the
eraser to $\nabla\psi_1,\dots ,\nabla\psi_M$.

\par We  now apply Proposition \ref{pregiade} to
(\ref{bigdaddy4}). We deduce that there is a linear combination of
acceptable $(\delta_{min}(M) -M+1)$-tensor fields, $\Sum_{h\in H}
a_h C^{h,i_1\dots i_{\delta_{min}(M) -M+1}}_{g}$, so that:

\begin{equation}
\label{yasser}
\begin{split}
&\Sum_{l\in L^{M}_\mu} a_l C^{\iota,l|i_1\dots
i_{\delta_{min}(M)-M}}_{g}(\psi_{M+1},\dots
,\psi_s)\nabla_{i_1}\upsilon\dots \nabla_{i_{\delta_{min}(M)
-M}}\upsilon-\Sum_{h\in H} a_h
\\& Xdiv_{i_{\delta_{min}(M)-M+1}}C^{h,i_1
\dots i_{\delta_{min}(M) -M+1}}_{g}(\psi_{M+1},\dots
,\psi_s)\nabla_{i_1}\upsilon\dots \nabla_{i_{\delta_{min}(M)
-M}}\upsilon=0.
\end{split}
\end{equation}

\par Therefore, since the above holds formally, we derive that
the linear combination of vector fields needed for (\ref{vasko2})
is precisely:
\begin{equation}
\begin{split} 
&\Sum_{h\in H} a_h C^{h,i_1\dots i_{\delta_{min}(M)
-M+1}}_{g}(\psi_{M+1},\dots ,\psi_s)\nabla_{i_1}\upsilon\dots
\nabla_{i_{\delta_{min}(M) -M}}\upsilon\\&\nabla^{k_1}
\upsilon\nabla_{k_1}\psi_1\dots\nabla^{k_M}
\upsilon\nabla_{k_M}\psi_M.
\end{split}\end{equation}
{\it The proof of (\ref{timgow}), (\ref{timgow2}),
(\ref{timgow3}):}
\newline

\par We start with (\ref{timgow}) and (\ref{timgow3}). We will
prove (\ref{timgow3}); this  proof applies to show  (\ref{timgow})
by just setting $z_{min}(M)=z_1$. We begin by noting an equation
analogous to (\ref{xarmosyno}): Let
$\alpha_1\ge 0$ be the smallest value of $\alpha$ for which 
 each $F^{q_1,z_{min}(M)}_{M,\alpha}$ with $\alpha
>\alpha_1$ is empty. We will then show (\ref{timgow3}) for $\alpha
=\alpha_1$. Clearly, if we can prove this, then 
(\ref{timgow3}) follows for every $\alpha$, by induction. 
We observe that each of the other complete contractions appearing
in the equation of Lemma \ref{corola} must either have less than
$M$ factors $\Delta\psi_h, R$ (in total), or will have $M$ such
factors in total but less than $\alpha_1$ factors $R$. This just
follows from the definition of $M$ and $\alpha_1$.

\par Again, using the
fact that the complete contractions are symmetric in the functions
$\psi_1,\dots ,\psi_s$, we may assume with no loss of generality
that these functions are all equal to $\psi$. We factor 
out the factors $\Delta\psi_h$, $R$ to write out:

\begin{equation}
\label{xarmosyno2}
\begin{split}
&\Sum_{f\in F^{q_1,z_{min}(M)}_{M,\alpha_1}} a_f
C^f_{g}(\psi_1,\dots ,\psi_s)+ \Sum_{z>z_{min}(M)} \Sum_{f\in
F^{q_1,z_{min}(M)}_{M,\alpha_1}} a_f C^f_{g}(\psi_1,\dots
,\psi_s)=
\\& \{\Sum_{f\in F^{q_1,z_{min}(M)}_{M,\alpha_1}} a_f
C^f_{g}(\psi_{M-\alpha_1+1},\dots ,\psi_s)
\\& +\Sum_{z>z_{min}(M)} \Sum_{f\in
F^{q_1,z_{min}(M)}_{M,\alpha_1}} a_f
C^f_{g}(\psi_{M-\alpha_1+1},\dots ,\psi_s)\}\cdot
R^{\alpha_1}\Delta\psi_1\dots \Delta\psi_{M-\alpha_1}.
\end{split}
\end{equation}

\par We now claim that modulo complete contractions 
of length $\ge\sigma -M-\alpha+1$:

\begin{equation}
\label{gatopardos}
\begin{split}
&Xdiv_{i_1}\dots Xdiv_{i_{z_{min}(M)-M}}\Sum_{f\in
F^{q_1,z_{min}(M)}_{M,\alpha_1}} a_f C^{f, i_1\dots
i_{z_{min}(M)-M}}_{g}(\psi_{M-\alpha_1+1},\dots
,\psi_s,\\&\Omega^{q_1-\alpha_1})+
 Xdiv_{i_1}\dots Xdiv_{i_{z-M}}\Sum_{z>z_{min}(M)} \Sum_{f\in
F^{q_1,z_{min}(M)}_{M,\alpha_1}} a_f C^{f,i_1\dots
i_{z-M}}_{g}(\psi_{M-\alpha_1+1},\\&\dots
,\psi_s,\Omega^{q_1-\alpha_1})=0.
\end{split}
\end{equation}

\par This follows by picking out the terms in 
(\ref{rhodesia}) with $\alpha_1$ factors $\nabla\Omega$, $M-\alpha$ factors 
$\nabla\psi$ (this sublinear combination must vanish separately) and
 then formally erasing these factors and the indices against which they 
contract.\footnote{A rigorous proof that this formal operation produces a true equation can be 
derived by virtue of the operation $Erase$--see the Appendix below.}

\par We may now apply Proposition \ref{pregiade} to (\ref{gatopardos}).
 We derive that
 there is a linear combination of $(z_{min}(M)-M+1)$-tensor fields,
$\Sum_{z\in Z}a_z C^{z,i_1\dots
i_{z_{min}(M)-M+1}}_{g}$ (written in dimension $-n+2M$), as
stated in Proposition \ref{pregiade}, so that modulo complete
contractions of length $\ge\sigma -M+z_{min}(M)+1$:

\begin{equation}
\label{explkill} \begin{split} & \Sum_{f\in
F^{q_1,z_{min}(M)}_{M,\alpha_1}} a_f C^{f, i_1\dots
i_{z_{min}(M)-M}}_{g}(\psi_{M-\alpha_1+1},\dots
,\psi_s,\Omega^{q_1-\alpha_1})
\nabla_{i_1}\upsilon\dots\nabla_{i_{z_{min}(M)-M}}\upsilon
\\&-Xdiv_{i_{z_{min}(M)-M+1}}\Sum_{t\in T}a_t C^{t,i_1\dots
i_{z_{min}(M)-M+1}}_{g}(\psi_{M-\alpha_1+1}, \dots
,\psi_s,\Omega^{q_1-\alpha_1})
\\&\nabla_{i_1}\upsilon\dots\nabla_{i_{z_{min}(M)-M}}\upsilon=0.
\end{split}
\end{equation}

\par We act on the above equation with the operation
$Riccify$. Observe that:

\begin{equation}
\label{president} \begin{split} & \Sum_{f\in
F^{q_1,z_{min}(M)}_{M,\alpha_1}} a_f Riccify[C^{f, i_1\dots
i_{z_{min}(M)-M}}_{g}(\psi_{M-\alpha_1+1},\dots ,\psi_s,\Omega)
\nabla_{i_1}\upsilon\\&\dots\nabla_{i_{z_{min}(M)-M}}\upsilon]
\cdot R^{\alpha_1} \Delta\psi_1,\dots\Delta\psi_{M-\alpha_1}= \Sum_{f\in
F^{q_1,z_{min}(M)}_{M,\alpha_1}} a_f C^f_{g}(\psi_1,\dots ,\psi_s).
\end{split}
\end{equation}

\par Thus, by virtue of the above equation and Lemma \ref{paparidis2},
 we see that the vector field required for equations
(\ref{timgow}), (\ref{timgow3}) is precisely:

\begin{equation}
 \label{r.r}
\begin{split}
 &\Sum_{t\in T} a_t Riccify[C^{t,i_1\dots
i_{z_{min}(M)-M+1}}_{g}(\psi_{M-\alpha_1+1}, \dots
,\psi_s,\Omega^{q_1-\alpha_1})
\\&\nabla_{i_1}\upsilon\dots\nabla_{i_{z_{min}(M)-M}}\upsilon]\cdot 
R^{\alpha_1}\Delta\psi_1\dots\Delta\psi_{M-\alpha_1}.
\end{split}
\end{equation}
 
{\it Proof of (\ref{timgow2}):}
\newline
\par The proof in this case
is very slightly different from the proof of (\ref{timgow}),
(\ref{timgow3}).

\par We recall that for each $f\in F^{q_1,\mu}_M$, $C^f_g(\psi_1,\dots,\psi_s)$  must have $M$ factors
$\Delta\psi_h$ (i.e. there are no factors $R$, by definition) and 
will furthermore contain 
factors $\nabla^{a_1\dots a_t}\nabla^{(p)}_{r_1\dots r_p}Ric_{ab}$
 and for
each such factor one of the indices ${}^{a_1},\dots ,{}^{a_t}$ is
contracting against the index ${}_b$ (this implies that for each
such factor we have $t>0$). With no loss of generality, we assume
that the index ${}^{a_1}$ is contracting against the index ${}_b$.
Thus, applying the contracted second Bianchi identity, we may
replace the factor $\nabla^{a_1\dots a_t}\nabla^{(p)}_{r_1\dots
r_p}Ric_{ab}$ by a factor $\frac{1}{2}\nabla^{a_2\dots
a_t}\nabla^{(p+1)}_{r_1\dots r_p a}R$, modulo introducing complete
contractions with more that $\sigma$ factors. Moreover, as in the
previous case we set $\psi_1,\dots ,\psi_s=\psi$, although we will
still write $\psi_1,\dots ,\psi_s$ for notational convenience.

\par We pick out the complete contractions indexed in
$\bigcup_{z>\mu} F^{q_1,z}$ that have exactly $M$ factors
$\Delta\psi_h$. By our notational conventions, they will be
indexed in $\bigcup_{z>\mu} F^{q_1,z}_{M,0}$. We again write out:

\begin{equation}
\label{staxeria}
\begin{split}
&\Sum_{f\in  F^{q_1,\mu}_M} a_f C^f_{g}(\psi_1,\dots ,\psi_s) +
\Sum_{f\in \bigcup_{z>\mu} F^{q_1,z}_{M,0}} a_f
C^f_{g}(\psi_1,\dots ,\psi_s)=
\\& (\frac{1}{2})^{q_1}\Sum_{f\in  F^{q_1,\mu}_M} a_f C^f_{g}(\psi_{M+1},\dots
,\psi_s)\Delta\psi_1\dots\Delta\psi_M 
\\&+ \Sum_{f\in \bigcup_{z>\mu}
F^{q_1,z}_M} a_f C^f_{g}(\psi_{M+1},\dots
,\psi_s)\Delta\psi_1\dots\Delta\psi_M,
\end{split}
\end{equation}
where each $C^f_{g}(\psi_{M+1},\dots ,\psi_s)$, $f\in
F^{q_1,\mu}_M$ is now in the form:

\begin{equation}
\label{gro9ies}
\begin{split}
&contr(\nabla^{f_1\dots f_y}\nabla^{(m_1)}R_{ijkl}\otimes\dots
\otimes \nabla^{g_1\dots g_p}\nabla^{(m_s)}R_{ijkl}
\\&\otimes \nabla^{y_1\dots
y_w}\nabla^{(d_1)}_{a_1\dots
a_{d_1}}R\otimes\dots\otimes\nabla^{x_1\dots
x_p}\nabla^{(d_{q_1})}_{b_1\dots b_{d_{q_1}}}R  \otimes
\\& \nabla^{a_1\dots a_{t_1}}\nabla^{(u_1)}\psi_1\otimes\dots\otimes \nabla^{c_1\dots
c_{t_s}}\nabla^{(u_s)}\psi_s),
\end{split}
\end{equation}
 while each
 $C^f_{g}(\psi_{M+1},\dots ,\psi_s)$, $f\in
F^{q_1,z}_{M,0}, z>\mu$ is still in the form:

\begin{equation}
\label{gro9ies2}
\begin{split}
&contr(\nabla^{f_1\dots f_y}\nabla^{(m_1)}R_{ijkl}\otimes\dots
\otimes \nabla^{g_1\dots g_p}\nabla^{m_s}R_{ijkl}
\\&\otimes \nabla^{y_1\dots
y_w}\nabla^{(d_1)}Ric_{ij}\otimes\dots\otimes\nabla^{x_1\dots
x_p}\nabla^{(d_{q_1})}Ric_{ij}  \otimes
\\& \nabla^{a_1\dots a_{t_1}}\nabla^{u_1}\psi_1\otimes\dots\otimes \nabla^{c_1\dots
c_{t_s}}\nabla^{u_s}\psi_s),
\end{split}
\end{equation}
(with $\delta>\mu$).

\par Now, for each $f\in F^{q_1,\mu}$ we denote by
$C^f_{g}(\psi_{M+1},\dots ,\psi_s,\Omega^{q_1})$ the complete
contraction that arises from $C^f_{g}(\psi_{M+1},\dots ,\psi_s)$
by replacing the factors $\nabla^{a_1\dots
a_t}\nabla^{(p)}_{r_1\dots r_p}R$ by $\nabla^{a_1\dots
a_th}\nabla^{(p)}_{r_1\dots r_p h}\Omega$. For each $f\in
F^{q_1,z}_{M,0},z>\mu$, \\$C^f_{g}(\psi_{M+1},\dots
,\psi_s,\Omega^{q_1})$ is the same as before. Applying Lemma
\ref{labrome} to the equation $\int_{M^n}I^s_{g}dV_{g}=0$
 and then the eraser to the $M$ factors $\nabla\psi$, we derive
 that modulo complete contractions of length $\ge\sigma -M+1$:

\begin{equation}
\label{xamenorouxo}
\begin{split}
&Xdiv_{i_1}\dots Xdiv_{\mu-M}\Sum_{f\in  F^{q_1,\mu}_M} a_f C^{f,
i_1\dots i_{\mu-M}}_{g}(\psi_{M+1},\dots ,\psi_s,\Omega^{q_1})
\\&+Xdiv_{i_1}\dots Xdiv_{z-M}\Sum_{f\in \bigcup_{z>\mu}
F^{q_1,z}_{M,0}} a_f C^{f,i_1,\dots i_{z-M}}_{g}(\psi_{M+1},\dots
,\psi_s,\Omega^{q_1})=0.
\end{split}
\end{equation}

\par Now, applying Proposition \ref{pregiade} to the above, we
deduce that there is a linear combination of acceptable
$(\mu-M+1)$-tensor fields, $\Sum_{t\in T} a_t C^{t,i_1\dots
i_{\mu-M+1}}_{g}(\psi_1,\dots ,\psi_s,\Omega)$, so that modulo
complete contractions of length $\ge\sigma -M+\mu+1$:

\begin{equation}
\label{kabalarhdes}
\begin{split}
&\Sum_{f\in  F^{q_1,\mu}_M} a_f C^{f, i_1\dots
i_{\mu-M}}_{g}(\psi_{M+1},\dots
,\psi_s,\Omega^{q_1})\nabla_{i_1}\upsilon\dots\nabla_{i_{\mu
-M}}\upsilon
\\&-Xdiv_{i_{\mu-M +1}} \Sum_{t\in T} a_t C^{t,i_1\dots
i_{\mu-M+1}}_{g}(\psi_1,\dots ,\psi_s,\Omega^{q_1})
\nabla_{i_1}\upsilon\dots\nabla_{i_{\mu -M}}\upsilon=0.
\end{split}
\end{equation}

\par Finally, an observation: The above equation holds formally. We observe
that by construction, each complete contraction
$$C^{f, i_1\dots
i_{\mu-M}}_{g}(\psi_{M+1},\dots
,\psi_s,\Omega^{q_1})\nabla_{i_1}\upsilon\dots\nabla_{i_{\mu
-M}}\upsilon$$ has at least one factor $\nabla\upsilon$
contracting against {\it each} factor $\nabla^{(p)}\Omega$.
Therefore, since the equation holds formally we may assume
with no loss of generality that the same must be true of each
vector field
$$C^{t,i_1\dots i_{\mu-M+1}}_{g}(\psi_1,\dots ,\psi_s,\Omega)
\nabla_{i_1}\upsilon\dots\nabla_{i_{\mu -M}}\upsilon.$$

\par Moreover, we know that for each $C^f$ and each of its factors $\nabla^{(p)}\Omega$,
 one factor $\nabla\upsilon$ is contracting against the
last index ${}_{r_p}$. Modulo introducing complete contractions of
length $\ge\sigma +\mu-M+1$, we may assume that the same is true
of each of the vector fields $C^{t,i_1\dots i_{\mu-M+1}}$. But then,
since the above holds formally, we may assume that when we apply
the permutations to make the left hand side of (\ref{kabalarhdes})
{\it formally} zero, the index ${}_{r_p}$ in each factor
$\nabla^{(p)}\Omega$ is {\it not} permuted. (One can prove this by applying 
the operation $Erase$ repeatedly).\footnote{See the appendix for 
the strict definition of the operation $Erase[\dots]$.}

\par Now, we define an operation $Riccify''$ which is slightly different
from the standard $Riccify$: We replace each of the expressions of
the form $$\nabla^{(r)}\psi\nabla^{t_1}\upsilon\dots
\nabla^{t_y}\upsilon,
\nabla^{(m)}R_{ijkl}\nabla^{t_1}\upsilon\dots
\nabla^{t_y}\upsilon$$ (where each of the factors $\nabla\upsilon$
is contracting against the factor $\nabla^{(r)}\psi_h$, $\nabla^{(m)}
R_{ijkl}$ respectively) as in the operation $Riccify$. But we also
replace each of the expressions
$$\nabla^{(p)}_{r_1\dots r_p}\Omega
\nabla^{r_{a_1}}\upsilon\dots \nabla^{r_{a_t}}\upsilon
\nabla^{r_p}\upsilon$$ by a factor
$$\nabla^{r_{a_1}\dots r_{a_t}}\nabla^{(p-1)}_{r_1\dots r_{p-1}}R$$
 We then observe that since (\ref{kabalarhdes}) holds formally {\it without}
permuting the last index ${}_{r_p}$ in each factor
$\nabla^{(p)}\Omega$ (and that index is contracting against a
factor $\nabla\upsilon$), we then have that:

\begin{equation}
\label{kabalarhdes2}
\begin{split}
&\Sum_{f\in  F^{q_1,\mu}_M} a_f Riccify'' [C^{f, i_1\dots
i_{\mu-M}}_{g}(\psi_{M+1},\dots
,\psi_s,\Omega^{q_1})\nabla_{i_1}\upsilon\dots\nabla_{i_{\mu
-M}}\upsilon]
\\&-Xdiv_{\mu-M +1} \Sum_{t\in T}
Riccify''[a_t C^{t,i_1\dots i_{\mu-M+1}}_{g}(\psi_1,\dots
,\psi_s,\Omega^{q_1}) \nabla_{i_1}\upsilon\dots\nabla_{i_{\mu
-M}}\upsilon]=0,
\end{split}
\end{equation}
modulo complete contractions of length $\ge\sigma -M+1$.

Hence, the vector field needed for (\ref{timgow2}) is precisely:

$$\Sum_{t\in T} Riccify'[a_t C^{t,i_1\dots
i_{\mu-M+1}}_{g}(\psi_1,\dots ,\psi_s,\Omega^{q_1})
\nabla_{i_1}\upsilon\dots\nabla_{i_{\mu
-M}}\upsilon]\Delta\psi_1\dots \Delta\psi_M.$$

{\it Note: Codification of the remaining cases of Lemmas
\ref{corola} and \ref{cancwant}.}
\newline

\par Here we codify what remains to be proven for Lemmas
\ref{corola} and \ref{cancwant}. We will then prove these claims
in the paper \cite{alexakis3} in this series.

\par {\it  Lemma \ref{corola}:} What remains to be
proven is the second claim in that Lemma: 

\par Recall the index set $F$ in Lemma \ref{corola} (this indexes the
complete contractions in $I^s_g(\psi_1,\dots,\psi_s)$, in the form
(\ref{linisymric}), with at least one factor $\nabla^{(p)}Ric$ or
$R$). Recall that for each $q, 1\le q\le \sigma-s$, $F^q\subset F$
stands for the index set of complete contractions with precisely
$q$ factors $\nabla^{(p)}Ric$ or $R$. Recall that for each
complete contraction in the form (\ref{linisymric}) we have
denoted by $|\Delta|$ the number of factors in one of the forms
$\Delta\psi_h$, $R$. For each index set $F^q$ above, let us denote
by $F^{q,*}\subset F^q$ the index set of complete contractions
with $|\Delta|\ge\sigma-2$, and $F^{*}=\bigcup _{q>0}F^{q,*}$.

\par Claim: There exists a linear combinations of vector fields
(indexed in $H$ below), each in the form (\ref{linisymric}) with
$\sigma$ factors, so that modulo complete contractions of length
$>\sigma$:

\begin{equation}
\label{alexstone} \sum_{f\in F^{*}} a_f
C^f_g(\psi_1,\dots,\psi_s)-div_i\sum_{h\in H} a_h
C^{h,i}_g(\psi_1,\dots,\psi_s)= \sum_{y\in Y} a_y
C^y_g(\psi_1,\dots,\psi_s),
\end{equation}
where the complete contractions indexed in $Y$ are in the form
(\ref{linisymric}) with length $\sigma$, and satisfy all the
properties of the sublinear combination $\sum_{f\in F}\dots$ but
in addition have $|\Delta|\le\sigma-3$.
\newline

{\it The remaining claims for Lemma \ref{cancwant}:} 
\begin{lemma}
\label{violi}
Denote 
by $L^{*}_\mu\subset L_\mu,J^{*}\subset J$ the index sets
of complete contractions in the hypothesis of Lemma \ref{cancwant}
with $|\Delta|\ge\sigma-2$, among the complete contractions
indexed in $L_\mu,J$ respectively. We claim that there
exists a linear combination of $(\mu+1)$-tensor fields fields (indexed in
$H$ below), with length $\sigma$, in the form (\ref{linisymric})
{\it without factors $\nabla^{(p)}Ric,R$} and with $\delta=\mu$ so
that:

\begin{equation}
\label{okounka}
\begin{split} 
 &\sum_{l\in L^{*}_\mu} a_l
C^{l,i_1\dots i_\mu}_g(\psi_1,\dots,\psi_s)\nabla_{i_1}\upsilon\dots\nabla_{i_\mu}\upsilon
-div_i \sum_{h\in H}
a_hC^{h,i|i_1\dots i_\mu}_g(\psi_1,\dots,\psi_s)
\\&\cdot \nabla_{i_1}\upsilon\dots\nabla_{i_\mu}\upsilon
=\sum_{l\in \overline{L}} a_l C^{l,i_1\dots i_\mu}_g(\psi_1,\dots,\psi_s)
\nabla_{i_1}\upsilon\dots\nabla_{i_\mu}\upsilon,
\end{split}
\end{equation}
where the $\mu$-tensor fields indexed in $\overline{L}$ are in
the form (\ref{linisymric}) with no factors $\nabla^{(p)}Ric$ or
$R$ and with $|\Delta|\le\sigma-3$.
\end{lemma}
 (Notice that if we can show
the above, then in proving Lemma \ref{cancwant} we can assume with
no loss of generality that $L^{*}_\mu=\emptyset$). In the setting
$L^{*}_\mu=\emptyset$, what remains to be shown to complete the
proof of Lemma \ref{cancwant} is the following:

\begin{lemma}
\label{postpone2}
 Assume that $L^{*}_\mu=\emptyset$ in the hypothesis of Lemma
\ref{cancwant}. Denote by $J^{*}\subset J$ the index set of the
complete contractions $C^j_g(\psi_1,\dots,\psi_s)$ with
$|\Delta|\ge\sigma-2$. We then claim that there exists a linear combination of
vector fields (indexed in $H$ below) so that:

\begin{equation}
\label{o3} \sum_{j\in J^{*}} a_j C^j_g(\psi_1,\dots,\psi_s)-div_i
\sum_{h\in H} a_h C^{h,i}_g(\psi_1,\dots,\psi_s)=\sum_{y\in Y'}
a_y C^y_g(\psi_1,\dots,\psi_s),
\end{equation}
where the complete contractions indexed in $Y'$ are in the form
(\ref{linisymric}) with length $\sigma$, with no factors
$\nabla^{(p)}Ric$ or $R$ and have $\delta\ge\mu+1$ and in addition
 $|\Delta|\le\sigma-3$. 
\end{lemma}
If we can show the above claims, then
we will have completely shown Lemmas \ref{corola} and
\ref{cancwant}, and hence also Proposition \ref{burj}, which
implies Proposition \ref{killtheta}.

\section{Appendix: Some Technical Tools.}

We prove here some technical claims, which will be useful in this series of papers.
\newline

{\bf The Eraser:}
We  consider complete contractions
$C^h_{g}(\Omega_1,\dots ,\Omega_p,\phi_1,\dots ,\phi_u)$ in
 the form:

\begin{equation}
\label{form1early}
\begin{split}
&pcontr(\nabla^{(m_1)}R_{ijkl}\otimes\dots\otimes\nabla^{(m_s)}R_{ijkl}
\otimes
\\& \nabla^{(b_1)}\Omega_1\otimes\dots\otimes \nabla^{(b_p)}\Omega_p
\otimes\nabla\phi_1\otimes\dots \otimes\nabla\phi_u),
\end{split}
\end{equation}
 with length $\sigma+u$. We fistly define a {\it formal} operation on
such complete contractions, which we call
the {\it eraser} operation:

\begin{definition}
\label{deferase}
\par Consider a set of complete
contractions $C^h_{g}(\Omega_1,\dots ,\Omega_p,\phi_1, \dots
,\phi_u)$, \\$h\in H$, each in the form (\ref{form1early}). 
Assume that for each $h\in H$, some particular factor
$\nabla\phi_b$ ($b$ is fixed, i.e. $b$ is independent of $h\in H$) is
contracting
 against a factor $\nabla^{(m)}R_{ijkl}$ and moreover
against a  derivative index in that factor.

\par We then define $Erase_{\nabla\phi_b}[C^h_{g}(\Omega_1,
\dots ,\Omega_p,\phi_1,\dots ,\phi_u)]$ to stand for the complete
contraction (of weight $-n+2$) that formally arises from
$C^h_{g}(\Omega_1,\dots ,\Omega_p,\phi_1,\dots ,\phi_u)$
 by erasing the factor $\nabla\phi_b$ and also erasing the
derivative index that it contracts against.\footnote{Note that we thus
obtain a complete contraction of length $\sigma +u-1$.} 
\end{definition}

\begin{lemma}
\label{erasing} Consider a set of complete contractions
$\{C^h_{g}(\Omega_1,\dots ,\Omega_p,\phi_1, \dots ,\phi_u)\}_{h\in
H}$ as in the above definition and assume that modulo complete contractions of
length $\ge\sigma +u+1$:

\begin{equation}
\label{olae3w} U_g(\Omega_1,\dots
,\Omega_p,\phi_1,\dots ,\phi_u)=\Sum_{h\in H} a_h C^h_{g}(\Omega_1,\dots
,\Omega_p,\phi_1,\dots ,\phi_u)=0.
\end{equation}

We claim
that modulo complete contractions of length $\ge\sigma +u$:

\begin{equation}
\label{olae3w2} \Sum_{h\in H} a_h
Erase_{\nabla\phi_b}[C^h_{g}(\Omega_1,\dots ,\Omega_p,\phi_1,\dots
,\phi_u)]=0.
\end{equation}
\end{lemma}

{\it Proof:} We call the factor $\nabla^{(m)}R_{ijkl}$ against which
$\nabla\phi_b$ contracts the {\it special factor}. We break the
index set $H$ into subsets $H_\mu$, where $h\in H_\mu$ if and only
if $C^h$ has $m=\mu$ derivatives on the special factor. Observe 
that since (\ref{olae3w}) holds formally, it follows that for each different $\mu$:

\begin{equation}
\label{olae3w2} \Sum_{h\in H_\mu} a_h C^h_{g}(\Omega_1,\dots
,\Omega_p,\phi_1,\dots ,\phi_u)=\Sum_{t\in T} a_t
C^t_{g}(\Omega_1,\dots ,\Omega_p,\phi_1,\dots , \phi_u),
\end{equation}
where each $C^t$ has length $\ge\sigma + u+1$. This holds because the
 linearized version\footnote{See the introduction in \cite{a:dgciI} for a discussion 
of linearized complete contractions.} of
(\ref{olae3w})  must hold formally (for the linearized
complete contractions), and also because under any of
the linearized permutations by which we can make the linearized
version of (\ref{olae3w}) formally zero, the number of derivatives
on the special factor remains invariant. Now, it would suffice to
show that for each $\mu$:

\begin{equation}
\label{olae3w3} \Sum_{h\in H_\mu} a_h
Erase_{\nabla\phi_b}[C^h_{g}(\Omega_1,\dots ,\Omega_p,\phi_1,\dots
,\phi_u)]=0,
\end{equation}
modulo complete contractions of length $\ge\sigma +u$.

\par In order to show this we write:

$$U^\mu_{g}:=\Sum_{h\in H} a_h
C^h_{g}(\Omega_1,\dots ,\Omega_p,\phi_1,\dots ,\phi_u)-\Sum_{t\in
T} a_t C^t_{g}(\Omega_1,\dots ,\Omega_p,\phi_1,\dots
,\phi_u)(=0).$$

\par Now, consider $Image^1_{\phi'}[U^\mu_{g}]$. We denote by
$Image^{1,A}_{\phi'}[U^\mu_{g}]$ the sublinear combination
 that arises in $Image^1_{\phi'}[U^\mu_{g}]$ by replacing one
 of the factors of the form $\nabla^{(m)}R_{ijkl}$ by one of the
 four linear terms on the right hand side of
(\ref{curvtrans}).  Now, let us denote by ${}_a$ the index in the
special factor that contracts against the factor $\nabla\phi_b$.
We denote by $Image^{1,B}_{\phi'}[U^\mu_{g}]$ the linear
combination that arises from $Image^1_{\phi'}[U^\mu_{g}]$ by
applying the transformation law (\ref{levicivita}) to the special
factor and bringing out a factor $\nabla_a\phi'$ (observe that for
every contraction in $Image^{1,B}_{\phi'}[U^\mu_{g}]$ 
the two factors $\nabla\phi_b,\nabla\phi'$ contract against
each other). We denote by $Image^{1,C}_{\phi'}[U^\mu_{g}]$ the
sublinear combination that arises in $Image^1_{\phi'} [U^\mu_{g}]$
when we apply the transformation law (\ref{levicivita}) to any
complete contraction $C^h_{g}(\Omega_1,\dots \Omega_p,
\phi_1,\dots ,\phi_u)$ and bring out a factor $\nabla_f\phi'$,
where $f\ne a$. We thus have that each complete contraction in
$Image^{1,C}_{\phi'}[U^\mu_{g}]$ has
 length $\sigma+u+1$ and a factor $\nabla\phi'$ but it
{\it does not} contract against a factor $\nabla\phi_b$.

\par Finally, we denote by $\Sum_{w\in W} a_w
C^w_{g}(\Omega_1,\dots ,\Omega_p,\phi_1,\dots ,\phi_u, \phi')$ a
generic linear combination of complete contractions with either
length $\sigma +u+1$ and a factor $\nabla^{(q)}\phi'$ ($q\ge 2$)
{\it or} with length $>\sigma +u+1$.

\par By virtue of (\ref{olae3w2}), we derive that

\begin {equation}
\label{benny} Image^1_{\phi'}[U^\mu_{g}(\Omega_1,\dots
,\Omega_p,\phi_1,\dots ,\phi_u, \phi')]=0.
\end{equation}

\par In addition, we deduce that:

\begin {equation}
\label{benny2} lin \{Image^{1,A}_{\phi'}[U^\mu_{g}(\Omega_1,\dots
,\Omega_p,\phi_1,\dots ,\phi_u, \phi')]\}=0.
\end{equation}
Hence, since the above holds formally, we may repeat the
permutations by which we make the above formally zero to the
linear combination \\ $Image^{1,A}_{\phi'}[U^\mu_{g}(\Omega_1,\dots
,\Omega_p,\phi_1,\dots ,\phi_u)]$; we deduce that:

\begin{equation}
\label{benny3}
\begin{split}
&Image^{1,A}_{\phi'}[U^\mu_{g}(\Omega_1,\dots
,\Omega_p,\phi_1,\dots ,\phi_u)]= \Sum_{y\in Y} a_y
C^y_{g}(\Omega_1,\dots ,\Omega_p,\phi_1,\dots ,\phi_u,\phi')+
\\& \Sum_{w\in W} a_w
C^w_{g}(\Omega_1,\dots ,\Omega_p,\phi_1,\dots ,\phi_u, \phi'),
\end{split}
\end{equation}
where each $C^y_{g}(\Omega_1,\dots ,\Omega_p,\phi_1,\dots ,\phi_u,
\phi')$ has length $\sigma +u+1$ and a factor $\nabla\phi_b$, but
that factor contracts against a factor $\nabla^{(m)}R_{ijkl}$.
This follows by virtue of the formula (\ref{curvature}).

\par Hence we deduce that, modulo complete contractions of length
$\ge\sigma +u+2$:

\begin{equation}
\label{9anatw}
\begin{split}
& Image^{1,C}_{\phi'}[U^\mu_{g}] +Image^{1,B}_{\phi'}[U^\mu_{g}]
+\Sum_{y\in Y} a_y C^y_{g}(\Omega_1,\dots ,\Omega_p,\phi_1,\dots
,\psi_u, \phi')
\\& +\Sum_{w\in W} a_w
C^w_{g}(\Omega_1,\dots ,\Omega_p,\phi_1,\dots ,\psi_u, \phi')=0.
\end{split}
\end{equation}

\par Then, since the above must hold formally, it follows that,
 modulo complete contractions of length $\ge\sigma
+u+1$:

\begin{equation}
\label{9anatw2}
 Image^{1,B}_{\phi'}[U^\mu_{g}]+Image^{1,C}_{\phi'}[U^\mu_{g}]
+\Sum_{y\in Y} a_y C^y_{g}(\Omega_1,\dots ,\Omega_p,\phi_1,\dots
,\psi_u, \phi')=0.
\end{equation}

\par Now, since the above must hold formally, and since
each complete contraction in $Image^{1,C}_{\phi'}[U^\mu_{g}]$ has
the factor $\nabla\phi'$ not contracting against the factor
$\nabla\phi_b$, we derive:

\begin{equation}
\label{9anatw3} Image^{1,B}_{\phi'}[U^\mu_{g}]=0,
\end{equation}
(modulo contractions of length $\ge\sigma +u+2$).

\par Lastly we observe, by virtue of the formula
(\ref{levicivita})
and by virtue of the factor $e^{2\phi(x)}$ in
(\ref{curvtrans}), that:

\begin{equation}
\label{9anatw4} Image^{1,B}_{\phi'}[U^\mu_{g}]=-(\mu
+1)\nabla^c\phi_b\nabla_c\phi'\Sum_{h\in H_\mu} a_h
Erase_{\nabla\phi_b}[C^h_{g}]
\end{equation}

\par (\ref{9anatw3}), (\ref{9anatw4}) imply our claim. $\Box$
\newline

{\it Note:} The analogous result is completely obvious if the
factor $\nabla\phi_b$ is contracting against a factor
$\nabla^{(y)}\Omega_f$: We just replace $\Omega_f$ by
$\Omega_f\cdot \phi'$ and pick out the sublinear combination with
an expression $\nabla^w\phi'\nabla_w\phi_b$. $\Box$
\newline

{\bf The operation $Sub_\omega$:} We state another 
useful tool. For future
 reference, we will introduce the complete contractions we
 will be studying.  We will be considering
 complete contractions  $C^{m}_{g}(\Omega_1,\dots
,\Omega_p,\phi_1,\dots ,\phi_u)$ in the form: 
\begin{equation}
\label{xeilh} 
contr(\nabla^{(m_1)}R\otimes\dots\otimes\nabla^{(m_a)}R\otimes\nabla^{(r_1)}
\Omega_1\otimes\dots\otimes\nabla^{(r_p)}\Omega_p\otimes\nabla\phi_1\otimes\dots\otimes\nabla\phi_u),
\end{equation}
 with length $\sigma +u$ and {\it with one internal
contraction}. The contractions in the form above all have
a given number $u$ of factors $\nabla\phi_1,\dots,\nabla\phi_u$
and a given number $p$ of factors $\nabla^{(y)}\Omega_1,\dots
,\nabla^{(w)}\Omega_p$ (here the functions $\phi_1,\dots ,\phi_u$
and $\Omega_1,\dots ,\Omega_p$ are understood to be different). We
assume an equation:

\begin{equation}
\label{LSE} \Sum_{d\in D} a_d C^d_{g}(\Omega_1,\dots
,\Omega_p,\phi_1,\dots ,\phi_u)=0,
\end{equation}
which holds modulo complete contractions of length $\ge\sigma
+u+1$. We divide the index set $D$ into
subsets $D_1,D_2\subset D$. We will say that $d\in D_1$ if the
internal contraction is between internal indices in a factor
$\nabla^{(m)}R_{ijkl}$  and $D_2=D\setminus D_1$. In other words,
$C^d_{g}(\Omega_1,\dots ,\Omega_p,\phi_1,\dots ,\phi_u)$, $d\in
D_1$  will have a factor $\nabla^{(m)}Ric_{ik}$, while
$C^d_{g}(\Omega_1,\dots ,\Omega_p,\phi_1,\dots ,\phi_u)$, $d\in
D_2$ will have an internal contraction between indices
$(\nabla^s,{}_s)$.

\par We then define an operation $Sub_\omega$ that acts on the
complete contractions $C^d_{g}$, $d\in D$ as follows: For $d\in
D_1$ $Sub_\omega[C^d_{g}]$ will stand for the complete
contraction that arises from $C^d_{g}$ by replacing the factor
$\nabla^{(p)}_{r_1\dots r_p}Ric_{ik}$ by a factor
$-\nabla^{(p+2)}_{r_1\dots r_pik}\omega$. If $d\in D_2$,
$Sub_\omega[C^d_{g}]$ will stand for the complete contraction that
arises from $C^d_{g}$ by picking out the internal contraction
$(\nabla^s,{}_s)$, then  erasing the derivative index $\nabla^s$ and
then adding a factor $\nabla^s\omega$ and contracting it against
the index ${}_s$ that has been left hanging.

\begin{lemma}
\label{subomega} Assuming (\ref{LSE}) we claim, in the
notation above:

\begin{equation}
\label{subomegaeq} \Sum_{d\in D} a_d
Sub_\omega[C^d_{g}(\Omega_1,\dots ,\Omega_p,\phi_1,\dots
,\phi_u)]+\Sum_{v\in V} a_v C^v_{g}(\Omega_1,\dots
,\Omega_p,\phi_1,\dots ,\phi_u,\omega)=0,
\end{equation}
where each $C^v_{g}$ is a contraction of length $\sigma +u+1$ in the form 
\begin{equation}
\label{xeilh'} 
contr(\nabla^{(m_1)}R\otimes\dots\otimes\nabla^{(m_a)}R\otimes\nabla^{(a)}\omega\otimes\nabla^{(r_1)}
\Omega_1\otimes\dots\otimes\nabla^{(r_p)}\Omega_p\otimes\nabla\phi_1\otimes\dots\otimes\nabla\phi_u),
\end{equation}
with
a factor $\nabla^{(a)}\omega, a\ge 2$.
\end{lemma}

{\it Proof:} The proof goes as follows: We re-write (\ref{LSE}) in
the form:

\begin{equation}
\label{LSE2} S_{g}=\Sum_{d\in D} a_d C^d_{g}(\Omega_1,\dots
,\Omega_p,\phi_1,\dots ,\phi_u)+\Sum_{h\in H} a_h
C^h_{g}(\Omega_1,\dots ,\Omega_p,\phi_1,\dots ,\phi_u)=0,
\end{equation}
where each $C^h_{g}$ has length $\ge\sigma +u+1$. We then re-write
this in a high dimension $N$ (we can do this since the equation
holds formally--see the discussion in the section on ``Trans-dimensional isomorphisms'' in \cite{a:dgciI}) 
and take $Image^1_\omega[S_{g}]$. We of course
have $Image^1_\omega[S_{g}]=0$. By virtue of the transformation
laws (\ref{curvtrans}), (\ref{levicivita}), we derive:

\begin{equation}
\label{paidi} \begin{split} & (0=)Image^1_\omega[S_{g^N}]= \Sum_{d\in
D} a_d N\cdot Sub_\omega[C^d_{g^N}(\Omega_1,\dots
,\Omega_p,\phi_1,\dots ,\phi_u)]+
\\&N\cdot \Sum_{v\in V} a_v
C^v_{g^N}(\Omega_1,\dots ,\Omega_p,\phi_1,\dots
,\phi_u,\omega)+
\\&N\Sum_{j\in J_1} a_j(N) C^j_{g^N}(\Omega_1,\dots
,\Omega_p,\phi_1,\dots ,\phi_u,\omega),
\end{split}
\end{equation}
where here the contractions $C^v_{g^N}$ are 
in the form (\ref{xeilh'}), have length $\sigma +u+1$
and a factor $\nabla^{(a)}\omega,a\ge 2$, while the contractions
$C^j_{g^N}$, $j\in J_1$ have length $\ge\sigma +u+2$; each coefficient 
$a_j(N)$ is apolynomial in $N$, of degree 0 or 1. Now,
re-writing the above in dimension $n$ and picking out the
sublinear combination of terms that are multiplied by
$N$ (notice this sublinear combination must vanish
separately) gives us our claim. $\Box$
\newline

{\bf $\nabla\upsilon$'s into $Xdiv$'s:} 
 We finally present a final technical Lemma which will be used on 
numerous occasions in this series of papers.

First some notation: We let $\sum_{f\in F} a_f
C^{f,i_1\dots i_\alpha}_g(\Omega_1,\dots,\Omega_p,\phi_1,\dots,\phi_u)$ stand for a
linear combination of $\alpha$-tensor fields, with each 
$C^{f,i_1\dots i_\alpha}_g(\Omega_1,\dots,\Omega_p,\phi_1,\dots,\phi_u)$
being  a partial contraction in the form:

\begin{equation}
\label{form2}
\begin{split}
&pcontr(\nabla^{(m_1)}R_{ijkl}\otimes\dots\otimes\nabla^{(m_s)}R_{ijkl}
\otimes
\\& \nabla^{(b_1)}\Omega_1\otimes\dots\otimes \nabla^{(b_p)}\Omega_p
\otimes\nabla\phi_1\otimes\dots\otimes\nabla\phi_u),
\end{split}
\end{equation}
 each having a given number $\sigma_1$ of factors 
$\nabla^{(m)}R_{ijkl}$, a given number $p$ of factors 
$\nabla^{(b)}\Omega_h$, $1\le h\le p$ and a given number $u$ of factors 
$\nabla\phi_y, 1\le y\le u$. We are also assuming that each $b_i\ge 2, 1\le i\le p$ 
and that each factors $\nabla\phi_h$ is 
contracting against one of the factors 
$\nabla^{(m)}R_{ijkl}, \nabla^{(b)}\Omega_h$.\footnote{In particular, 
no free index belongs to one of the factors $\nabla\phi_y$.} 
Furthermore, we assume that none of these tensor 
fields has an internal contraction.

\par We assume an equation:

\begin{equation}
\label{beginnow}
\begin{split}
&\sum_{f\in F} a_f C^{f,i_1\dots
i_\alpha}_g(\Omega_1,\dots,\Omega_p,\phi_1,\dots,\phi_u)
\nabla_{i_1}\upsilon\dots\nabla_{i_\alpha}\upsilon=
\\&\sum_{y\in Y} a_y C^{y,i_1\dots i_\mu}_g(\Omega_1,\dots,\Omega_p,\phi_1,\dots,\phi_u)
\nabla_{i_1}\upsilon\dots\nabla_{i_\alpha}\upsilon
\\&+ \sum_{z\in Z}
a_z C^{z,i_1\dots i_\alpha}_g(\Omega_1,\dots,\Omega_p,\phi_1,\dots,\phi_u)
\nabla_{i_1}\upsilon\dots\nabla_{i_\alpha}\upsilon,
\end{split}
\end{equation}
where the tensor fields in the RHS have length $\sigma+u+1$;
furthermore, the ones indexed in $Y$ have a factor
$\nabla\Omega_p$ (with only one derivative), while the ones
indexed in $Z$ have a factor $\nabla^{(b)}\Omega_p$ with $b\ge 2$.

\par We recall that for the tensor fields indexed in $F$, $Xdiv_{i_1}\dots Xdiv_{i_\alpha}
[C^{f,i_1\dots i_\alpha}_g]$ stands for the sublinear combination
in $div_{i_1}\dots div_{i_\alpha} [C^{f,i_1\dots i_\alpha}_g]$
where neither of the derivatives $\nabla^{i_r}$ is allowed to hit
any factor $\nabla\phi_t$, nor the factor $T$ to which ${}_{i_h}$
belongs. For the tensor fields indexed in $Y$, $Xdiv_{i_1}\dots
Xdiv_{i_\alpha} [C^{y,i_1\dots i_\alpha}_g]$ stands for the
sublinear combination in $div_{i_1}\dots div_{i_\alpha}
[C^{y,i_1\dots i_\alpha}_g]$ where neither of the derivatives
$\nabla^{i_r}$ is allowed to hit any factor $\nabla\phi_t$, nor
the factor $T$ to which ${}_{i_h}$ belongs, {\it nor} the factor
$\nabla\Omega_p$.

\par Our claim is the following:

\begin{lemma}
\label{Xdivext} Assume the equation (\ref{beginnow}). We then
claim that:

\begin{equation}
\label{beginnow'}
\begin{split}
&\sum_{f\in F} a_f Xdiv_{i_1}\dots Xdiv_{i_\alpha} C^{f,i_1\dots
i_\alpha}_g(\Omega_1,\dots,\Omega_p,\phi_1,\dots,\phi_u)=
\\&\sum_{y\in Y} a_y Xdiv_{i_1}\dots Xdiv_{i_\alpha}
C^{y,i_1\dots i_\alpha}_g(\Omega_1,\dots,\Omega_p,\phi_1,\dots,\phi_u)
\\&+ \sum_{z\in Z}
a_z C^{z}_g(\Omega_1,\dots,\Omega_p,\phi_1,\dots,\phi_u),
\end{split}
\end{equation}
here $\sum_{z\in Z}\dots$ stands for a generic linear combination of complete
 contractions in the form (\ref{form2}) 
 with length $\sigma+u+1$ and with a factor $\nabla^{(A)}\Omega_p, A\ge 2$
\end{lemma}

{\it Proof of Lemma \ref{Xdivext}:} We consider (\ref{beginnow})
and immediately derive an equation:

\begin{equation}
\label{beginnow2}
\begin{split}
&\sum_{f\in F} a_f div_{i_1}\dots div_{i_\alpha} C^{f,i_1\dots
i_\alpha}_g(\Omega_1,\dots,\Omega_p,\phi_1,\dots,\phi_u)=
\\&\sum_{y\in Y} a_y div_{i_1}\dots div_{i_\alpha}
C^{y,i_1\dots i_\alpha}_g(\Omega_1,\dots,\Omega_p,\phi_1,\dots,\phi_u)
\\&+ \sum_{z\in Z}
a_z div_{i_1}\dots div_{i_\alpha}C^{z,i_1\dots
i_\mu}_g(\Omega_1,\dots,\Omega_p,\phi_1,\dots,\phi_u).
\end{split}
\end{equation}

\par Now, we divide the LHS of the above into three linear combinations: $L^1$
is the sublinear combination which consists of terms with no
internal contractions and with one derivative on each function
$\phi_h$;\footnote{Observe that $L^1=\sum_{f\in F} a_f
Xdiv_{i_1}\dots Xdiv_{i_\alpha} C^{f,i_1\dots
i_\mu}_g(\Omega_1,\dots,\Omega_p)$.} $L^2$ is the sublinear
combination which consists of terms with at least one internal
contraction in some factor and and with one derivative on
each function $\phi_h$; $L^3$ stands for the sublinear combination
of terms with at least one function $\phi_h$ differentiated more than once
(and $\nabla^{(B)}\Omega_p$ still satisfies $B\ge 2$ by construction).

\par It easily follows that each of these three sublinear combinations
must vanish separately {\it at the linearized level}.\footnote{See the section 
``Background material'' in \cite{a:dgciI} for a 
strict definition of linearized complete contractions.} We denote by 
$lin\{L^1\}, lin\{L^2\}, lin\{L^3\}$ the linear combinations of 
linearized complete contractions that arise from $L^1,L^2,L^3$ 
by replacing each complete contraction 
$C_g(\Omega_1,\dots,\Omega_p,\phi_1,\dots,\phi_u)$ by its linearization 
$lin\{C_g(\Omega_1,\dots,\Omega_p,\phi_1,\dots,\phi_u)\}$.

\par Then, repeating the permutations by which we make
 the equation $lin\{L^2\}=0$ and $lin\{L^3\}=0$ formally zero to the non-linear setting, we
 derive that:

$$L^2=\sum_{z\in Z\bigcup Z'} a_z C^{z}_g(\Omega_1,\dots,\Omega_p,\phi_1,\dots,\phi_u),$$
where the terms indexed in $Z$ above are a generic linear
combinations with the properties described above. The terms
indexed in $Z'$ have length $\sigma+u+1$ and have only factors
$\nabla\phi_h, \nabla\Omega_p$ but also have at least one internal
contraction. By the same reasoning we derive an equation:

$$L^3=\sum_{z\in Z\bigcup Z'\bigcup Z''} a_z C^{z}_g(\Omega_1,\dots,\Omega_p,\phi_1,\dots,\phi_u)$$
where the tensor fields indexed in $Z''$ have length $\sigma+u+1$
and at least one factor $\nabla^{(B)}\phi_{h}$, $B\ge 2$.

\par Thus, replacing the above into (\ref{beginnow2}) we derive:

\begin{equation}
\label{beginnow3}
\begin{split}
&\sum_{f\in F} a_f Xdiv_{i_1}\dots Xdiv_{i_\alpha} C^{f,i_1\dots
i_\alpha}_g(\Omega_1,\dots,\Omega_p,\phi_1,\dots,\phi_u)=
\\&\sum_{y\in Y} a_y Xdiv_{i_1}\dots Xdiv_{i_\alpha}
C^{y,i_1\dots i_\alpha}_g(\Omega_1,\dots,\Omega_p,\phi_1,\dots,\phi_u)
\\&+ \sum_{z\in
Z\bigcup Z'\bigcup Z''} a_z C^{z}_g(\Omega_1,\dots,\Omega_p,\phi_1,\dots,\phi_u),
\end{split}
\end{equation}
(and the above holds modulo terms with length $\ge\sigma+u+2$).

\par Now, using the above we derive that we can write:

\begin{equation}
\label{beginnow4}
\begin{split}&\sum_{f\in F} a_f Xdiv_{i_1}\dots Xdiv_{i_\alpha}
C^{f,i_1\dots i_\mu}_g(\Omega_1,\dots,\Omega_p,\phi_1,\dots,\phi_u)
\\&=\sum_{m\in M} a_m
C^m_g(\Omega_1,\dots,\Omega_p,\phi_1,\dots,\phi_u),
\end{split}
\end{equation}
where the terms indexed in $M$
have length $\sigma+u+1$ and {\it no internal contractions}, and
also have {\it one derivative on each function $\phi_h$}.

\par Therefore, substituting the above into (\ref{beginnow3})
(and using the fact that (\ref{beginnow3}) holds modulo complete
contractions of length $\ge\sigma+u+2$), we derive that in
(\ref{beginnow3}):

$$\sum_{z\in Z'\bigcup Z''} a_z C^{z}_g(\Omega_1,\dots,\Omega_p,\phi_1,\dots,\phi_u)=0,$$
modulo complete contractions of length $\ge\sigma+u+2$. This
completes the proof of our claim. $\Box$

\end{document}